\documentclass[10pt,reqno,twoside]{amsart}

\usepackage{amsfonts}
\usepackage{paralist}
  \usepackage{graphics}
  \usepackage{epsfig}
 \usepackage[colorlinks=true]{hyperref}
\usepackage{hyperref}

\numberwithin{equation}{section}
\usepackage{latexsym}
\usepackage{amsmath}
\usepackage{amssymb}
\usepackage{mathrsfs}
\usepackage{graphicx}
\usepackage{ifthen}
\usepackage[T1]{fontenc}
\usepackage[latin1]{inputenc}
\usepackage{color}
\newtheorem{Satz}{Theorem}[section]

\newtheorem{Def}[Satz]{Definition}

\newcommand{\oo}{\overline \Omega}

\newcommand{\po}{\partial\Omega}

\newcommand{\dist}{\text {dist}}

\newcommand{\cat}{\text {cat}}

\begin{document}

\title[Sign-changing bubbling solutions]
{
Sign-changing bubbling solutions for
an exponential nonlinearity in $\mathbb{R}^2$
}



\author[Yibin Zhang]{Yibin Zhang}
\address{(Yibin Zhang)   College of Sciences, Nanjing Agricultural University,
Nanjing 210095, China}
\email{yibin10201029@njau.edu.cn}



\begin{abstract}
Very differently from those perturbative
techniques of  Deng-Musso in \cite{DM},
we use  the assumption of
a $C^1$-stable critical point to construct
positive or
 sign-changing   solutions with arbitrary $m$
isolated  bubbles
to the boundary value problem
$-\Delta u=\lambda u|u|^{p-2}e^{|u|^p}$  under
homogeneous Dirichlet boundary
condition in   a   bounded,   smooth planar domain  $\Omega$,
when $0<p<2$ and $\lambda>0$ is a small but free parameter.
We build a  vanishing identity of first order and
an identity of second order to  prove that for any $0<p<1$
the delicate energy expansion of these bubbling solutions
always converges to $4\pi m$ from below,
but for any $1<p<2$ the energy always converges to
$4\pi m$ from above, where the latter case sharply
recurs
a result of De Marchis-Malchiodi-Martinazzi-Thizy in \cite{DMMT}
involving    concentration and compactness
properties  at any critical energy level
$4\pi m$  only for  positive bubbling solutions.
A sufficient condition on the intersection between
 the nodal line of these sign-changing  solutions
and the boundary of the domain  is founded.
Moreover,
for $\lambda$ small enough,
we prove that
when $\Omega$ is an  arbitrary bounded domain,
this problem has not only at least two pairs of bubbling solutions which
change sign exactly once and whose nodal lines intersect the boundary,
but also a bubbling solution which changes sign exactly twice or three times;
when $\Omega$ has an axial symmetry,
this problem has a bubbling solution which
alternately changes sign arbitrarily many times
along the axis of symmetry through the domain.

\end{abstract}

\maketitle

\section{Introduction}
This paper deals with the existence and   asymptotic profile when the positive
but free parameter $\lambda$ tends to $0$ of positive or sign-changing
bubbling solutions in the distributional sense
for the following  problem
\begin{equation}\label{1.1}
\left\{\aligned
&-\Delta u=\lambda u|u|^{p-2}e^{|u|^p}
\,\,\,\,\,
\textrm{in}\,\,\,\,\,
\Omega,\\[1mm]
&u=0\,\,
\qquad\ \,\,
\  \ \qquad
\qquad\,
\textrm{on}\,\,\,
\partial\Omega,
\endaligned\right.
\end{equation}
where $0<p<2$ and $\Omega\subset\mathbb{R}^2$ is a $C^{2}$,  bounded and connected   domain.
This problem is related, but not equivalent,  to
a \textit{mean-field type equation}, namely the corresponding
Euler-Lagrange equation with the seemingly same form  as  in (\ref{1.1})
for the
  functional involving Moser-Trudinger critical and subcritical nonlinearities
 (see \cite{A,DMMT,F,T})
\begin{equation}\label{1.2}
\aligned
I_{p,\beta}(u)=\frac{2-p}2\left(
\frac{p\|u\|_{H_0^1(\Omega)}^2}{2\beta}
\right)^{\frac{p}{2-p}}
-\log
\int_{\Omega}\left(e^{|u|^p}-1
\right)dx,\,\,\ \ \,\,u\in H_0^1(\Omega),
\endaligned
\end{equation}
for any real number $\beta>0$,
whereas
$\lambda>0$ is not a free parameter but
a number $\lambda=\lambda_\beta$
defined by  the relation of the energy
\begin{equation}\label{1.2a}
\aligned
\frac{\lambda p}{2}
\left(
\int_{\Omega}
\big(e^{|u|^p}-1\big)dx
\right)^{\frac{2-p}{p}}
\left(
\int_{\Omega}
|u|^p
e^{|u|^p}dx
\right)^{\frac{2(p-1)}{p}}=\beta.
\endaligned
\end{equation}


For the precise statement,
it is useful  to recall some
well-known definitions.
Let $G(x,y)$  denotes
the
Green's function  of the   problem
\begin{equation}\label{1.10}
\left\{\aligned
&-\Delta_xG(x,y)=\delta_y(x),\,\,\,\,\,\,\,
x\in\Omega,\\
&G(x,y)=0,
\qquad\ \,\,
\qquad\ \,\,
x\in\partial\Omega,
\endaligned\right.
\end{equation}
and $H(x,y)$  its regular part  defined uniquely as
\begin{equation}\label{1.11}
\aligned
H(x,y)=
G(x,y)-\frac{1}{2\pi}\log\frac{1}{|x-y|}.
\endaligned
\end{equation}
In this way,
 for any $y\in\Omega$,
$H(\cdot,y)\in C^{1,\alpha}(\oo)\cap C^{\infty}(\Omega)$
and
$G(\cdot,y)\in C^{1,\alpha}_{loc}(\oo\setminus\{y\})\cap C^{\infty}_{loc}(\Omega\setminus\{y\})$
for any $0<\alpha<1$.
Moreover, the  Robin function $y\mapsto H(y,y)\in C^1(\Omega)$
(see \cite{GT}).
For any integer $m\geq1$, we introduce the classical
Kirchhoff-Routh path function $\varphi_m:\mathcal{F}_m(\Omega)\rightarrow \mathbb{R}$
of the form
\begin{equation}\label{1.3}
\aligned
\varphi_m(\xi)=\varphi_m(\xi_1,\ldots,\xi_m)
=
\sum\large_{i=1}^m
H(\xi_i,\xi_i)
+
\sum\large_{i,k=1,\,
i\neq k}^m
  a_i
  a_k
  G(\xi_i,\xi_k),
\endaligned
\end{equation}
where  $a_i\in\{-1,1\}$   and
$$
\aligned
\xi=(\xi_1,\ldots,\xi_m)\in
\mathcal{F}_m(\Omega):=\left\{\,\xi=(\xi_1,\ldots,\xi_m)\in\Omega^m
:\,\,\,\,
\xi_i\neq \xi_j
\quad
\textrm{if}
\,\,\,\,
i\neq j
\right\}.
\endaligned
$$
Additionally, for any $0<p<2$,
let $\gamma$ and $\varepsilon$ be two positive parameters  defined  uniquely  by the relations
\begin{equation}\label{1.5}
\aligned
p\lambda\gamma^{2(p-1)}\varepsilon^2e^{\gamma^p}=1,
\endaligned
\end{equation}
and
\begin{equation}\label{2.7}
\aligned
p\gamma^{p}=-4\log\varepsilon.
\endaligned
\end{equation}
Here, we take $\gamma$ and $\varepsilon$ such that
$\lambda\rightarrow0$ if and only if $\gamma\rightarrow+\infty$
and  $\varepsilon\rightarrow0$.
Moreover, $\lambda=\varepsilon^2$
and $\lambda^2e^{\gamma}=1$ if
$p=1$.

Clearly, for $p=1$,  by using the maximum principle we can write
 problem (\ref{1.1})
in terms of positive solutions   as
\begin{equation}\label{1.12}
\aligned
\left\{\aligned
&-\Delta
u=\varepsilon^2e^u\,\,\,\,\,
\,\,\,\,
\textrm{in}\,\,\,\,\,\,
\Omega,\\[2mm]
&u=0\qquad
\,\,\,\,\,
\quad\quad
\textrm{on}\,\,\,\,\po,
\endaligned\right.
\endaligned
\end{equation}
where $\varepsilon>0$ is a small parameter.
This is called
{\it Liouville equation}
after \cite{L}, which occurs in various context such
as the prescribed Gaussian curvature problem in conformal geometry
\cite{A1},
the mean-field limit of vortices in two dimensional  turbulent Euler flows
\cite{CLMP1,CLMP2},
and several other fields  of applied mathematics
\cite{BE,CI,CY}.
The asymptotic behavior of family of blow-up solutions of  equation (\ref{1.12})
has been  founded in \cite{BM,LS,MW,NS}.  Namely if $u_\varepsilon$ is a family of
positive solutions of equation (\ref{1.12}) satisfying
$$
\aligned
\lim_{\varepsilon\rightarrow0}
\big\|u_\varepsilon\big\|_{L^\infty(\Omega)}=+\infty
\,\,\,\,\,\quad\,\,\,\,
\textrm{and}
\,\,\,\,\quad\,\,\,\,
\lim_{\varepsilon\rightarrow0}\frac{\varepsilon^2}{2}\int_{\Omega}e^{u_\varepsilon}dx=\beta<+\infty,
\endaligned
$$
then, up to subsequences,
$\beta=4\pi m$, $m\in\mathbb{N}^*$ and
$u_\varepsilon$ makes
$m$ different blow-up
points  $\xi_1,\ldots,\xi_m\in\Omega$
 such that, as $\varepsilon\rightarrow0$,
 $$
\aligned
u_\varepsilon=8\pi\sum\limits_{i=1}^{m}G(x,\xi_i)+o(1)
\,\,\,\ \quad\ \ \,\,\,
\textrm{local uniformly in}
\,\,\,\,
\overline{\Omega}\setminus\{\xi_1,\ldots,\xi_m\},
\endaligned
$$
and
$$
\aligned
\varepsilon^2e^{u_\varepsilon}\rightharpoonup8\pi\sum_{i=1}^m \delta_{\xi_i}
\,\,\,\ \quad\ \ \,\,\,
\textrm{weakly in the sense of measure in}
\,\,\,\,
\overline{\Omega}.
\endaligned
$$
Moreover, the location of these bubbling points
$\xi=(\xi_1,\ldots,\xi_m)$
can be characterized as a critical point
of the  functional
\begin{equation}\label{1.13}
\aligned
\sum\large_{i=1}^m
H(\xi_i,\xi_i)
+
\sum\large_{i,k=1,\,
i\neq k}^m
  G(\xi_i,\xi_k).
\endaligned
\end{equation}
Reciprocally, the existence of positive solutions for equation
(\ref{1.12}) with exactly the asymptotic  profile above has been addressed in
\cite{BP,CL,DKM,EGP2,W}. In particular,
in the spirit of some perturbative and reductional techniques,
the construction of positive solutions with
arbitrary  $m$ distinct blow-up points
is achieved  respectively under  the three  different assumptions:
for any $m\geq1$ if
the functional defined in (\ref{1.13}) has a non-degenerate
critical point in $\mathcal{F}_m(\Omega)$ (\cite{BP}),
for any $m\geq1$ if
$\Omega$ is not simply connected (\cite{DKM}),
and for any  $m\in\{1,\ldots,h\}$ provided that $\Omega$ is an
$h$-dumbbell with thin handles (\cite{EGP2}).

For $p=2$,  problem (\ref{1.1}) becomes
\begin{equation}\label{1.14}
\aligned
\left\{\aligned
&-\Delta
u=\lambda ue^{u^2}
\,\,\,\,\,
\,\,\,\,\,
\textrm{in}\,\,\,\,
\Omega,\\[2mm]
&u=0
\,\,\,\,
\quad\quad\quad\qquad
\textrm{on}\,\,\,\po,
\endaligned\right.
\endaligned
\end{equation}
where $\lambda>0$ is a small but free parameter.
This problem is related, but not equivalent, to
the corresponding Euler-Lagrange equation
just as in (\ref{1.14})
for   the classical
Moser-Trudinger functional
\begin{equation}\label{1.12a}
\aligned
F(u)=
\int_{\Omega}
\big(e^{u^2}-1\big)dx,
\,\quad\ \,\,
\forall
\,
u\in H^1_0(\Omega),
\endaligned
\end{equation}
under the constraint of the Dirichlet energy
\begin{equation}\label{1.12b}
\aligned
\|u\|^2_{H^1_0(\Omega)}
=\beta,
\endaligned
\end{equation}
for any  $\beta>0$,  because for the latter equation,  $\lambda>0$ is not a free
parameter but a number $\lambda=\lambda_\beta$ defined by
\begin{equation}\label{}
\aligned
\lambda=\frac{2\beta}{
\,\langle
DF(u),\,u
\rangle\,
}
=\frac{\beta}{
\,\int_{\Omega}
u^2
e^{u^2}dx\,}.
\endaligned
\end{equation}
Indeed,  problem  (\ref{1.14})
is the Euler-Lagrange equation for  the  free energy functional
\begin{equation}\label{1.20}
\aligned
J_\lambda^{2}(u)=\frac12\int_{\Omega} |\nabla u|^2-\frac{\lambda}{2}
\int_{\Omega}e^{u^2},\,\,\ \ \,\,\forall\,
u\in H_0^1(\Omega).
\endaligned
\end{equation}
Due to  the critical Moser-Trudinger inequality, the semilinear Moser-Trudinger
functional $J_\lambda^{2}$
is well defined and smooth on $H_0^1(\Omega)$ and its positive critical points
are  smooth solutions of problem  (\ref{1.14}). The classification of
asymptotic  behavior of energy
bounded family of positive solutions of  (\ref{1.14}) has been described in \cite{A2,AD,AS,D,DT}.
Especially by the technique of the pointwise exhaustion of concentration points, it is proved
 in \cite{D,DT} that  if $u_{\varepsilon}$ is a family of positive solutions of problem (\ref{1.14}) with
$\lambda=\lambda_{\varepsilon}$  satisfying the free  or Dirichlet energy bound
$$
\aligned
\lim_{\varepsilon\rightarrow0}
J_{\lambda_\varepsilon}^{2}(u_{\varepsilon})=\frac{1}{2}\beta<+\infty
\,\,\,\,\,\ \quad\,\,\,\,
\textrm{or}
\,\,\,\,\ \quad\,\,\,\,
\lim_{\varepsilon\rightarrow0}\|u_\varepsilon\|^2_{H^1_0(\Omega)}=\beta<+\infty,
\endaligned
$$
then, up to subsequences, $\lambda_\varepsilon\rightarrow\lambda_0$ as $\varepsilon\rightarrow0$ for some $\lambda_0\geq0$,
and  if $\lambda_0>0$, then there exists a positive solution $u_0\in\mathcal{C}^2(\oo)$ of
\begin{equation*}
\aligned
\left\{\aligned
&-\Delta
u_{0}=\lambda_{0} u_{0}e^{u^2_{0}}
\,\ \,\,\,
\,\,\,\,\,
\textrm{in}\,\,\,\,
\Omega,\\[2mm]
&u_{0}=0
\,\,\,\,
\quad\quad\quad\quad\qquad
\textrm{on}\,\,\,\po,
\endaligned\right.
\endaligned
\end{equation*}
such that $u_\varepsilon\rightarrow u_0$ in $C^2(\oo)$ as $\varepsilon\rightarrow0$,
while if $\lambda_0=0$, then $u_\varepsilon\rightharpoonup0$ weakly in $H^1_{0}(\Omega)$ and
$\beta=4\pi m$, $m\in\mathbb{N}^*$, and
$u_\varepsilon$ makes
$m$ different blow-up
points  $\xi_1,\ldots,\xi_m\in\Omega$
and $m$ sequences $(\xi_{i,\varepsilon})_{\varepsilon}$ of points in $\Omega$
such that for each $i=1,\ldots,m$,
$\nabla u_{\varepsilon}(\xi_{i,\varepsilon})=0$,
$\gamma_{i,\varepsilon}:=u_{\varepsilon}(\xi_{i,\varepsilon})\rightarrow+\infty$,
$\xi_{i,\varepsilon}\rightarrow\xi_i$ and
$u_\varepsilon\rightarrow 0$ in $C_{loc}^2(\oo\setminus\{\xi_1,\ldots,\xi_m\})$
as $\varepsilon\rightarrow0$, and  there exists $\theta_i>0$ such that,
as $\varepsilon\rightarrow0$,
\begin{equation}\label{1.16}
\aligned
\sqrt{\lambda_\varepsilon}\gamma_{i,\varepsilon}\rightarrow\frac{2}{\theta_i},
\endaligned
\end{equation}
and
\begin{equation}\label{1.17}
\aligned
\gamma_{i,\varepsilon}\big[
u_{\varepsilon}(\xi_{i,\varepsilon}+\mu_{i,\varepsilon}y)-\gamma_{i,\varepsilon}
\big]\rightarrow
\log\frac{1}{1+|y|^2}
\qquad
\textrm{in}
\quad
C_{loc}^1(\mathbb{R}^2),
\endaligned
\end{equation}
  where
\begin{equation}\label{1.18}
\aligned
\mu_{i,\varepsilon}
=\left(\frac{4}{\,\lambda_{\varepsilon}\gamma_{i,\varepsilon}^{2}e^{\gamma_{i,\varepsilon}^2}\,}\right)^{\frac{1}{2}}
\rightarrow0.
\endaligned
\end{equation}
Moreover, the precise characterization  of these concentration  points
$\xi=(\xi_1,\ldots,\xi_m)$
and positive numbers
$\theta=(\theta_1,\ldots,\theta_m)$
can be seen  as a critical point
of the  functional
\begin{equation}\label{1.15}
\aligned
2\pi\sum\large_{i=1}^m
\theta_i^2H(\xi_i,\xi_i)
+
2\pi\sum\large_{i,k=1,\,
i\neq k}^m
\theta_i\theta_k
  G(\xi_i,\xi_k)
+\sum\large_{i=1}^m
\left(\theta_i^2-\theta_i^2\log\theta_i\right).
\endaligned
\end{equation}
Conversely, the existence of positive solutions of problem
(\ref{1.14})  has been  obtained by using  variational  (\cite{A2,DMR3})
or  perturbative  (\cite{DMR1}) method.
In particular,
positive solutions $u_\lambda$
of  problem
(\ref{1.14}) with
arbitrary  $m$ distinct blow-up points and
the energy quantization property
$\|u_\lambda\|^2_{H^1_0(\Omega)}\rightarrow 4m\pi$
or  $J_{\lambda}^{2}(u_{\lambda})\rightarrow 2m\pi$
as $\lambda\rightarrow0$
are  constructed  in \cite{DMR1}  under  the   assumption that
the functional defined in (\ref{1.15}) has a stable
critical point situation. Furthermore, if
$\Omega$ is not simple connected, at least one of
such  solutions  with $m=2$ always exists.
While if $\Omega$ has $d\geq1$ holes, additional $d+1$
bubbling solutions  with $m=1$ always exist.
Also,
similar results to \cite{D} and \cite{DMR1}
on compact Riemann surfaces
have been achieved  in \cite{Y} and \cite{FM}, respectively.

Recently for $p\in[1,2]$, motivated by these  quantization  results of the
possible blow-up energy levels in  \cite{D,DT,Y},
De Marchis-Malchiodi-Martinazzi-Thizy  have  obtained a precise  quantization
of  energy levels of  positive blow-up solutions
 on a closed Riemann surface,
and have proven  compactness at each critical energy level
$\beta=4\pi m$, $m\in\mathbb{N}^*$ for $p\in(1,2]$ by exhibiting
that delicate energy expansions
of blowing-up sequences of positive solutions always converge to $4\pi m$ from above
(see Sections $4$-$5$ in \cite{DMMT}). Furthermore, following closely  arguments  of
the pointwise exhaustion of concentration points
in
\cite{DMMT,D,DT}, it is stated in \cite{MMT} that if $u_{\varepsilon}$ is a sequence
of the pointwise blow-up of
positive solutions of equation (\ref{1.1}) with
$\lambda=\lambda_{\varepsilon}>0$ and $p=p_{\varepsilon}\in[1,2]$   satisfying
the energy bound
\begin{equation}\label{1.19}
\aligned
\lim_{\varepsilon\rightarrow0}
\beta_{\varepsilon}
:=
\lim_{\varepsilon\rightarrow0}
\frac{\lambda_{\varepsilon} p_{\varepsilon}}{2}
\left(
\int_{\Omega}
\big(e^{u_{\varepsilon}^{p_{\varepsilon}}}-1\big)dx
\right)^{\frac{2-p_{\varepsilon}}{p_{\varepsilon}}}
\left(
\int_{\Omega}
u_{\varepsilon}^{p_{\varepsilon}}
e^{u_{\varepsilon}^{p_{\varepsilon}}}dx
\right)^{\frac{2(p_{\varepsilon}-1)}{p_{\varepsilon}}}=\beta<+\infty,
\endaligned
\end{equation}
then, up to a subsequence,   $\lambda_\varepsilon\rightarrow0$ as $\varepsilon\rightarrow0$, and
$\beta=4\pi m$, $m\in\mathbb{N}^*$, and
$u_\varepsilon$ makes
$m$ different blow-up
points  $\xi_1,\ldots,\xi_m\in\Omega$
and $m$ sequences $(\xi_{i,\varepsilon})_{\varepsilon}$ of points in $\Omega$
such that for each $i=1,\ldots,m$,
$\nabla u_{\varepsilon}(\xi_{i,\varepsilon})=0$,
$\gamma_{i,\varepsilon}:=u_{\varepsilon}(\xi_{i,\varepsilon})\rightarrow+\infty$,
$\xi_{i,\varepsilon}\rightarrow\xi_i$ and
$u_\varepsilon\rightarrow u_0$ in $C_{loc}^1(\oo\setminus\{\xi_1,\ldots,\xi_m\})$
as $\varepsilon\rightarrow0$, and  as $\varepsilon\rightarrow0$,
\begin{equation}\label{1.21}
\aligned
\frac{1}{2}p_{\varepsilon}\gamma_{i,\varepsilon}^{p_{\varepsilon}-1}\big[
u_{\varepsilon}(\xi_{i,\varepsilon}+\mu_{i,\varepsilon}y)-\gamma_{i,\varepsilon}
\big]\rightarrow
\log\frac{1}{1+|y|^2}
\qquad
\textrm{in}
\quad
C_{loc}^1(\mathbb{R}^2),
\endaligned
\end{equation}
 where
\begin{equation}\label{1.22}
\aligned
\mu_{i,\varepsilon}
=\left(\frac{8}{\,\lambda_{\varepsilon}p_{\varepsilon}\gamma_{i,\varepsilon}^{2(p_{\varepsilon}-1)}
e^{\gamma_{i,\varepsilon}^{p_{\varepsilon}}}\,}\right)^{\frac{1}{2}}
\rightarrow0.
\endaligned
\end{equation}
Differently from the description (\ref{1.16}) of those bubble
heights $\gamma_{i,\varepsilon}$ under $p=2$, it is  easily known  from
Page $1225$ of \cite{DMMT}  that  if $p\in[1,2)$,
all the bubble heights have the same growth rate, namely
\begin{equation}\label{1.23}
\aligned
\lim_{\varepsilon\rightarrow0}
\frac{\gamma_{j,\varepsilon}}{\gamma_{i,\varepsilon}}=1
\qquad
\textrm{for all}
\,\,\,\,
i,j=1,\ldots,m,
\,\,\,
i\neq
j.
\endaligned
\end{equation}
Reciprocally for $p\in(0,2)$, the existence of the pointwise blow-up of
positive solutions for equation
(\ref{1.1}) with only some parts of
 asymptotic  properties  above has been proven by Deng-Musso in
\cite{DM}. Roughly speaking, by applying purely those techniques and
arguments on
the Lyapunov-Schmidt finite-dimensional reduction method from \cite{DKM} they
build positive solutions of equation (\ref{1.1}) with
arbitrary  $m$ distinct blow-up points under the assumption that
the functional defined in (\ref{1.13}) has a topologically non-trivial
critical level in an open set compactly contained in the configuration space $\mathcal{F}_m(\Omega)$.
Hence if $\Omega$ is not simply connected,
just like  \cite{DKM,EMP} they construct positive solutions of equation (\ref{1.1})
with arbitrary $m$ isolated bubbles.

One of our goals  is  to construct positive or sign-changing bubbling
solutions of equation (\ref{1.1}) not only for $p\in(0,2)$
but also partly for $p=2$, by trying to apply the \textit{standard}
 finite-dimensional  Lyapunov-Schmidt reduction method
 as in del Pino-Kowalczyk-Musso \cite{DKM} and Esposito-Musso-Pistoia \cite{EMP,EMP1}.
Based on the above asymptotic properties (\ref{1.21})-(\ref{1.23}) of positive solutions
of equation (\ref{1.1}) under $p\in[1,2)$, just like \cite{DM} we build the approximate form
of positive  or sign-changing solutions of equation (\ref{1.1})
which  not only have  the approximation
of up to third order, but also have the key  ingredients $\omega^j_{\mu_{i}}$,
$j=1,2,3$ (see (\ref{2.8}), (\ref{2.16}) and (\ref{2.17})).
Combining  relations (\ref{1.5})-(\ref{2.7}) with (\ref{1.22})-(\ref{1.23}), we set
each  bubble height $\gamma_{i,\varepsilon}=\gamma$ and
take the parameter  $\mu_{i,\varepsilon}=\varepsilon\mu_{i}$, where
the scaling parameter $\varepsilon$ is defined in (\ref{2.7}) and
the concentration parameter
$\mu_i$ is chosen by system (\ref{2.22})
 such that  we can obtain a good approximation for the solution
of equation (\ref{1.1}) near each blow-up point $\xi_i$.
Thanks to some ODE theory in  \cite{CI,EMP,MM}, we get the precise $C^1$-expansion (\ref{2.13})
of   $\omega^j_{\mu_{i}}$.
Hence for the scaling equation (\ref{2.26})  (or (\ref{2.29})), by applying
a  $L^\infty$-norm $\|\cdot\|_{*}$ involving
the new weighted function (\ref{2.32}) we estimate
its error term $E_{\xi'}$ as  $O\left(1/|\log\varepsilon|^4\right)$ and
obtain  the  delicate expansion (\ref{3.2}) for its potential $W_{\xi'}$
near each scaling point $\xi'_i$ which is absolutely  different from  the original  expansion
of $W_{\xi'}$ in Page $2274$ of \cite{DM}  because
expansion (\ref{3.2}) of $W_{\xi'}$ does not only depend on $e^{\omega_{\mu_i}}$, but  also
explicitly relates to  the term $\omega^1_{\mu_{i}}$. This hints that  when we perform a
finite-dimensional variational reduction in which the main ingredient is an analysis of
bounded invertibility of the linearized operator $\mathcal{L}$ in a suitable $L^\infty$-weighted space,
we should try our best to avoid adopting any perturbative and reductional
techniques from \cite{DM} again.  Indeed, if we carefully check the validity of the last equality
in Page $2286$ of \cite{DM},  we hardly yield the same result so that we never conclude bounded invertibility of
$\mathcal{L}$. Meanwhile if  we rigidly check the validity of the $C^1$-expansion of the finite-dimensional restriction
$F_{\lambda}(\xi)$ in Pages $2294$-$2299$ of \cite{DM}, we find that it is irrational to derive the claimed lemma
only from some very rough estimates at the price of losing  the key ingredient $\omega^1_{\mu_{i}}$.
Regrettably,  all these giant  and obvious gaps has also appeared in several other open literatures of
Deng-Musso  regarding the same nonlinearity $\lambda u^{p-1}e^{u^p}$ with
$p\in(0,2)$ and $\lambda>0$ free but small enough (see \cite{CDF,D1,D2, DM1,DM2,DM3,DM4,DGM}).
Whereas  observing that
expansion (\ref{3.2}) of potential $W_{\xi'}$ has  almost the same  structure as that of the corresponding
potential  in Page $43$ of \cite{EMP}, along those technical approach in \cite{EMP}
we successfully carry out an invertible analysis
of the linearized operator $\mathcal{L}$ and compute the precise $C^1$-expansion of the
reduction function $F_{\lambda}(\xi)$ involving the small $C^1$-perturbation of
the Kirchhoff-Routh path function $\varphi_{m}(\xi)$ under a suitable equivalent form.
With the help of the assumption of a $C^1$-stable critical point of $\varphi_{m}(\xi)$,
in this article we can recover
the main results for positive solutions in \cite{BP,DKM,DM,EGP2}
and construct sign-changing  solutions of equation (\ref{1.1})
with arbitrary $m$ isolated bubbles.

Finally,  it is necessary to point out that for equation (\ref{1.1})
with $p=2$, namely for equation (\ref{1.14})
by adopting  our  fundamental approximate form (\ref{2.8}) of positive or
sign-changing solutions and these technical arguments of the finite-dimensional variational reduction
in this article
we also have enough space to
recover the main result for positive solutions  in
del Pino-Musso-Ruf \cite{DMR1} and construct its  sign-changing  solutions
with arbitrary $m$ isolated bubbles for which the initial novelty
of the technique lies in the new choices of each
bubble height $\gamma_{i}$, scaling parameter $\varepsilon_{i}$ and concentration parameter  $\mu_{i}$
based on the difference of the growth rate of bubble heights between
 Moser-Trudinger nonlinear critical case (\ref{1.16}) and  subcritical case (\ref{1.23}). Roughly speaking, by (\ref{1.16}) we first
set each bubble height $\gamma_{i}=2/(\theta_i\sqrt{\lambda})$, where $\theta_i>0$
represents  the different bubble weight. Then similar to   (\ref{1.5})  we
take each scaling parameter $\varepsilon_i$ by  relation  $2\lambda\gamma_{i}^{2}\varepsilon_{i}^2e^{\gamma_{i}^2}=1$.
In order to obtain a good approximation for the solution
of equation (\ref{1.1}) under $p=2$ near each blow-up point $\xi_i$,
similar to system (\ref{2.22}) we adjust
each concentration parameter
$\mu_i$  through the following nonlinear system
\begin{eqnarray*}\label{}
\log\big(8\mu_i^2\big)=
2\log\frac{8}{\,\theta_i^2\,}
+
\left[
1-
\frac{1}{4}
\sum_{j=1}^3\frac{D^j_{\mu_i}}{(2\gamma_i^2)^{j}}
\right]
8\pi  H(\xi_i,\xi_i)
+\left[\sum_{j=1}^3\frac{D^j_{\mu_i}}{(2\gamma_i^2)^{j}}
\right]\log(\varepsilon_i\mu_i)
\qquad
\quad
\nonumber\\
+\sum_{k=1,\,k\neq i}^m
\frac{\theta_k}{\,\theta_i\,}
\left[
1-
\frac14\sum_{j=1}^3\frac{D^j_{\mu_k}}{(2\gamma_k^2)^{j}}
\right]
8\pi a_i a_k G(\xi_i,\xi_k),
\qquad\,\,\,
\textrm{for all}\,\,
i=1,\ldots,m,
\end{eqnarray*}
with $a_i\in\{-1,1\}$ and $D^j_{\mu_i}$ defined in (\ref{2.14})-(\ref{2.15}).
As a consequence, it  directly checks that for $p=2$,
the corresponding  potential of equation (\ref{1.1}) at
each approximate solution will have almost the same structure with  expansion (\ref{3.2})
of $W_{\xi'}$.
But for equation (\ref{1.1}) with $p=2$,
restrained by the too long pages of this present article we have to leave the
detailed discussion for the construction of its positive or sign-changing bubbling solutions   in the future manuscript.

In this paper, our  first result concerns the existence of positive or sign-changing
solutions of equation $(\ref{1.1})$ under $p\in(0,2)$ that blow up at
$m$ distinct points. Moreover, we compute the  delicate energy expansion
of second order  for
these bubbling solutions near the critical level $4\pi m$  and prove that
for any  $p\in(0,1)$, the  energy always converges
to $4\pi m$ from below, but for any $p\in(1,2)$, the  energy always converges
 to $4\pi m$ from above.  For the latter case  we sharply
recur the result of De Marchis-Malchiodi-Martinazzi-Thizy in Page $1237$ of \cite{DMMT}
regarding   concentration and compactness
properties  at any critical energy level
$4\pi m$  only for  positive bubbling solutions. In addition, we obtain
the pointwise $C^1$-limiting behavior of these solutions uniformly away from each bubbling point and
hence  find a sufficient condition on
the intersection between the nodal line of the sign-changing bubbling solutions and the boundary
of the domain.

\vspace {1mm}
\vspace {1mm}
\vspace {1mm}
\vspace {1mm}

\noindent{\bf Theorem 1.1.}\,\,\,{\it
Let $0<p<2$ and
$m$ be an integer with $m\geq1$.
Assume that
$\lambda>0$ is a  small but free parameter and
 $\xi^*=(\xi^*_1,\ldots,\xi^*_m)$ is a $C^1$-stable
 critical point for
the function $\varphi_m:\mathcal{F}_m(\Omega)\rightarrow \mathbb{R}$ in
the sense of Definition 6.1. Then, there exists
$\lambda_0>0$
such that for any $\lambda\in(0,\lambda_0)$,
problem
 {\upshape (\ref{1.1})} has a  pair of solutions $\pm u_\lambda$ such that, as $\lambda\rightarrow0$,
\begin{equation}\label{1.4}
\aligned
p\gamma^{p-1}
u_{\lambda}(x)\rightarrow
8\pi
\sum\limits_{i=1}^{m}
 a_i G(x,\xi_i^*)
\,\quad\textrm{in}\,\,\,\,
C^{1}_{loc}\big(\overline{\Omega}\setminus\{\xi^*_1,\ldots,\xi^*_m\}\big),
\endaligned
\end{equation}
and
\begin{equation}\label{1.6}
\aligned
p\gamma^{p-1}\lambda u_\lambda|u_\lambda|^{p-2}e^{|u_\lambda|^p}\rightharpoonup
8\pi\sum_{i=1}^ma_i\delta_{\xi_i^*}
\qquad
\textrm{weakly in the sense of measure in}
\,\,\,
\overline{\Omega},
\endaligned
\end{equation}
and
\begin{equation}\label{1.7}
\aligned
\beta_\lambda
=
\frac{\lambda p}{2}
\left(
\int_{\Omega}
\big(e^{|u_\lambda|^p}-1\big)dx
\right)^{\frac{2-p}{p}}
\left(
\int_{\Omega}
|u_\lambda|^p
e^{|u_\lambda|^p}dx
\right)^{\frac{2(p-1)}{p}}
=
4\pi m
+
O\left(\frac{1}{|\log\varepsilon|^2}\right),
\endaligned
\end{equation}
but for any $0<p<1$,
\begin{eqnarray}\label{1.7a}
\beta_\lambda\leq
4\pi m
\left\{1+
\frac{4(p-1)}{\,p^2\gamma^{2p}\,}
\left[1+
O\left(\frac{1}{|\log\varepsilon|}\right)
\right]\right\}
<4\pi m,
\end{eqnarray}
and for any $1<p<2$,
\begin{eqnarray}\label{1.7b}
\beta_\lambda\geq
4\pi m
\left\{1+
\frac{4(p-1)}{\,p^2\gamma^{2p}\,}
\left[1+
O\left(\frac{1}{|\log\varepsilon|}\right)
\right]\right\}
>4\pi m,
\end{eqnarray}
where  $a_i\in\{-1,1\}$,
 $\gamma$ and $\varepsilon$  are defined in  {\upshape(\ref{1.5})}-{\upshape(\ref{2.7})}.
More precisely,
\begin{eqnarray}\label{1.8}
u_\lambda(x)=\frac1{p\gamma^{p-1}}
\sum\limits_{i=1}^{m}a_i\left\{\,
\log
\frac1{\,((\varepsilon\mu^{\varepsilon}_i)^2+|x-\xi_i^\varepsilon|^2)^2\,}
+8\pi H(x,\xi_i^\varepsilon)
+\sum_{j=1}^3\left(\frac{p-1}{p}\right)^j\frac{1}{\gamma^{jp}}
\right.
\quad\qquad
\nonumber\\[1mm]
\left.
\times\left[\omega^j_{\mu^{\varepsilon}_i}\left(\frac{x-\xi_i^\varepsilon}{\varepsilon}\right)+\,D^j_{\mu^{\varepsilon}_i}\log(\varepsilon\mu^{\varepsilon}_i)
-2\pi D^j_{\mu^{\varepsilon}_i}H(x,\xi_i^\varepsilon)
\right]
+O\left(\frac{1}{|\log\varepsilon|^3}\right)\right\}
\quad\,
\textrm{as}\ \,
\lambda\rightarrow0
\end{eqnarray}
uniformly in $C(\oo)\cap H_0^1(\Omega)$,
where $\omega^j_{\mu^{\varepsilon}_i}$ and
$D^j_{\mu^{\varepsilon}_i}$ are defined in {\upshape(\ref{2.9})} and {\upshape(\ref{2.14})}, respectively,
$\mu_i^{\varepsilon}$ satisfies
$1/C\leq\mu_i^{\varepsilon}\leq C$
for some $C>0$,
and
$\xi^\varepsilon=(\xi^\varepsilon_1,\ldots,\xi^\varepsilon_{m})\in\mathcal{F}_m(\Omega)$
converges along a subsequence towards
$\xi^*$.
In addition, if
\begin{equation}\label{1.9a}
\aligned
a_1+\cdots+a_m=0,
\endaligned
\end{equation}
 then, for any $\lambda>0$ small enough,
\begin{equation}\label{1.9}
\aligned
\overline{\{x\in\Omega:\,u_{\lambda}(x)=0\}}\cap\po\neq\emptyset.
\endaligned
\end{equation}
}

\vspace {1mm}

When the free parameter $\lambda\rightarrow0$,  according to  Theorem 1.1 we can construct
a family of positive or sign-changing \textit{weak} solutions
$u_\lambda$ of
problem (\ref{1.1}) whose arbitrary $m$  bubbling points locate near a
$C^1$-stable critical point of  the path function $\varphi_{m}(\xi)$ and
whose energy $\beta_\lambda\rightarrow 4m\pi$
for $p=1$, $\beta_\lambda\uparrow 4m\pi$
for  $p\in(0,1)$, but $\beta_\lambda\downarrow 4m\pi$
for $p\in(1,2)$.  Observe that critical points $u_\beta$ of the aforementioned functional $I_{p,\beta}$ defined in (\ref{1.2})
exactly solve (\ref{1.1}), but $\lambda>0$ is a Lagrange multiplier
constrained  by the energy relation (\ref{1.2a}) of the parameter $\beta$.
So, Theorem 1.1 hints  that if the free parameter
$\beta\rightarrow 4m\pi$
for $p=1$,
$\beta \uparrow 4m\pi$
for $p\in(0,1)$, but $\beta \downarrow 4m\pi$
for $p\in(1,2)$, and
the path function $\varphi_{m}(\xi)$ has  a $C^1$-stable critical point
 $\xi^*=(\xi^*_1,\ldots,\xi^*_m)$ in $\mathcal{F}_m(\Omega)$, then,
the  functional $I_{p,\beta}$  will have a family of positive or sign-changing critical points $u_\beta$
in $H_0^1(\Omega)$
blowing up at each point $\xi_i^*$, $i=1,\ldots,m$ and $\lambda=\lambda_\beta\rightarrow0+$.
In this direction, it has been testified partly  for
the existence of positive bubbling critical points of $I_{p,\beta}$, which can be found in \cite{EGP2}
for $p=1$ and in \cite{DMMT,MMT} for  $p\in(1,2)$.
As for every $\beta>0$ and $p=2$,
the classical  Moser-Trudinger functional  (\ref{1.12a})
with the Dirichlet energy constraint (\ref{1.12b}) always has
a positive critical point, which is proven respectively in \cite{MMT} if $\Omega$ is non-simply connected
and  in \cite{DMMT} if the domain $\Omega$ is replaced by a closed Riemann surface.
Also,  it is proven  in  \cite{DMR2} that
the  Moser-Trudinger functional  (\ref{1.12a}) always has  a
positive 1-bubble critical point whose Dirichlet energy
$\beta\downarrow 4\pi$ if the smooth domain $\Omega$ is  bounded,
a
positive 2-bubble critical point whose Dirichlet energy
$\beta\downarrow 8\pi$ if $\Omega$ is non-simply connected,
a
positive $m$-bubble critical point whose Dirichlet energy
$\beta\downarrow 4m\pi$ if $\Omega$ is an annulus.
Additionally,  when $\Omega$ is a sufficiently symmetric domain,  it is proven
 in \cite{MTV} that for any integer $m\geq1$ and any $\beta\geq4m\pi$,
the  Moser-Trudinger functional  (\ref{1.12a}) has a family  of
sign-changing  critical points $u_\varepsilon$
with $m$  positive  spherical bubbles  clustered  at $0\in\Omega$
such that $\|u_\varepsilon\|^2_{H^1_0(\Omega)}
\rightarrow\beta$ and $u_\varepsilon\rightharpoonup u_0$ in $H^1_0(\Omega)$ as $\varepsilon\rightarrow0$,
where $u_0$ is  a sign-changing critical point of the  Moser-Trudinger functional  (\ref{1.12a})  under
the constraint of the Dirichlet energy at $\beta-4m\pi$ if $\beta>4m\pi$, but $u_0=0$ if $\beta=4m\pi$.

The proof of expansions
(\ref{1.7})-(\ref{1.7b}) of critical energy of  bubbling solutions
 is based on
two key identities (\ref{8.18}) and (\ref{9.10}) of first order and
second order.
When  $p\in(0,2)$, the  expansion  (\ref{1.7}) of second order
of energy $\beta_\lambda$
of  positive or sign-changing $m$-bubble
solutions arises from the delicate Taylor expansion and
the vanishing identity
(\ref{8.18}) of first order. For $m=1$, by using the identity
(\ref{9.10}) of second order we can compute
\begin{eqnarray}\label{6.14a}
\beta_\lambda
=4\pi
\left\{1+
\frac{4}{\,p^2\gamma^{2p}\,}
\left[p-1+
O\left(\frac{1}{|\log\varepsilon|}\right)
\right]\right\}.
\end{eqnarray}
For general $m\geq1$, by applying  H\"{o}lder's inequality for vectors in $\mathbb{R}_{+}^m$
and using identities (\ref{8.18}) and (\ref{9.10})
again  we conclude (\ref{1.7a})-(\ref{1.7b})
and hence observe that $\beta_\lambda\uparrow 4\pi m$
for  $p\in(0,1)$, but $\beta_\lambda\downarrow 4\pi m$
for $p\in(1,2)$.  Additionally, noticing that from Corollaries  A.8 and B.6
the identities (\ref{8.18}) and (\ref{9.10}) still hold for $p=2$,
along the analogous proof in (\ref{1.7}) and (\ref{1.7b}) it is possible
for the aforementioned construction of positive or sign-changing $m$-bubble
solutions  $u_\lambda$ of equation (\ref{1.1}) under $p=2$
to prove that the Dirichlet energy
\begin{eqnarray}
\beta_\lambda
=
\lambda
\int_{\Omega}
u^2_\lambda
e^{u^2_\lambda}dx
=\|u_\lambda\|^2_{H^1_0(\Omega)}
\geq
4\pi
\left\{m+
\sum_{i=1}^m\frac{1}{\,\gamma_i^{4}\,}
\big[1+
o\left(1\right)
\big]\right\}
\qquad
\textrm{and}
\qquad
\beta_\lambda\downarrow4\pi m,
\end{eqnarray}
and hence sharply
recur the corresponding inequality  in Page $1237$ of \cite{DMMT}, in particular for $m=1$ and $p=2$,
\begin{eqnarray}\label{6.14b}
\beta_\lambda
=
\lambda
\int_{\Omega}
u^2_\lambda
e^{u^2_\lambda}dx
=4\pi
+
\frac{4\pi}{\,\gamma_1^{4}\,}
+
o\left(\frac{1}{\,\gamma_1^{4}\,}\right),
\end{eqnarray}
where $\gamma_i$,  $i=1,\ldots,m$, are the aforementioned bubble heights defined by (\ref{1.16}).
Here, it easily checks that the second $4\pi$ appearing in (\ref{6.14b}) can be derived from
the identity  (\ref{9.10}) of second order under $p=2$ and hence obtain its solid interpretation
of mathematical structure
for $1$-bubble solution of equation (\ref{1.14}).
This, together with  asymptotic analysis  of $1$-bubble positive solutions of  equation (\ref{1.14})
in \cite{DMMT,MM}, implies that the
 identity  (\ref{9.10}) of second order under  $p=2$
gives an explicit geometric meaning to the term
$4\pi/\gamma_1^{4}$ in (\ref{6.14b}) and hence positively answers
the \textbf{Open Problem $3$} raised  by Mancini-Martinazzi
in Page $13$ of \cite{MM}.

If $\Omega$ is $C^2$,
the asymptotic behavior (\ref{1.4}) of solutions $u_{\lambda}$  uniformly away from each bubbling point
can be derived from  an \textit{a priori} $C^{1,\alpha}$ estimate in
Lemma 6.2 and the Ascoli-Arzel\'{a} Theorem. Furthermore, we prove that
if condition (\ref{1.9a}) holds, namely the total numbers of changing
signs  of bubbling solutions  $u_{\lambda}$ are even, then, the nodal line of these solutions
always intersects the boundary of $\Omega$. Moreover, we find that as for the construction of nodal
bubbling solutions with the intersection property between the nodal line of these solutions
and the boundary of $\Omega$,
condition  (\ref{1.9a}) is not only suitable for the Liouville type equation (\ref{1.1}) under $p\in(0,2)$, but also for
a two-dimensional sinh-Poisson equation in \cite{BP1} and a
two-dimensional  Lane-Emden equation with large exponent in \cite{EMP1}.
Additionally, motivated by \cite{BP1} with respect to the construction of sign-changing bubbling solutions on
a nonsmooth domain, it follows that
if $\Omega$ is a convex polygon with a finite number of corner points
$\{\varsigma_1,\ldots,\varsigma_n\}\subset\po$ and inner angles $\{\theta_1,\ldots,\theta_n\}$, $\theta_i\in(0,\pi)$,
then,  by using Grisvard's
elliptic regularity theory in nonsmooth domains \cite{G}
and Adams's standard Sobolev space theory  \cite{A} we believe  that the pointwise $C^1$-limiting behavior (\ref{1.4}) still holds
and hence Theorem 1.1 can be also proven to be true in this nonsmooth domain
although we apply the finite-dimensional  Lyapunov-Schmidt reduction method in the scheme of
a  $L^\infty$-space involving
the new weighted function (\ref{2.32}) instead of the classical  Hilbert space $H_0^1(\Omega)$ just as  \cite{BP1,EGP2}.

Remarks that the construction of positive or sign-changing bubbling solutions of  problem (\ref{1.1}) strongly depends
on the  assumption   that  the Kirchhoff-Routh  path function $\varphi_{m}(\xi)$ has  a $C^1$-stable critical point
 $\xi^*=(\xi^*_1,\ldots,\xi^*_m)$ in $\mathcal{F}_m(\Omega)$.
If $a_i=1$ for all $i=1,\ldots,m$, from \cite{DKM,EGP2} it follows that whether
$\Omega$ is non-simply connected
or  an
$h$-dumbbell with thin handles and $h\geq m$, such a stable critical point
always exists and hence Theorem 1.1 implies that   problem  (\ref{1.1})
admits  a positive  solution  blowing up at each point $\xi_i^{*}$, $i=1,\ldots,m$.
While if $a_ia_j=-1$ for some $i$ and $j$, we observe that such a stable critical point also always
exists in the following three cases:\\
(C1) \,
if $m=2$, $a_1=1$ and $a_2=-1$ (see \cite{BP1,BP2,EMP1});\\
(C2) \,
if $m=3$ or $m=4$ and $a_i=(-1)^{i+1}$, $i=1,\ldots,m$ (see \cite{BP2});\\
(C3) \,
if  $\Omega$ is symmetric with respect a line through it, then for  all $m\geq1$  and $a_i=(-1)^{i+1}$,
$i=1,\ldots,m$ (see  \cite{BPW});\\
hence by Theorem 1.1 we get the existence of a sign-changing solution of   (\ref{1.1}) blowing up at each point $\xi_i^{*}$.
More exactly, let
$$
\aligned
\mathcal{C}_2(\Omega):=
\mathcal{F}_2(\Omega)/(\xi_1,\xi_2)\sim(\xi_2,\xi_1)=
\left\{(\xi_1,\xi_2)\in\Omega\times\Omega
:\,\,
\xi_1\neq \xi_2
\right\}/(\xi_1,\xi_2)\sim(\xi_2,\xi_1)
\endaligned
$$
be the quotient manifold of $\mathcal{F}_2(\Omega)$
modulo  the equivalence $(\xi_1,\xi_2)\sim(\xi_2,\xi_1)$
 and
define  $\cat\big(\mathcal{C}_2(\Omega)\big)$
as the Lusternik-Schnirelmann category of $\mathcal{C}_2(\Omega)$.
Here $\cat\big(\mathcal{C}_2(\Omega)\big)\geq2$
(see \cite{BMP}). Then we have the following theorems.

\vspace {1mm}
\vspace {1mm}
\vspace {1mm}
\vspace {1mm}

\noindent{\bf Theorem 1.2.}\,\,\,{\it
Fix $m=2$. Then for any $p\in(0,2)$,
there exists
$\lambda_0>0$ such that for any $\lambda\in(0,\lambda_0)$, \\
\indent {\upshape (i)}\,\,
  problem
 {\upshape (\ref{1.1})} has at least
 $k:=\cat\big(\mathcal{C}_2(\Omega)\big)$
pairs of sign-changing solutions $\pm u_\lambda^i$ with $i=1,\ldots,k$
such that  as $\lambda\rightarrow0$,
$$
\aligned
p\gamma^{p-1}
\lambda u_\lambda^i|u_\lambda^i|^{p-2}e^{|u_\lambda^i|^p}\rightharpoonup
8\pi(\delta_{\xi_1^i}
-
\delta_{\xi_2^i}
)
\qquad
\,\,
\textrm{weakly in the sense of measure in}
\,\,\,\,
\overline{\Omega},
\endaligned
$$
\indent\,\,\indent where  the blow-up point
$\xi^i=(\xi_{1}^i,\,\xi_{2}^i)\in\mathcal{F}_2(\Omega)$
is a critical point of $\varphi_2$ with $a_1=1$ and $a_2=-1$;\\
\indent {\upshape (ii)}\,\,the set $\Omega\setminus\{x\in\Omega:\,u_{\lambda}^i(x)=0\}$
has exactly two connected components;\\
\indent {\upshape (iii)}\,$\overline{\{x\in\Omega:\,u_{\lambda}^i(x)=0\}}\cap\po\neq\emptyset$.
}

\vspace {1mm}
\vspace {1mm}
\vspace {1mm}
\vspace {1mm}

\noindent{\bf Theorem 1.3.}\,\,\,{\it
Fix
 $m=3$ or $m=4$.
Then for any $p\in(0,2)$,
there exists
$\lambda_0>0$ such that for any $\lambda\in(0,\lambda_0)$,
  problem
 {\upshape (\ref{1.1})} has a pair of sign-changing solutions $\pm u_\lambda$
such  that
as $\lambda\rightarrow0$,
$$
\aligned
p\gamma^{p-1}
\lambda u_\lambda|u_\lambda|^{p-2}e^{|u_\lambda|^p}\rightharpoonup
8\pi\sum_{i=1}^m
(-1)^{i+1} \delta_{\xi_i^*}
\qquad
\,\,
\textrm{weakly in the sense of measure in}
\,\,\,\,
\overline{\Omega},
\endaligned
$$
where the blow-up point
$\xi^*=(\xi_{1}^*,\ldots,\xi_{m}^*)\in\mathcal{F}_m(\Omega)$ is
a critical point of $\varphi_m$ with $a_i=(-1)^{i+1}$,
$i=1,\ldots,m$.
}

\vspace{1mm}
\vspace{1mm}
\vspace{1mm}
\vspace{1mm}

\noindent{\bf Theorem 1.4.}\,\,\,{\it
Assume that
$\Omega\cap(\mathbb{R}\times\{0\})\neq \emptyset$
and $\Omega$ is
symmetric with respect to the reflection at
$\mathbb{R}\times\{0\}$.
Then  for any $p\in(0,2)$ and any fixed integer $m\geq1$,
 there exists
$\lambda_0=\lambda_0(m)>0$ such that for any $\lambda\in(0,\lambda_0)$, problem
 {\upshape (\ref{1.1})} has a pair of sign-changing solutions $\pm u_\lambda$ such that
as $\lambda\rightarrow0$,
$$
\aligned
p\gamma^{p-1}\lambda u_\lambda|u_\lambda|^{p-2}e^{|u_\lambda|^p}\rightharpoonup
8\pi\sum_{i=1}^m
(-1)^{i+1} \delta_{\xi_i^*}
\qquad
\,\,
\textrm{weakly in the sense of measure in}
\,\,\,\,
\overline{\Omega},
\endaligned
$$
where the blow-up point
$\xi^*=(\xi_{1}^*,\ldots,\xi_{m}^*)\in\mathcal{F}_m(\Omega)$ is a critical point of $\varphi_m$ with $a_i=(-1)^{i+1}$,
$i=1,\ldots,m$, and it satisfies
$\xi_{i}^*=(t_i,\,0)$,
$t_1<t_2<\cdots<t_m$.}

\vspace{1mm}
\vspace{1mm}
\vspace{1mm}
\vspace{1mm}

%
%
%

This paper is organized as follows:
In Section $2$ we  describe
a good approximation for the solution of problem (\ref{1.1}) and estimate the error.
Then we rewrite  problem (\ref{1.1}) in terms of a linearized operator $\mathcal{L}$
for which
a solvability theory
is performed  through solving a linearized  problem in Section $3$. In Section $4$ we
solve  a  nonlinear projected problem.
Then  we reduce (\ref{1.1}) to
solve a finite system $c_{ij}=0$ and compute the $C^1$-expansion of
the  reduction function $F_\lambda$, as we will see in Section $5$.
Section $6$ gives the proof of Theorems 1.1-1.4.
In the final two sections we apply a precise expression of
the key ingredient $\omega^1_{\mu_{i}}$ to
give the proof of the vanishing identity
(\ref{8.18}) of first order and  the  identity  (\ref{9.10}) of second order, respectively.

Notation:
In this paper the letters $C$ and $D$ will always denote
some  universal positive  constants  independent of $\lambda$ and $\varepsilon$,
which could be changed from one line to another.
The symbol $o(t)$ (respectively $O(t)$) will denote a quantity for which
$\frac{o(t)}{|t|}$ tends to zero
(respectively, $\frac{O(t)}{|t|}$ stays bounded )
as parameter $t$ goes to zero.
Moreover, we will use the notation
$o(1)$ (respectively $O(1)$)
to stand for a quantity which tends to zero
(respectively, which remains uniformly bounded) as $\lambda$  and $\varepsilon$
tend to zero.

\vspace{1mm}
\vspace{1mm}
\vspace{1mm}

\section{An approximation for the solution}
The basic cells to obtain  an  approximate solution of problem (\ref{1.1})
are given by  the four-parameter family of functions
\begin{equation}\label{2.1}
\aligned
\omega_{\varepsilon,\mu,\xi}(z)=\log\frac{8\mu^2}{(\varepsilon^2\mu^2+|z-\xi|^2)^2},\quad
\,\,\varepsilon>0,\,\,\,\,\mu>0,\,\,\,\,\xi\in\mathbb{R}^2,
\endaligned
\end{equation}
which exactly solve
$$
\aligned
-\Delta
\omega=\varepsilon^2e^\omega\,\ \,\,\,\,
\textrm{in}\,\,\,\,
\mathbb{R}^2,
\qquad\quad\qquad
\int_{\mathbb{R}^2}\varepsilon^2e^\omega=8\pi.
\endaligned
$$
Set
\begin{equation}\label{2.3}
\aligned
\omega_{\mu}(z)=\omega_{1,\mu,(0,0)}(|z|)\equiv\log\frac{8\mu^2}{(\mu^2+|z|^2)^2}.
\endaligned
\end{equation}
The configuration space for $m$ concentration points $\xi=(\xi_1,\ldots,\xi_m)$ we
try to look for  is the following
\begin{eqnarray}\label{2.4}
\mathcal{O}_d:=\left\{\,\xi=(\xi_1,\ldots,\xi_m)\in\Omega^m
:\,
|\xi_i-\xi_j|\geq 4d,
\,\,\,\,\,\,\,
\dist(\xi_i,\po)\geq 4d,
\quad
i,j=1,\ldots,m,\,\,\,i\neq j
\right\},
\end{eqnarray}
where
$d>0$ is a small but fixed number.
Let  us fix
$\xi=(\xi_1,\ldots,\xi_m)\in\mathcal{O}_d$.
For numbers $\mu_i$,
$i=1,\ldots,m$, yet to be chosen, but  we always consider
\begin{equation}\label{2.6}
\aligned
d\leq\mu_i\leq1/ d,
\quad\,\,\,
\,\,\,
i=1,\ldots,m.
\endaligned
\end{equation}
Let
\begin{equation}\label{2.8}
\aligned
U_i(x)=\frac{1}{p\gamma^{p-1}}\left[
\omega_{\varepsilon,\mu_i,\xi_i}(x)
+\sum_{j=1}^3\left(\frac{p-1}{p}\right)^j\frac{1}{\gamma^{jp}}\omega^j_{\mu_{i}}\left(\frac{x-\xi_i}{\varepsilon}\right)
\right],
\quad\
i=1,\ldots,m.
\endaligned
\end{equation}
Here,  $\omega^j_{\mu_i}$,  $j=1,2,3$,
 are radial solutions of
\begin{equation}\label{2.9}
\aligned
\Delta\omega^j_{\mu_i}+e^{\omega_{\mu_i}(|z|)}\omega^j_{\mu_i}=e^{\omega_{\mu_i}(|z|)}f^j_{\mu_i}
\qquad
\textrm{in}\,\,\ \,\,\mathbb{R}^2,
\endaligned
\end{equation}
 with
\begin{equation}\label{2.10}
\aligned
f^1_{\mu_i}=-\left[
\omega_{\mu_i}+\frac{1}{2}(\omega_{\mu_i})^2
\right],
\endaligned
\end{equation}
and
\begin{equation}\label{2.11}
\aligned
f^2_{\mu_i}=-\left\{
\left[
\omega^{1}_{\mu_i}+\frac{p-2}{2(p-1)}(\omega_{\mu_i})^2\right]
+\omega_{\mu_i}\left[\omega^{1}_{\mu_i}+\frac{1}{2}(\omega_{\mu_i})^2\right]
+\omega_{\mu_i}\omega^{1}_{\mu_i}+\frac{p-2}{6(p-1)}(\omega_{\mu_i})^3
+\frac12\left[\omega^{1}_{\mu_i}+\frac{1}{2}(\omega_{\mu_i})^2\right]^2
\right\},
\endaligned
\end{equation}
and
\begin{eqnarray}\label{2.12}
f^3_{\mu_i}=-\left\{
\left[
\omega^{2}_{\mu_i}+\frac{p-2}{p-1}\omega_{\mu_i}\omega^1_{\mu_i}
+\frac{(p-2)(p-3)}{6(p-1)^2}(\omega_{\mu_i})^3
\right]
+
\left[
\omega^{1}_{\mu_i}+\frac{p-2}{2(p-1)}(\omega_{\mu_i})^2\right]
\left[\omega^{1}_{\mu_i}+\frac{(\omega_{\mu_i})^2}{2}\right]
+\omega_{\mu_i}\Big[
\omega_{\mu_i}\omega^{1}_{\mu_i}
\right.
\nonumber\\
\left.\left.
+\,\,\omega^{2}_{\mu_i}
+
\frac{p-2}{6(p-1)}(\omega_{\mu_i})^3+
\frac12\left(\omega^{1}_{\mu_i}+\frac{1}{2}(\omega_{\mu_i})^2\right)^2\right]
+
\omega_{\mu_i}\omega^{2}_{\mu_i}
+\frac{p-2}{2(p-1)}(\omega_{\mu_i})^2\omega^{1}_{\mu_i}
+\frac{(p-2)(p-3)}{24(p-1)^2}(\omega_{\mu_i})^4
\ \,\
\right.
\nonumber\\
\left.
+\,
\frac{1}2(\omega^{1}_{\mu_i})^2
+
\left[\omega^{1}_{\mu_i}+\frac{1}{2}(\omega_{\mu_i})^2\right]\left[
\omega^{2}_{\mu_i}+\omega_{\mu_i}\omega^{1}_{\mu_i}+\frac{p-2}{6(p-1)}(\omega_{\mu_i})^3\right]
+\frac16\left[\omega^{1}_{\mu_i}+\frac{1}{2}(\omega_{\mu_i})^2\right]^3
\right\},
\quad\quad\,\,\,
\textrm{for}\,\,\,
p\neq1,\
\end{eqnarray}
having asymptotic (see \cite{CI,EMP,MM})
\begin{equation}\label{2.13}
\aligned
\left\{
\aligned
&\omega^j_{\mu_i}(z)=\frac{D^j_{\mu_i}}{2}\log\left(1+\frac{|z|^2}{\mu_i^2}\right)+O\left(\frac{\mu_i}{\mu_i+|z|}\right)
\,\,\ \ \quad\,\,
\textrm{as}\,\,\,
|z|\rightarrow+\infty,\\
&\nabla\omega^j_{\mu_i}(z)=D^j_{\mu_i}
\cdot
\frac{ z}{\mu_i^2+|z|^2}+O\left(\frac{\mu_i}{\mu_i^2+|z|^2}\right)
\,\,
\qquad
\,\,\,\textrm{for all}\,\,\,z\in\mathbb{R}^2,
\endaligned
\right.
\endaligned
\end{equation}
for $j=1, 2, 3$,   where
\begin{equation}\label{2.14}
\aligned
D^j_{\mu_i}=\frac{1}{2\pi}\int_{\mathbb{R}^2}
\Delta\big[\omega^j_{\mu_i}(\mu_iy)\big]dy
\qquad\quad
\textrm{and}
\qquad\quad
D^j_{\mu_i}=8\int_{0}^{+\infty}t\frac{t^2-1}{(t^2+1)^3}f^j_{\mu_i}(\mu_i t)dt,
\endaligned
\end{equation}
in particular, $D^2_{\mu_i}$ is   explicitly computed by (\ref{9.8}) and
\begin{equation}\label{2.15}
\aligned
D^1_{\mu_i}=4\log8-8-8\log\mu_i.
\endaligned
\end{equation}

We now  approximate the solution of problem (\ref{1.1}) by
\begin{equation}\label{2.16}
\aligned
U_\xi(x):=\sum_{i=1}^m
a_i
PU_i(x)=\sum_{i=1}^m
a_i
\big[U_i(x)+H_i(x)\big],
\endaligned
\end{equation}
where  $a_i\in\{-1,1\}$  and   $H_i$  is a correction term defined as the solution of
\begin{equation}\label{2.17}
\aligned
-\Delta
H_i=0
\,\,\,\,\,\,
\textrm {in}\,\,\,\,\,\,\Omega,
\qquad\qquad
H_i
=-U_i
\,\,\,\,\,\,
\textrm{on}\,\,\,\,\po.
\endaligned
\end{equation}

\vspace{1mm}
\vspace{1mm}

\noindent{\bf Lemma 2.1.}\,\,{\it
For any $i=1,\ldots,m$ and for any $\varepsilon$ small enough,
\begin{eqnarray}\label{2.19}
H_i(x)=\frac{1}{p\gamma^{p-1}}\left\{
\left[
1-
\frac{1}{4}
\sum_{j=1}^3\left(\frac{p-1}{p}\right)^j\frac{D^j_{\mu_i}}{\gamma^{jp}}
\right]
8\pi H(x,\xi_i)-\log(8\mu_i^2)
+\left[\sum_{j=1}^3\left(\frac{p-1}{p}\right)^j\frac{D^j_{\mu_i}}{\gamma^{jp}}
\right]
\log(\varepsilon\mu_i)
+O\left(\frac{\varepsilon}{|\log\varepsilon|}\right)
\right\},
\end{eqnarray}
where the convergence    holds in $C^{1}(\oo)\cap C^{\infty}(\Omega)$
uniformly for any $\xi=(\xi_1,\ldots,\xi_m)\in\mathcal{O}_d$
and for any  $\mu=(\mu_1,\ldots,\mu_m)$ satisfying assumption
  {\upshape(\ref{2.6})}.
}

\vspace{1mm}

\begin{proof}
Since $\dist(\xi_i,\po)\geq 4d$ for any
$i=1,\ldots,m$, by
  (\ref{2.1}), (\ref{2.8}) and (\ref{2.13})
 we readily get
 $$
\aligned
H_i(x)=&\,
\frac{1}{p\gamma^{p-1}}\left\{
\left[2
-\frac{1}{2}\sum_{j=1}^3
\left(\frac{p-1}{p}\right)^j\frac{D^j_{\mu_i}}{\gamma^{jp}}
\right]
\log\big(\varepsilon^2\mu_i^2+|x-\xi_i|^2\big)
-\log(8\mu_i^2)
+\left[
\sum_{j=1}^3
\left(\frac{p-1}{p}\right)^j\frac{D^j_{\mu_i}}{\gamma^{jp}}
\right]\log(\varepsilon\mu_i)
\right.\\[1mm]
&
\left.
+
\frac{p-1}{p}\frac{1}{\gamma^p}
O
\left(\frac{\varepsilon\mu_i}{\varepsilon\mu_i+|x-\xi_i|}
\right)
\right\}
\endaligned
$$
 in $C^1(\po)$ as $\varepsilon\rightarrow0$.
Consider the harmonic function
$$
\aligned
Z_i(x)=p\gamma^{p-1}H_i(x)
-\left[
1-
\frac{1}{4}
\sum_{j=1}^3\left(\frac{p-1}{p}\right)^j\frac{D^j_{\mu_i}}{\gamma^{jp}}
\right]8\pi
H(x,\xi_i)
+\log(8\mu_i^2)
-\left[\sum_{j=1}^3\left(\frac{p-1}{p}\right)^j\frac{D^j_{\mu_i}}{\gamma^{jp}}
\right]\log(\varepsilon\mu_i).
\endaligned
$$
From (\ref{1.10})-(\ref{1.11})  we  have clearly that for any $\xi=(\xi_1,\ldots,\xi_m)\in\mathcal{O}_d$
and for any  $\mu=(\mu_1,\ldots,\mu_m)$ satisfying  assumption
  {\upshape(\ref{2.6})},
\begin{equation}\label{2.18}
\aligned
Z_i(x)
=O\left(\frac{\varepsilon}{|\log\varepsilon|}\right)
\qquad\
\textrm{uniformly}
\ \,
\textrm{in}
\ \
C^1(\po).
\endaligned
\end{equation}
According to the  maximum principle
and the Green's representation formula for harmonic function
we derive the $C^\infty(\Omega)$-convergence in expansion (\ref{2.19}).
By  Theorem 6.13 in \cite{GT}
we deduce
the $C(\oo)\cap C^{\infty}(\Omega)$-convergence in
 (\ref{2.19}).
Furthermore, using the boundary gradient estimate in Proposition 2.20 of \cite{HL},
by (\ref{2.18}) we  obtain
the $C^{1}(\oo)\cap C^{\infty}(\Omega)$-convergence in
(\ref{2.19}).
\end{proof}

\vspace{1mm}

From Lemma 2.1, uniformly far away from each point  $\xi_i$,
namely  $|x-\xi_i|\geq d$ for all $i=1,\ldots,m$,
one has
\begin{equation}\label{2.20}
\aligned
U_{\xi}(x)=\frac{1}{p\gamma^{p-1}}
\sum\limits_{i=1}^{m}
\left\{
\left[
1
-\frac14
\sum_{j=1}^3\left(\frac{p-1}{p}\right)^j\frac{D^j_{\mu_i}}{\gamma^{jp}}
\right]
8\pi a_i G(x,\xi_i)
+O\left(
\varepsilon
\right)
\right\}.
\endaligned
\end{equation}
But for
$|x-\xi_i|< d$  with some $i\in\{1,\ldots,m\}$,
$$
\aligned
PU_i(x)=&\,\frac{1}{p\gamma^{p-1}}
\left\{
p\gamma^{p}+
\omega_{\mu_i}\left(\frac{x-\xi_i}{\varepsilon}\right)
+\sum_{j=1}^3\left(\frac{p-1}{p}\right)^j\frac{1}{\gamma^{jp}}\omega^j_{\mu_{i}}\left(\frac{x-\xi_i}{\varepsilon}\right)
+\left[
1
-\frac14
\sum_{j=1}^3\left(\frac{p-1}{p}\right)^j\frac{D^j_{\mu_i}}{\gamma^{jp}}
\right]
8\pi  H(\xi_i,\xi_i)
\right.\\[1mm]
&\left.
-\log(8\mu_i^2)
+\left[
\sum_{j=1}^3\left(\frac{p-1}{p}\right)^j\frac{D^j_{\mu_i}}{\gamma^{jp}}
\right]\log(\varepsilon\mu_i)
+O\left(
|x-\xi_i|
+\varepsilon\right)
\right\},
\endaligned
$$
and for any $k\neq i$,
$$
\aligned
PU_k(x)=\frac{1}{p\gamma^{p-1}}
\left\{
\left[
1
-\frac14
\sum_{j=1}^3\left(\frac{p-1}{p}\right)^j\frac{D^j_{\mu_k}}{\gamma^{jp}}
\right]
8\pi  G(\xi_i,\xi_k)
+O\left(|x-\xi_i|
+
\varepsilon
\right)
\right\}.
\endaligned
$$
Hence  for
$|x-\xi_i|< d$,
\begin{eqnarray}\label{2.21}
U_{\xi}(x)=\frac{a_i}{p\gamma^{p-1}}
\left[
p\gamma^{p}+
\omega_{\mu_i}\left(\frac{x-\xi_i}{\varepsilon}\right)
+\sum_{j=1}^3\left(\frac{p-1}{p}\right)^j\frac{1}{\gamma^{jp}}\omega^j_{\mu_{i}}\left(\frac{x-\xi_i}{\varepsilon}\right)
+\,O\left(
|x-\xi_i|
+\varepsilon
\right)
\right]
\end{eqnarray}
will be a good  approximation for the  solution of problem (\ref{1.1})
near the point $\xi_i$
provided that  for each $i=1,\ldots,m$, the
concentration
parameter $\mu_i$ satisfies the nonlinear system
\begin{eqnarray}\label{2.22}
\log\big(8\mu_i^2\big)=\left[
1-
\frac{1}{4}
\sum_{j=1}^3\left(\frac{p-1}{p}\right)^j\frac{D^j_{\mu_i}}{\gamma^{jp}}
\right]
8\pi  H(\xi_i,\xi_i)
+\left[\sum_{j=1}^3\left(\frac{p-1}{p}\right)^j\frac{D^j_{\mu_i}}{\gamma^{jp}}
\right]\log(\varepsilon\mu_i)
\nonumber\\
+\sum_{k=1,\,k\neq i}^m
\left[
1-
\frac14\sum_{j=1}^3\left(\frac{p-1}{p}\right)^j\frac{D^j_{\mu_k}}{\gamma^{jp}}
\right]
8\pi a_ia_k G(\xi_i,\xi_k).
\qquad\qquad\qquad\qquad\quad\,
\end{eqnarray}
From  {\upshape(\ref{2.7})},
{\upshape(\ref{2.14})},
{\upshape(\ref{2.15})}
 and the Implicit Function Theorem
it follows that for
any  sufficiently small $\varepsilon$
and any points $\xi=(\xi_1,\ldots,\xi_m)\in\mathcal{O}_d$,
system {\upshape(\ref{2.22})} has a unique solution
$\mu=(\mu_1,\ldots,\mu_m)$ satisfying
  {\upshape(\ref{2.6})}. Moreover,
for any $i=1,\ldots,m$,
\begin{equation}\label{2.24}
\aligned
\log\big(8\mu_i^2\big)=\left\{
\frac{2(p-1)}{2-p}(1-\log8)+\frac{8\pi}{2-p}\left[ H(\xi_i,\xi_i)
+\sum\large_{k=1,\,k\neq i}^m
  a_ia_k  G(\xi_i,\xi_k)\right]
\right\}
\left[\,1+O\left(\frac{1}{|\log\varepsilon|}\right)\right].
\endaligned
\end{equation}

%
%
%
%
%
%
%

We  make  the change of variables
\begin{equation}\label{2.25}
\aligned
\upsilon(y)=p\gamma^{p-1}u(\varepsilon y)-p\gamma^{p},\,\,\,
\,\quad\forall\,\,y\in\Omega_{\varepsilon}:=\varepsilon^{-1}\Omega.
\endaligned
\end{equation}
From the definitions of
$\varepsilon$ and $\gamma$ in relations (\ref{1.5})-(\ref{2.7})
we can rewrite equation (\ref{1.1}) in the following form
\begin{equation}\label{2.26}
\begin{array}{ll}
\left\{\aligned
&-\Delta \upsilon =f(\upsilon)
\,\ \ \,\,
\,\,
\textrm{in}\,\,\,\,\,\Omega_\varepsilon,\\[1mm]
&\upsilon =-p\gamma^{p}
\quad\qquad
\textrm{on}\,\,\,\partial\Omega_\varepsilon,
\endaligned\right.
\end{array}
\end{equation}
where
\begin{equation}\label{2.27}
\aligned
f(\upsilon)=\left(1+\frac{\upsilon}{p\gamma^p}\right)
\left|1+\frac{\upsilon}{p\gamma^p}\right|^{p-2}
e^{\gamma^p\left(\left|1+\frac{\upsilon}{p\gamma^p}\right|^p-1\right)}.
\endaligned
\end{equation}
For equation (\ref{2.26}) we
write $\xi_i'=\xi_i/\varepsilon$, $i=1,\ldots,m$
and define its corresponding approximate solution  as
\begin{equation}\label{2.28}
\aligned
V_{\xi'}(y)=p\gamma^{p-1}U_\xi(\varepsilon y)-p\gamma^{p},
\endaligned
\end{equation}
with
$\xi'=(\xi_1',\ldots,\xi_m')$
and
$U_\xi$  defined in (\ref{2.16}).  What remains of this paper is to
 look for solutions of problem (\ref{2.26}) in the form
$\upsilon=V_{\xi'}+\phi$, where $\phi$ will represent a higher order correction.
In terms of $\phi$,
problem (\ref{2.26}) becomes
\begin{equation}\label{2.29}
\aligned
\left\{\aligned
&\mathcal{L}(\phi)=-\big[
E_{\xi'}+N(\phi)
\big]
\quad\textrm{in}\,\,\,\,\,\,\Omega_\varepsilon,\\[1mm]
& \phi =0
\qquad\qquad\qquad\qquad\,\,
\textrm{on}\,\,\,\,
\partial\Omega_\varepsilon,
\endaligned\right.\endaligned
\end{equation}
where
\begin{equation*}\label{2.30}
\aligned
\mathcal{L}(\phi):=-\Delta\phi-W_{\xi'}\phi\,\qquad
\ \,\textrm{with}\,\quad\,
W_{\xi'}:=f'(V_{\xi'}),
\endaligned
\end{equation*}
and
\begin{equation}\label{2.31}
\aligned
E_{\xi'}:=-\Delta V_{\xi'} -f(V_{\xi'}),
\quad\quad\qquad
N(\phi):=-\big[f(V_{\xi'}+\phi)-f(V_{\xi'})-f'(V_{\xi'})\phi\big].
\endaligned
\end{equation}
A key step in solving (\ref{2.29}), or equivalently (\ref{1.1}),  is that of
a solvability theory for the linear operator $\mathcal{L}$ under the configuration
space $\mathcal{O}_d$ of concentration points $\xi_i$. In developing this
theory, we will take into account the invariance, under
translations and dilations, of the problem
$\Delta e^v+e^v=0$ in $\mathbb{R}^2$. We will perform the
solvability theory for the linear operator $\mathcal{L}$ in a new weighted
$L^\infty$ space, following \cite{DKM,EMP}.
For any  $\xi=(\xi_1,\ldots,\xi_m)\in\mathcal{O}_d$ and $h\in L^\infty(\Omega_\varepsilon)$,
let us  introduce a  $L^\infty$-norm
$\|h\|_{*}:=\sup_{y\in\overline{\Omega}_\varepsilon}\big|\mathbf{H}_{\xi'}(y)h(y)\big|$ involving
 the new weighted
function
\begin{equation}\label{2.32}
\aligned
\mathbf{H}_{\xi'}(y)=\left[
\sum\limits_{i=1}^m\frac{\mu_i^\sigma}{(\mu_i^2+|y-\xi'_i|^2)^{(2+\sigma)/2}}
\right]^{-1},
\endaligned
\end{equation}
where $\sigma$ is   small   but  fixed, independent of $\varepsilon$,
such that $0<\sigma<\min\{(2-p)/p,1/2\}$.
With respect to the $\|\cdot\|_{*}$-norm,
the error term $E_{\xi'}$ defined in (\ref{2.31}) can be estimated as follows.

\vspace{1mm}
\vspace{1mm}
\vspace{1mm}
\vspace{1mm}

\noindent{\bf Proposition 2.2.}\,\,{\it
There exists a constant $C>0$ such that
for  any $\xi=(\xi_1,\ldots,\xi_m)\in\mathcal{O}_d$
and for  any $\varepsilon$ small enough,
\begin{equation}\label{2.33}
\aligned
\|E_{\xi'}\|_{*}\leq\frac{C}{\gamma^{4p}}=O\left(\frac{1}{|\log\varepsilon|^4}\right).
\endaligned
\end{equation}
}\noindent\begin{proof}
Observe that
\begin{eqnarray}\label{2.34}
-\Delta V_{\xi'}(y)
=-p\gamma^{p-1}\varepsilon^{2}\sum_{i=1}^m a_i\Delta \big(U_i+H_i\big)(\varepsilon y)
=-p\gamma^{p-1}\varepsilon^{2}\sum_{i=1}^ma_i\Delta U_i(\varepsilon y)
\quad\,
\nonumber\\
=\sum_{i=1}^m
a_i\,
e^{\omega_{\mu_i}\left(y-\xi'_i\right)}\left[
1
+
\sum_{j=1}^3\left(\frac{p-1}{p}\right)^j\frac{1}{\gamma^{jp}}
\big(
\omega^j_{\mu_i}-f^j_{\mu_i}
\big)
\right]\left(y-\xi'_i\right).
\end{eqnarray}
From (\ref{2.3}),  (\ref{2.6}) and (\ref{2.13})
we get that if $|y-\xi'_i|\geq d/\varepsilon$ for all $i=1,\ldots,m$,
$$
\aligned
\omega_{\mu_i}(y-\xi'_i)=4\log\varepsilon+O\left(1\right),\,\qquad\qquad\,
\omega_{\mu_i}^j(y-\xi'_i)=-D^{j}_{\mu_i}\log\varepsilon+O\left(1\right),
\quad
j=1, 2, 3,
\endaligned
$$
and then, by (\ref{2.10})-(\ref{2.12}),
\begin{equation}\label{2.35}
\aligned
\left|
-\Delta V_{\xi'}(y)
\right|
=\left[\sum_{i=1}^me^{\omega_{\mu_i}(y-\xi'_i)}
\right]
O\left(|\log\varepsilon|^3
\right).
\endaligned
\end{equation}
On the other hand, in the same region, by (\ref{2.20}) and (\ref{2.28}) we obtain
\begin{equation}\label{2.36}
\aligned
\left|
1+\frac{V_{\xi'}(y)}{p\gamma^p}
\right|
=
\left|
\frac{p\gamma^{p-1}U_{\xi}(\varepsilon y)}{p\gamma^p}
\right|
=O\left(\frac{1}{|\log\varepsilon|}
\right),
\endaligned
\end{equation}
then
\begin{equation}\label{2.37}
\aligned
\left|f(V_{\xi'})\right|
=
\left|1+\frac{V_{\xi'}}{p\gamma^p}\right|^{p-1}
e^{\gamma^p\left(\left|1+\frac{V_{\xi'}}{p\gamma^p}\right|^p-1\right)}
=\frac{O(\varepsilon^{\frac{2+p}{p}})}{\,|\log\varepsilon|^{p-1}\,}
\exp\left[
-\frac{2-p}{p}|\log\varepsilon|+
O\left(\frac{1}{|\log\varepsilon|^{p-1}}\right)
\right],
\endaligned
\end{equation}
which, together with    (\ref{2.32}) and (\ref{2.35}), implies that for any  $0<p<2$,
\begin{eqnarray}\label{2.38}
\left|\mathbf{H}_{\xi'}(y)
E_{\xi'}(y)
\right|
\leq C\left\{
\sum\limits_{i=1}^m\frac{|\log\varepsilon|^3}{(\mu_i^2+|y-\xi'_i|^2)^{(2-\sigma)/2}}
+
\frac{\varepsilon^{\frac{2-p}{p}-\sigma}}{|\log\varepsilon|^{p-1}}
\exp\left[
-\frac{2-p}{p}|\log\varepsilon|+
O\left(
\frac{1}{|\log\varepsilon|^{p-1}}
\right)
\right]
\right\}
=o\left(\frac{1}{\gamma^{4p}}
\right).
\end{eqnarray}
Let us  fix an index $i\in\{1,\ldots,m\}$
and the region $|y-\xi'_i|\leq d/\varepsilon^\theta$ with
any $\theta<1$ but close enough to $1$.
From (\ref{2.13}), (\ref{2.21}), (\ref{2.28}) and the Taylor expansion we have that
in the ball  $|y-\xi'_i|<\mu_i|\log\varepsilon|^\tau$ with  $\tau\geq 10$ sufficiently large but fixed,
$$
\aligned
\left(1+\frac{V_{\xi'}}{p\gamma^p}\right)
\left|
1+\frac{V_{\xi'}}{p\gamma^p}
\right|^{p-2}
=a_i&\left\{1+\frac{p-1}{p}\frac1{\gamma^p}\underbrace{\omega_{\mu_i}(y-\xi'_i)}\limits_{A_1}
+\left(\frac{p-1}{p}\right)^2\frac1{\gamma^{2p}}\underbrace{
\left[
\omega^{1}_{\mu_i}+\frac{p-2}{2(p-1)}(\omega_{\mu_i})^2\right]
(y-\xi'_i)}\limits_{A_2}
\right.\\
&\,\,\,\,
+\left(\frac{p-1}{p}\right)^3\frac1{\gamma^{3p}}\underbrace{\left[
\omega^{2}_{\mu_i}+\frac{p-2}{p-1}\omega_{\mu_i}\omega^1_{\mu_i}
+\frac{(p-2)(p-3)}{6(p-1)^2}(\omega_{\mu_i})^3
\right]
(y-\xi'_i)}\limits_{A_3}
\\
&\left.\,\,\,\,
+\,O\left(\frac{\log^4(\mu_i+|y-\xi'_i|)}{\gamma^{4p}}\right)
\right\},
\endaligned
$$
and
\begin{eqnarray}\label{2.380}
\gamma^p\left(\left|1+\frac{V_{\xi'}}{p\gamma^p}\right|^{p}-1\right)
=\underbrace{\omega_{\mu_i}(y-\xi'_i)}\limits_{A_1}+\frac{p-1}{p}\frac1{\gamma^p}
\underbrace{\left[\omega^{1}_{\mu_i}+\frac{1}{2}(\omega_{\mu_i})^2\right](y-\xi'_i)}\limits_{B_1}
\qquad\qquad\qquad\qquad\qquad\qquad\qquad
\qquad\qquad\qquad\ \,\,
&&
\nonumber
\\
+\left(\frac{p-1}{p}\right)^2\frac1{\gamma^{2p}}\underbrace{
\left[
\omega^{2}_{\mu_i}+\omega_{\mu_i}\omega^{1}_{\mu_i}+\frac{p-2}{6(p-1)}(\omega_{\mu_i})^3\right]
(y-\xi'_i)}\limits_{B_2}
\qquad\qquad\qquad\qquad\qquad\qquad\qquad\qquad\,\,
&&
\nonumber
\\
+\left(\frac{p-1}{p}\right)^3\frac1{\gamma^{3p}}\underbrace{
\left[
\omega^{3}_{\mu_i}
+\frac12(\omega^{1}_{\mu_i})^2
+
\omega_{\mu_i}\omega^{2}_{\mu_i}
+\frac{p-2}{2(p-1)}(\omega_{\mu_i})^2\omega^{1}_{\mu_i}
+\frac{(p-2)(p-3)}{24(p-1)^2}(\omega_{\mu_i})^4
\right]
(y-\xi'_i)}\limits_{B_3}
&&
\nonumber
\\
+\,
O
\left(\frac{\log^5(\mu_i+|y-\xi'_i|)}{\gamma^{4p}}\right).
\qquad\qquad\qquad\qquad\qquad\qquad\qquad\quad
\qquad\qquad\qquad\qquad\qquad\qquad\qquad\ \,\,
&&
\end{eqnarray}
Then
%
%
%
%
\begin{eqnarray}\label{2.381}
e^{\gamma^p\left[\left(1+\frac{V_{\xi'}}{p\gamma^p}\right)^{p}-1\right]}
=e^{\omega_{\mu_i}\left(y-\xi'_i\right)}
\left\{1+\frac{p-1}{p}\frac1{\gamma^p}
B_1
+\left(\frac{p-1}{p}\right)^2\frac1{\gamma^{2p}}
\left[B_2+\frac12(B_1)^2\right]
\right.
\quad\qquad\quad
\nonumber\\[1mm]
\left.
+\left(\frac{p-1}{p}\right)^3\frac1{\gamma^{3p}}
\left[B_3+B_1B_2+\frac16(B_1)^3\right]
+O\left(\frac{\log^5(\mu_i+|y-\xi'_i|)}{\gamma^{4p}}\right)
\right\}.
\end{eqnarray}
Thus by  the definition of $f(\cdot)$ in (\ref{2.27}) and  the definitions of $f^j_{\mu_i}$,
$j=1, 2, 3$ in (\ref{2.10})-(\ref{2.12}),
%
%
\begin{eqnarray}\label{2.39}
f(V_{\xi'})=a_i\,
e^{\omega_{\mu_i}\left(y-\xi'_i\right)}
\left\{\,1\,+\,\frac{p-1}{p}\frac1{\gamma^p}
\underbrace{\left(A_1+B_1\right)}\limits_{=\big(
\omega^1_{\mu_i}-f^1_{\mu_i}
\big)\left(y-\xi'_i\right)}
+\,\,\left(\frac{p-1}{p}\right)^2\frac1{\gamma^{2p}}
\underbrace{\left[A_2+A_1B_1+
B_2+\frac12(B_1)^2\right]}\limits_{=\big(
\omega^2_{\mu_i}-f^2_{\mu_i}
\big)\left(y-\xi'_i\right)}
\right.
&&
\nonumber\\[1mm]
+\,\left(\frac{p-1}{p}\right)^3\frac1{\gamma^{3p}}
\underbrace{\left[
A_3+
A_2B_1
+A_1\left(B_2+\frac12(B_1)^2\right)+
B_3+B_1B_2
+\frac16(B_1)^3\right]}\limits_{=\big(
\omega^3_{\mu_i}-f^3_{\mu_i}
\big)\left(y-\xi'_i\right)}
\ \,
&&
\nonumber\\[1mm]
\left.
+\,O\left(\frac{\log^8(\mu_i+|y-\xi'_i|)}{\gamma^{4p}}\right)
\right\}.
\qquad\qquad\qquad\qquad\qquad\qquad
\qquad\qquad\qquad\qquad\quad\,\
&&
\end{eqnarray}
From (\ref{2.34}) and (\ref{2.39})
  we obtain that
in the region
$|y-\xi'_i|<\mu_i|\log\varepsilon|^\tau$,
\begin{equation*}\label{2.40}
\aligned
E_{\xi'}=-\Delta V_{\xi'} -f(V_{\xi'})
=a_i\,e^{\omega_{\mu_i}\left(y-\xi'_i\right)}
O\left(\frac{\log^8(\mu_i+|y-\xi'_i|)}{\gamma^{4p}}\right),
\endaligned
\end{equation*}
then  by   (\ref{2.32}),
\begin{eqnarray}\label{2.41}
\left|\mathbf{H}_{\xi'}(y)
E_{\xi'}(y)
\right|
\leq \frac{C}{\gamma^{4p}}\frac{\,\log^8\big(\mu_i+|y-\xi'_i|\big)\,}
{(\mu_i^2+|y-\xi'_i|^2)^{(2-\sigma)/2}}
=O\left( \frac{1}{\gamma^{4p}}\right).
\end{eqnarray}
In the remaining region
$\mu_i|\log\varepsilon|^\tau\leq
|y-\xi'_i|\leq d/\varepsilon^\theta$ with
any $\theta<1$ but close enough to $1$
and any   $\tau\geq10$ sufficiently large but fixed,
by  (\ref{2.10})-(\ref{2.13}) and (\ref{2.34}) we  get
that there exists a constant $D>0$, independent of every $\theta<1$, such that
\begin{equation}\label{2.42}
\aligned
\left|-\Delta V_{\xi'}(y)
\right|\leq D |\log\varepsilon|^3
e^{\omega_{\mu_i}(y-\xi'_i)}.
\endaligned
\end{equation}
In  this   region,
by  (\ref{2.6}),  (\ref{2.13}), (\ref{2.21}) and (\ref{2.28}) we can compute
$$
\aligned
1+\frac{V_{\xi'}}{p\gamma^p}=&a_i+\frac{a_i}{p\gamma^p}\left\{
\left[
1-
\frac{1}{4}
\sum_{j=1}^3\left(\frac{p-1}{p}\right)^j\frac{D^j_{\mu_i}}{\gamma^{jp}}
\right]\omega_{\mu_i}\left(y-\xi'_i\right)
+\left[\sum_{j=1}^3\left(\frac{p-1}{p}\right)^j\frac{D^j_{\mu_i}}{\gamma^{jp}}
\right]\frac14\log\left(\frac{8}{\mu_i^2}\right)
+O\left(\frac{\mu_i}{\mu_i+|y-\xi'_i|}\right)\right\},
\endaligned
$$
then
\begin{equation}\label{2.45}
\aligned
1-\theta+
\frac{1}{4|\log\varepsilon|}
\left[
\log\left(\frac{8\mu_i^2}{d^4}\right)
+\frac{(p-1)\theta D_{\mu_i}^1}{4}
\right]
+
O\left(\frac{1}{|\log\varepsilon|^2}\right)
\leq
\left|
1+\frac{V_{\xi'}}{p\gamma^p}
\right|
\leq
1-\tau\frac{\log|\log\varepsilon|}{|\log\varepsilon|}
+
O\left(\frac{1}{|\log\varepsilon|}\right),
\endaligned
\end{equation}
which, together with the Taylor expansion,  implies that  there exists a constant $D>0$, independent of every $\theta<1$, such that
\begin{equation}\label{2.43}
\aligned
\left|1+\frac{V_{\xi'}}{p\gamma^p}\right|^{p-1}
\leq D\left(1+\frac{1}{|\log\varepsilon|^{p-1}}\right),
\endaligned
\end{equation}
and
\begin{equation}\label{2.46}
\aligned
e^{\gamma^p\left(\left|1+\frac{V_{\xi'}}{p\gamma^p}\right|^p-1\right)}\leq
De^{\left[
1-
\frac{1}{4}
\sum\limits_{j=1}^3\left(\frac{p-1}{p}\right)^j\frac{D^j_{\mu_i}}{\gamma^{jp}}
\right]
\omega_{\mu_i}\left(y-\xi'_i\right)
}
=O\left(
e^{
\omega_{\mu_i}\left(y-\xi'_i\right)
}
\right).
\endaligned
\end{equation}
Hence in the region $\mu_i|\log\varepsilon|^\tau\leq
|y-\xi'_i|\leq d/\varepsilon^\theta$, by   (\ref{2.32}), (\ref{2.42}), (\ref{2.43}) and (\ref{2.46}),
$$
\aligned
\left|\mathbf{H}_{\xi'}(y)
E_{\xi'}(y)
\right|
\leq C\left(
|\log\varepsilon|^3
+\frac{1}{|\log\varepsilon|^{p-1}}\right)
\frac{1}{(\mu_i^2+|y-\xi'_i|^2)^{(2-\sigma)/2}}
=o\left(\frac{1}{\gamma^{4p}}
\right),
\endaligned
$$
which, together with (\ref{2.38}) and (\ref{2.41}), establishes  the validity of
 estimate  (\ref{2.33}).
\end{proof}

\vspace{1mm}

\section{Analysis of the linearized operator}
In this section we give the solvability theory for the linear operator $\mathcal{L}$,
uniformly on $\xi\in\mathcal{O}_d$,
under a $L^{\infty}$-norm  involving the weighted function (\ref{2.32}).
Recall that
$\mathcal{L}(\phi)=-\Delta \phi -W_{\xi'}\phi$
with $W_{\xi'}=f'(V_{\xi'})$.
As in Proposition 2.2, we have for  $W_{\xi'}$
and $f''(V_{\xi'})$  the following  expansions.

\vspace{1mm}
\vspace{1mm}
\vspace{1mm}
\vspace{1mm}

\noindent{\bf Proposition 3.1.}\,\,{\it
There exists a constant $D_0>0$ such that
for any $\xi=(\xi_1,\ldots,\xi_m)\in\mathcal{O}_d$ and
for any $\varepsilon$ small enough,
\begin{equation}\label{3.1}
\aligned
\big|W_{\xi'}(y)\big|
\leq D_0\sum_{i=1}^me^{
\omega_{\mu_i}\left(y-\xi'_i\right)
}
\qquad\qquad
\textrm{and}
\qquad\qquad
\big|f''(V_{\xi'})\big|\leq D_0\sum_{i=1}^me^{
\omega_{\mu_i}\left(y-\xi'_i\right)
},
\endaligned
\end{equation}
uniformly in each region   $\mu_i|\log\varepsilon|^\tau\leq|y-\xi'_i|\leq d/\varepsilon^\theta$ with
any $\theta<1$ but close enough to $1$ and any   $\tau\geq10$ sufficiently large but fixed.
Moreover,  if $|y-\xi'_i|<\mu_i|\log\varepsilon|^\tau$,
\begin{equation}\label{3.2}
\aligned
W_{\xi'}(y)=\frac{8\mu_i^2}{(\mu_i^2+|y-\xi'_i|^2)^2}
\left\{1+\frac{p-1}{p}\frac1{\gamma^p}
\left[1+\omega^{1}_{\mu_i}+\frac{1}{2}(\omega_{\mu_i})^2+2\omega_{\mu_i}\right]\left(y-\xi'_i\right)
+O\left(\frac{\log^4(\mu_i+|y-\xi'_i|)}{\gamma^{2p}}\right)
\right\}.
\endaligned
\end{equation}
In addition,
\begin{equation}\label{3.3}
\aligned
\|W_{\xi'}\|_{*}\leq C
\,\,\qquad\quad\ \,\,
\textrm{and}
\,\,\qquad\quad\ \,\,
\|f''(V_{\xi'})\|_{*}\leq C.
\endaligned
\end{equation}}\noindent\begin{proof}
For the sake of simplicity, we prove the estimates for the potential $W_{\xi'}$ only.
From (\ref{2.27}) we can get
$$
\aligned
W_{\xi'}=\frac{p-1}{p}\frac1{\gamma^p}\left|1+\frac{V_{\xi'}}{p\gamma^p}\right|^{p-2}
e^{\gamma^p\left(\left|1+\frac{V_{\xi'}}{p\gamma^p}\right|^p-1\right)}
+
\left|1+\frac{V_{\xi'}}{p\gamma^p}\right|^{2(p-1)}
e^{\gamma^p\left(\left|1+\frac{V_{\xi'}}{p\gamma^p}\right|^p-1\right)}
:=I+J.
\endaligned
$$
If $|y-\xi'_i|<\mu_i|\log\varepsilon|^\tau$ for some $i=1,\ldots,m$ and  $\tau\geq10$ sufficiently large but fixed,
by  (\ref{2.381})
and the Taylor expansion we get
$$
\aligned
I
=& e^{\omega_{\mu_i}\left(y-\xi'_i\right)}\left\{1+\frac{p-1}{p}\frac1{\gamma^p}\left[\omega^{1}_{\mu_i}+\frac{1}{2}(\omega_{\mu_i})^2\right]\left(y-\xi'_i\right)
+O\left(\frac{\log^4(\mu_i+|y-\xi'_i|)}{\gamma^{2p}}\right)
\right\}\\[1mm]
&\times\frac{p-2}{p}\frac1{\gamma^p}\left[\frac{p-1}{p-2}+\frac{p-1}{p}\frac1{\gamma^p}\omega_{\mu_i}\left(y-\xi'_i\right)
+O\left(\frac{\log^2(\mu_i+|y-\xi'_i|)}{\gamma^{2p}}\right)
\right],
\endaligned
$$
and
$$
\aligned
J
=&e^{\omega_{\mu_i}\left(y-\xi'_i\right)}\left\{1+\frac{p-1}{p}\frac1{\gamma^p}
\left[\omega^{1}_{\mu_i}+\frac{1}{2}(\omega_{\mu_i})^2\right]\left(y-\xi'_i\right)
+O\left(\frac{\log^4(\mu_i+|y-\xi'_i|)}{\gamma^{2p}}\right)
\right\}
\\
&
\times
\left[1+\frac{p-1}{p}\frac2{\gamma^p}\omega_{\mu_i}\left(y-\xi'_i\right)
+O\left(\frac{\log^2(\mu_i+|y-\xi'_i|)}{\gamma^{2p}}\right)
\right],
\endaligned
$$
and hence
\begin{equation}\label{3.4}
\aligned
W_{\xi'}(y)=e^{\omega_{\mu_i}\left(y-\xi'_i\right)}
\left\{1+\frac{p-1}{p}\frac1{\gamma^p}
\left[1+\omega^{1}_{\mu_i}+\frac{1}{2}(\omega_{\mu_i})^2+2\omega_{\mu_i}\right]\left(y-\xi'_i\right)
+O\left(\frac{\log^4(\mu_i+|y-\xi'_i|)}{\gamma^{2p}}\right)
\right\}.
\endaligned
\end{equation}
But if  $\mu_i|\log\varepsilon|^\tau\leq
|y-\xi'_i|\leq d/\varepsilon^\theta$
with  any $\theta<1$ but close enough to $1$, by  (\ref{2.45})
we find
$$
\aligned
\left|1+\frac{V_{\xi'}(y)}{p\gamma^p}\right|^{p-2}
=O\left(1\right)
\qquad\qquad
\textrm{and}
\qquad\qquad
\left|1+\frac{V_{\xi'}(y)}{p\gamma^p}\right|^{2(p-1)}
=O\left(1\right),
\endaligned
$$
and  by (\ref{2.46}),
\begin{equation}\label{3.5}
\aligned
\big|W_{\xi'}(y)\big|\leq
Ce^{\gamma^p\left(\left|1+\frac{V_{\xi'}}{p\gamma^p}\right|^p-1\right)}
=
O\left(
e^{
\omega_{\mu_i}\left(y-\xi'_i\right)
}
\right).
\endaligned
\end{equation}
Finally, if
$|y-\xi'_i|\geq d/\varepsilon$
for all $i=1,\ldots,m$, by  (\ref{2.36}) we give
\begin{equation}\label{3.6}
\aligned
\big|W_{\xi'}(y)\big|=
\left(\frac{1}
{|\log\varepsilon|^{p-1}}
+\frac{1}
{|\log\varepsilon|^{2(p-1)}}\right)O(\varepsilon^{\frac{2+p}p})
\exp\left[
-\frac{2-p}{p}|\log\varepsilon|+
O\left(\frac{
1
}{|\log\varepsilon|^{p-1}}\right)\right].
\endaligned
\end{equation}
From   (\ref{3.4})-(\ref{3.6}) and the definition of $\left\|\cdot\right\|_{*}$
involving the  weighted
function
(\ref{2.32}), we easily prove the first estimate in  (\ref{3.3}).
\end{proof}

\vspace{1mm}
\vspace{1mm}

Set
\begin{equation}\label{3.7}
\aligned
Z_{0}(z)=\frac{|z|^2-1}{|z|^2+1},
\,\,\quad\qquad\quad\,\,
Z_{j}(z)=\frac{4z_j}{|z|^2+1},\,\,\,\,j=1,\,2.
\endaligned
\end{equation}
It is well known that any bounded solution of
\begin{equation}\label{3.8}
\aligned
\Delta\phi+\frac{8}{(1+|z|^2)^2}\phi=0
\,\quad\,\textrm{in}\,\,\,\,\mathbb{R}^2
\endaligned
\end{equation}
is a linear combination of $Z_j$, $j=0,1,2$ (see \cite{BP}).
Given $h\in C(\overline{\Omega}_\varepsilon)$ and points
$\xi=(\xi_1,\ldots,\xi_m)\in\mathcal{O}_d$,
we will solve
the following linear projected problem of finding a function $\phi\in
H^{2}(\Omega_\varepsilon)$
such that
\begin{equation}\label{3.9}
\left\{\aligned
&\mathcal{L}(\phi)=-\Delta \phi -W_{\xi'}\phi=h
+ \sum\limits_{i=1}^m\sum\limits_{j=1}^{2}c_{ij}e^{
\omega_{\mu_i}\left(y-\xi'_i\right)}Z_{ij}\,\,\ \,
\,\textrm{in}\,\,\,\,\,\,\Omega_\varepsilon,\\
&\phi=0\,\,\,\,\,\,\,\,
\ \ \ \ \ \ \ \ \ \ \,\,\,
\qquad\qquad\qquad\quad\qquad\qquad\qquad\qquad
\ \,\textrm{on}\,\,\,\,\partial\Omega_{\varepsilon},\\[1mm]
&\int_{\Omega_\varepsilon}
e^{
\omega_{\mu_i}\left(y-\xi'_i\right)}
Z_{ij}\phi=0
\,\qquad\qquad
\qquad
\forall\,\,i=1,\ldots,m,\,\,\,j=1, 2,
\endaligned\right.
\end{equation}
for some  coefficients
$c_{ij}\in\mathbb{R}$, $i=1,\ldots,m$ and $j=1,2$.
Here and in the sequel,
for any  $i=1,\ldots,m$ and $j=0,1,2$, we define
\begin{equation}\label{3.10}
\aligned
Z_{ij}(y):=Z_j\left(\frac{y-\xi'_i}{\mu_i}\right)=\left\{\aligned
&\frac{|y-\xi'_i|^2-\mu_i^2}{|y-\xi'_i|^2+\mu_i^2}
\quad\,\,\,\,\textrm{if}\,\,\,j=0,\\[1mm]
&\frac{4\mu_i(y-\xi'_i)_j}{|y-\xi'_i|^2+\mu_i^2}
\quad\,\,\,\,\textrm{if}\,\,\,j=1,2.
\endaligned\right.
\endaligned
\end{equation}

\vspace{1mm}
\vspace{1mm}

\noindent {\bf Proposition 3.2.}\,\,{\it
There exist constants $C>0$  and $\varepsilon_0>0$ such
that for any  $0<\varepsilon<\varepsilon_0$,    any points
$\xi=(\xi_1,\ldots,\xi_m)\in\mathcal{O}_d$ and
any $h\in C(\overline{\Omega}_\varepsilon)$,
problem {\upshape(\ref{3.9})}  admits a unique solution
$\phi=\mathcal{T}(h)\in
H^{2}(\Omega_\varepsilon)$ for
some coefficients
 $c_{ij}\in\mathbb{R}$,  $i=1,\ldots,m$, $j=1,2$,
which defines a linear operator of $h$ and satisfies the estimate
\begin{equation}\label{3.11}
\aligned
\|\mathcal{T}(h)\|_{L^{\infty}(\Omega_\varepsilon)}\leq C|\log\varepsilon|\,\|h\|_{*}.
\endaligned
\end{equation}

}\noindent{\it Proof.}\,
The proof of this result  will be split into six steps which we state and prove next.

{\bf Step 1:} The operator $\mathcal{L}$ satisfies the maximum principle in
$\widetilde{\Omega}_\varepsilon:=\Omega_\varepsilon\setminus\bigcup_{i=1}^mB(\xi'_i,R\mu_i)$
for $R$ large but independent of $\varepsilon$.
Specifically,
$$
\aligned
\textrm{if}
\,\,\,
\mathcal{L}(\psi)=
-\Delta\psi-W_{\xi'}\psi\geq0
\,\ \,\textrm{in}\,\,\,\widetilde{\Omega}_\varepsilon
\qquad
\textrm{and}
\qquad
\psi\geq0
\,\ \,\textrm{on}\,\,\,\partial\widetilde{\Omega}_\varepsilon,
\quad\quad
\textrm{then}
\quad
\psi\geq0
\,\ \,\textrm{in}\,\,\,\widetilde{\Omega}_\varepsilon.
\endaligned
$$
To  prove this, it is sufficient  to give a positive function $Z$
in $\widetilde{\Omega}_\varepsilon$ such that
$\mathcal{L}(Z)>0$.
Indeed, let
$$
\aligned
Z(y)=\sum_{i=1}^m\left[2-
\frac{\mu_i^{\sigma}}{(a|y-\xi'_i|)^\sigma}\right],\,\quad\,a>0.
\endaligned
$$
Clearly, if $|y-\xi'_i|\geq R\mu_i$ for all
 $i=1,\ldots,m$ and $R>1/a$,
then $m<Z(y)<2m$. Moreover,
in each region  $R\mu_i\leq
|y-\xi'_i|\leq d/\varepsilon^\theta$
with  any $\theta<1$ but close enough to $1$, by (\ref{2.4}),  (\ref{3.1}) and (\ref{3.2})
we find
$$
\aligned
\mathcal{L}(Z)=&-\Delta Z-W_{\xi'}Z
\geq\sum_{i=1}^m
\frac{\sigma^2\mu_i^{\sigma}}{a^\sigma|y-\xi'_i|^{2+\sigma}}
-2mD_0
\sum_{i=1}^m
\frac{8\mu_i^2}{(\mu_i^2+|y-\xi'_i|^2)^2}
\geq
\sum_{i=1}^m
\frac{\sigma^2\mu_i^{\sigma}}{a^\sigma|y-\xi'_i|^{2+\sigma}}
 -2mD_0\sum_{i=1}^m\frac{8\mu_i^2}{|y-\mu_i|^4}\\
\geq&
\sum_{i=1}^m
\frac{\sigma^2\mu_i^{\sigma}}{a^\sigma|y-\xi'_i|^{2+\sigma}}
\left[1
 -\frac{2mD_0a^\sigma}{\sigma^2}\frac{8}{R^{2-\sigma}}
 \right]>0
\endaligned
$$
provided $R>(16 m D_0 a^\sigma/\sigma^2)^{1/(2-\sigma)}$,
where $D_0$ is the constant in  Proposition 3.1.
As in the remaining region
$|y-\xi'_i|\geq d/\varepsilon$
for all $i=1,\ldots,m$,  by  (\ref{3.6}) we have that there exist
positive constants $C_1$ and $C_2$ such that
for any $\varepsilon$ small enough,
\begin{eqnarray}\label{3.41}
\mathcal{L}(Z)=-\Delta Z-W_{\xi'}Z
\geq\sum_{i=1}^m
\frac{\sigma^2\mu_i^{\sigma}}{a^\sigma|y-\xi'_i|^{2+\sigma}}
-2mC_{1}
\left(\frac{\varepsilon^{\frac4p}}
{|\log\varepsilon|^{p-1}}
+\frac{\varepsilon^{\frac4p}}
{|\log\varepsilon|^{2(p-1)}}\right)
\exp\left[O\left(\frac{
1
}{|\log\varepsilon|^{p-1}}\right)\right]
\ \,\quad
&&
\nonumber
\\[1mm]
\geq
\,
\varepsilon^{2+\sigma}
\left\{
\sum_{i=1}^m
\frac{\sigma^2\mu_i^\sigma}{a^\sigma C_{2}^{2+\sigma}}
-2mC_{1}\left(\frac{\varepsilon^{\frac{2-p}{p}-\sigma}}
{|\log\varepsilon|^{p-1}}
+\frac{\varepsilon^{\frac{2-p}{p}-\sigma}}
{|\log\varepsilon|^{2(p-1)}}\right)
\exp\left[-\frac{2-p}{p}
|\log\varepsilon|
+O\left(\frac{
1
}{|\log\varepsilon|^{p-1}}\right)\right]
\right\}
&&\nonumber
\\[1mm]
\geq
\,
\frac{1}{2}
\varepsilon^{2+\sigma}
\left[
\sum_{i=1}^m
\frac{\sigma^2\mu_i^\sigma}{a^\sigma C_{2}^{2+\sigma}}
\right]>0
\qquad\qquad\qquad\qquad\qquad\qquad\qquad
\qquad\qquad\qquad\qquad\qquad\qquad\qquad
\qquad\quad\
&&
\end{eqnarray}
because of $0<p<2$ and
$0<\sigma<\min\{(2-p)/p,1/2\}$.
The function $Z(x)$ is what we are looking for.

{\bf Step 2:} Let $R$ be as before. We define the ``inner norm'' of $\phi$ as
\begin{equation}\label{3.12}
\aligned
\|\phi\|_{**}=\sup_{y\in\bigcup_{i=1}^m\overline{B(\xi'_i,R\mu_i)}}|\phi(y)|
\endaligned
\end{equation}
and claim that there is a constant $C>0$ independent of $\varepsilon$ such that if $\mathcal{L}(\phi)=h$
in $\Omega_\varepsilon$, $\phi=0$ on $\po_\varepsilon$,
then
\begin{equation}\label{3.13}
\aligned
\|\phi\|_{L^{\infty}(\Omega_\varepsilon)}\leq C\left(\|\phi\|_{**}+
\|h\|_{*}
\right)
\endaligned
\end{equation}
for any $h\in C^{0,\alpha}(\overline{\Omega}_\varepsilon)$.
We will establish this estimate with the help of suitable barriers.
Let $M$ be a large number such that
$\Omega_\varepsilon\subset B(\xi'_i,M/\varepsilon)$ for all $i=1,\ldots,m$. Consider
the solution $\psi_i$ of the
problem
$$
\left\{\begin{aligned}
&-\Delta\psi_i=
\frac{2\mu_i^\sigma}{|y-\xi'_i|^{2+\sigma}}
\,\,\quad\quad\,\,
\ \ \,\,\textrm{in}\,\,
\quad\,\,
R\mu_i<\left|y-\xi'_i\right|
<\frac{M}{\varepsilon},\\[1mm]
&\psi_i(y)=0\,\ \ \ \ \
\ \,\textrm{on}
\quad
\left|y-\xi'_i\right|=R\mu_i
\quad
\textrm{and}
\quad
\left|y-\xi'_i\right|=\frac{M}{\varepsilon}.
\end{aligned}
\right.
$$
Namely,
the function  $\psi_i(y)$ is
the positive function  given by
$$
\aligned
\psi_i(y)=&-\frac2{\sigma^2}
\left(\frac{\mu_i^\sigma}{|y-\xi'_i|^{\sigma}}
-
\frac1{R^{\sigma}}
\right)
+
\frac2{\sigma^2}
\left(
\frac{(\varepsilon\mu_i)^{\sigma}}{M^{\sigma}}
-\frac1{R^{\sigma}}
\right)
\frac{\,1\,}{\,\log\frac{M}{\,R\varepsilon\mu_i}\,}\log\left|\frac{y-\xi'_i}{R\mu_i}\right|.
\endaligned
$$
Clearly,  the function  $\psi_i$ is
uniformly bounded from above by the constant  $2/(\sigma^2R^\sigma)$.
We take the barrier
$$
\aligned
\widetilde{\phi}(y)=\|\phi\|_{**}Z(y)+
\|h\|_{*}\sum_{i=1}^m\psi_i(y),
\endaligned
$$
where $Z$ was defined in the previous step.
First of all, observe that
by the definition of $Z$,
choosing $R$ larger if necessary
$$
\aligned
\widetilde{\phi}(y)\geq
\|\phi\|_{**}Z(y)\geq m\|\phi\|_{**}\geq|\phi(y)|
\qquad
\textrm{for}
\,\,\,\,
|y-\xi'_i|=R\mu_i,
\,\,\,
i=1,\ldots,m,
\endaligned
$$
and by  the positivity of $Z(y)$ and $\psi_i(y)$,
$$
\aligned
\widetilde{\phi}(y)
\geq
0=|\phi(y)|
\quad\,
\textrm{for}
\,\,\,\,
y\in\partial\Omega_\varepsilon.
\endaligned
$$
From the definition of $\|\cdot\|_*$ involving the weighted function
 (\ref{2.32}) we know that
\begin{equation}\label{3.14}
\aligned
\|h\|_{*}\left[
\,\sum_{i=1}^m
\frac{\mu_i^\sigma}{|y-\xi'_i|^{2+\sigma}}\right]
\geq
\|h\|_{*}
\left[
\,
\sum\limits_{i=1}^m\frac{\mu_i^\sigma}{(\mu_i^2+|y-\xi'_i|^2)^{(2+\sigma)/2}}
\right]
\geq|h(y)|,
\endaligned
\end{equation}
then, if $R\mu_i\leq
|y-\xi'_i|\leq d/\varepsilon^\theta$
with  any $\theta<1$ but close enough to $1$, by
the expansions of $W_{\xi'}$  in (\ref{3.1})-(\ref{3.2}) we conclude
$$
\aligned
L(\widetilde{\phi})
&\geq
\|h\|_{*}
\sum_{i=1}^m \mathcal{L}(\psi_i)(y)=
\|h\|_{*}
\sum_{i=1}^m\left[
\frac{2\mu_i^\sigma}{|y-\xi'_i|^{2+\sigma}}
-W_{\xi'}(y)\psi_i(y)
\right]\\
&\geq
\|h\|_{*}
\sum_{i=1}^m\left[
\frac{2\mu_i^\sigma}{|y-\xi'_i|^{2+\sigma}}
-\frac{2mD_0}{\sigma^2R^\sigma}\frac{8\mu_i^2}{(\mu_i^2+|y-\xi'_i|^2)^2}
\right]\\
&\geq
\|h\|_{*}
\left[\,\sum_{i=1}^m
\frac{\mu_i^\sigma}{|y-\xi'_i|^{2+\sigma}}\right]
\geq|h(y)|\geq|\mathcal{L}(\phi)(y)|
\endaligned
$$
provided $R>\sqrt{16mD_0}/\sigma$ and $\varepsilon$ small enough,
while if
$|y-\xi'_i|\geq d/\varepsilon$
for all $i=1,\ldots,m$, similar to (\ref{3.41}),  by  (\ref{3.6})  we get
$$
\aligned
L(\widetilde{\phi})
&\geq
\|h\|_{*}
\sum_{i=1}^m \mathcal{L}(\psi_i)(y)=
\|h\|_{*}
\sum_{i=1}^m\left[
\frac{2\mu_i^\sigma}{|y-\xi'_i|^{2+\sigma}}
-W_{\xi'}(y)\psi_i(y)
\right]\\
&\geq
\|h\|_{*}
\sum_{i=1}^m\left\{
\frac{2\mu_i^\sigma}{|y-\xi'_i|^{2+\sigma}}
-\frac{2C_1}{\sigma^2R^\sigma}
\left(\frac{\varepsilon^{\frac4p}}
{|\log\varepsilon|^{p-1}}
+\frac{\varepsilon^{\frac4p}}
{|\log\varepsilon|^{2(p-1)}}\right)
\exp\left[O\left(\frac{
1
}{|\log\varepsilon|^{p-1}}\right)\right]
\right\}\\
&\geq
\|h\|_{*}
\left[\,
\sum_{i=1}^m
\frac{\mu_i^\sigma}{|y-\xi'_i|^{2+\sigma}}\right]
\geq|h(y)|\geq|\mathcal{L}(\phi)(y)|.
\endaligned
$$
Therefore, by the maximum principle in Step 1, we obtain
$$
\aligned
|\phi(y)|\leq\widetilde{\phi}(y)
\quad\textrm{for}\,\,\,y\in\widetilde{\Omega}_\varepsilon.
\endaligned
$$
Since $Z(y)\leq 2m$ and $\psi_i(y)\leq 2/(\sigma^2R^\sigma)$ in $\widetilde{\Omega}_\varepsilon$, we arrive at
$$
\aligned
\|\phi\|_{L^{\infty}(\Omega_\varepsilon)}\leq C\left(\|\phi\|_{**}+
\|h\|_{*}
\right).
\endaligned
$$\indent{\bf Step 3:}
We prove uniform a priori estimates for solutions
$\phi$ of problem  $\mathcal{L}(\phi)=h$
in $\Omega_\varepsilon$,
$\phi=0$ on $\po_\varepsilon$,
where $h\in C^{0,\alpha}(\overline{\Omega}_\varepsilon)$
and in addition
$\phi$ satisfies the orthogonality conditions:
\begin{equation}\label{3.15}
\aligned
\int_{\Omega_\varepsilon}
e^{
\omega_{\mu_i}\left(y-\xi'_i\right)}
Z_{ij}\phi=0
\,\,\,\,\,\quad\,
\textrm{for}\,\,\,i=1,\ldots,m,
\,\,\,j=0,1,2.
\endaligned
\end{equation}
Namely, we prove that there exists a positive constant $C$ such that for any
points $\xi=(\xi_1,\ldots,\xi_m)\in\mathcal{O}_d$ and
$h\in C^{0,\alpha}(\overline{\Omega}_\varepsilon)$,
$$
\aligned
\|\phi\|_{L^{\infty}(\Omega_\varepsilon)}\leq C
\|h\|_{*},
\endaligned
$$
for $\varepsilon$ small enough.
By contradiction, assume the existence of  sequences
$\varepsilon_n\rightarrow0$,
points $\xi^n=(\xi_1^n,\ldots,\xi_m^n)\in\mathcal{O}_d$,
parameters $\mu^n=(\mu_1^n,\ldots,\mu_m^n)$
functions $h_n$,  $W_{(\xi^n)'}$  and associated solutions $\phi_n$
such that
\begin{equation}\label{3.16}
\aligned
\|\phi_n\|_{L^{\infty}(\Omega_{\varepsilon_n})}=1
\,\,\,\,\quad\,\,\,\,
\textrm{but}
\,\,\,\,\quad\,\,\,\,
\|h_n\|_{*}\rightarrow0,
\,\,\quad\,\,\textrm{as}\,\,\,\,n\rightarrow+\infty.
\endaligned
\end{equation}
Let us set $\widehat{\phi}^n_i(z)=\phi_n\left(\mu_i^nz+(\xi_i^n)'\right)$
for $i=1,\ldots,m$. By (\ref{3.13}) and  the expansion of $W_{(\xi^n)'}$
in (\ref{3.2}), elliptic regularity theory readily implies that
for each $i=1,\ldots,m$, $\widehat{\phi}^n_i$
converges uniformly over compact subsets of $\mathbb{R}^2$ to a  bounded solution
$\widehat{\phi}^{\infty}_i$ of equation (\ref{3.8})
and hence
$\widehat{\phi}^{\infty}_i$ must be a linear combination of the functions
 $Z_j$, $j=0,1,2$.
Moreover, in view of
 $\|\widehat{\phi}^n_i\|_{\infty}\leq1$,
 the corresponding orthogonality condition  of (\ref{3.15})
over $\widehat{\phi}^n_i$
passes to the limit and  by  Lebesgue's  theorem, it follows that
$$
\aligned
\int_{\mathbb{R}^2}\frac{8}{(1+|z|^2)^2}Z_j(z)\widehat{\phi}_i^{\infty}dy=0
\,\,\,\,\,\quad\,
\textrm{for}\,\,\,j=0,1,2.
\endaligned
$$
So, $\widehat{\phi}_i^{\infty}\equiv0$ for any $i=1,\ldots,m$.
By definition (\ref{3.12}) we conclude
$\lim_{n\rightarrow+\infty}\|\phi_n\|_{**}=0$. But
(\ref{3.13}) and (\ref{3.16}) tell us
$\liminf_{n\rightarrow+\infty}\|\phi_n\|_{**}>0$,
which is a contradiction.

{\bf Step 4:} We prove that for any solution $\phi$ of problem
$\mathcal{L}(\phi)=h$ in $\Omega_\varepsilon$,
$\phi=0$
on $\po_\varepsilon$, which  in addition satisfies the orthogonality conditions:
\begin{equation}\label{3.17}
\aligned
\int_{\Omega_\varepsilon}
e^{
\omega_{\mu_i}\left(y-\xi'_i\right)}
Z_{ij}\phi=0
\,\,\,\,\,\quad\,
\textrm{for}\,\,\,i=1,\ldots,m,
\,\,\,j=1,2,
\endaligned
\end{equation}
there exists a positive constant $C>0$
such that
$$
\aligned
\|\phi\|_{L^{\infty}(\Omega_\varepsilon)}\leq C|\log\varepsilon|
\|h\|_{*},
\endaligned
$$
for $h\in C^{0,\alpha}(\overline{\Omega}_\varepsilon)$.
Proceeding  by contradiction as in Step 3, we
can suppose further that
\begin{equation}\label{3.18}
\aligned
\|\phi_n\|_{L^{\infty}(\Omega_{\varepsilon_n})}=1
\,\,\,\,\quad\,\,\,\,
\textrm{but}
\,\,\,\,\quad\,\,\,\,
|\log\varepsilon_n|\|h_n\|_{*}\rightarrow0
\,\,\ \quad\,\,\textrm{as}\,\,n\rightarrow+\infty.
\endaligned
\end{equation}
but  we loss  the condition
$\int_{\mathbb{R}^2}\frac{8}{(1+|z|^2)^2}Z_0(z)\widehat{\phi}_i^{\infty}=0$
in the limit.
Therefore,
\begin{equation}\label{3.19}
\aligned
\widehat{\phi}_i^{n}\rightarrow
\widehat{\phi}_i^{\infty}=C_i\frac{|z|^2-1}{|z|^2+1}
\,\,\,\,\,
\,\,\textrm{in}\,\,\,C_{loc}^0(\mathbb{R}^2),
\endaligned
\end{equation}
with some constants $C_i$.
To give a contradiction, we have to  show that
$C_i=0$ for all $i=1,\ldots,m$, which will be achieved
by the stronger assumption
of $h_n$ in (\ref{3.18}).

For this aim, we use functions $Z_{i0}$ to build
suitable test functions. Let us consider
radial solutions
$\omega$ and $\theta$ respectively of
$$
\aligned
\Delta\omega+\frac{8}{(1+|z|^2)^2}\omega=
\frac{8}{(1+|z|^2)^2}Z_0(z)
\quad\ \quad\textrm{and}
\quad\ \quad
\Delta\theta+\frac{8}{(1+|z|^2)^2}\theta=
\frac{8}{(1+|z|^2)^2}
\quad
\quad\textrm{in}\,\,\,\mathbb{R}^2,
\endaligned
$$
having asymptotic (see \cite{CI})
$$
\aligned
\omega(z)=
\frac23\log\left(1+|z|^2\right)+O\left(\frac1{1+|z|}
\right),
\quad\quad\quad
\quad
\theta(z)=
O\left(\frac1{1+|z|}
\right)
\ \,\,
\,\quad\textrm{as}\,\,\,
|z|\rightarrow+\infty,
\endaligned
$$
and
$$
\aligned
\nabla\omega(z)
=\frac43\cdot\frac{z}{1+|z|^2}+
O\left(\frac1{1+|z|^2}\right),
\quad\,\,\,\qquad
\nabla\theta(z)=
O\left(\frac1{1+|z|^2}
\right)\,\,
\,\quad\textrm{for all}\,\,\,
z\in\mathbb{R}^2,
\endaligned
$$
because
$$
\aligned
8\int_0^{+\infty}r\frac{(r^2-1)^2}{(r^2+1)^4}dr=\frac43
\qquad\ \qquad
\textrm{and}
\qquad\ \qquad
8\int_0^{+\infty}r\frac{r^2-1}{(r^2+1)^3}dr=0.
\endaligned
$$
Obviously, we can take the following function as an explicit solution
$$
\aligned
\theta(z)=
1-Z_0(z)=\frac{2}{|z|^2+1}
\qquad\quad
\textrm{such that}
\quad
\theta(z)=
O\left(\frac1{1+|z|^2}
\right)
\quad
\textrm{as}\,\,\,
|z|\rightarrow+\infty.
\endaligned
$$
For the sake of simplicity, from now on we omit the dependence on $n$.
For $i=1,\ldots,m$, we define
\begin{equation}\label{3.20}
\aligned
u_i(y)=\omega\left(\frac{y-\xi'_i}{\mu_i}\right)
+\frac{4\log(\varepsilon\mu_i)}{3}Z_{i0}(y)
+\frac{8\pi}{3}H (\xi_i,\xi_i)
\theta\left(\frac{y-\xi'_i}{\mu_i}\right)
\endaligned
\end{equation}
and denote its projection $Pu_i=u_i+\widetilde{H}_i$
on the space $H_{0}^1(\Omega_\varepsilon)$, where
$\widetilde{H}_i$ is a correction term defined as the solution of
$$
\aligned
-\Delta\widetilde{H}_i=0
\quad\,\,
\textrm {in}\,\,\,\,\Omega_\varepsilon,
\quad\quad\quad\quad\quad
\widetilde{H}_i=
-u_i
\quad\,\,
\textrm{on}\,\,\,\,\po_\varepsilon.
\endaligned
$$
Observe that
$$
\aligned
\widetilde{H}_i(y)
=-\frac{4}{3}\log|\varepsilon y-\xi_i|
+O\left(\varepsilon\right)
\qquad
\textrm{in}
\,\,\,\,\,
C^1(\po_\varepsilon).
\endaligned
$$
Similar to the argument in the proof of Lemma 2.1,
we can easily derive  that
\begin{equation}\label{3.21}
\aligned
Pu_i=u_i-\frac{8\pi}{3}H (\varepsilon y,\xi_i)
+O\left(
\varepsilon
\right)
\,\,\,\quad\quad
\textrm{uniformly in}
\,\,\,
C^1(\overline{\Omega}_\varepsilon)\cap
C^{\infty}(\Omega_\varepsilon)
\,\,\,
\textrm{as}\,\,\,
\varepsilon\rightarrow0.
\endaligned
\end{equation}
Moreover, the test function $Pu_i$ solves
\begin{equation}\label{3.22}
\aligned
\Delta
Pu_i+W_{\xi'} Pu_i=e^{
\omega_{\mu_i}\left(y-\xi'_i\right)}Z_{i0}
+\left(
W_{\xi'}-e^{
\omega_{\mu_i}\left(y-\xi'_i\right)}
\right)
Pu_i+E_i
\,\,\quad\,\,
\textrm {in}\,\,\,\,\,\,\Omega_\varepsilon,
\endaligned
\end{equation}
where
\begin{equation}\label{3.23}
\aligned
E_i(y)=\left(Pu_i-u_i
+\frac{8\pi}{3}H(\xi_i,\xi_i)
\right)
e^{
\omega_{\mu_i}\left(y-\xi'_i\right)
}.
\endaligned
\end{equation}
Multiply (\ref{3.22}) by $\phi$ and integrate by parts to obtain
\begin{equation}\label{3.24}
\aligned
\int_{\Omega_\varepsilon}
e^{
\omega_{\mu_i}\left(y-\xi'_i\right)}
Z_{i0}\phi+
\int_{\Omega_\varepsilon}
\left(
W_{\xi'}-e^{
\omega_{\mu_i}\left(y-\xi'_i\right)}
\right)
Pu_i\phi
=-\int_{\Omega_\varepsilon}Pu_ih-
\int_{\Omega_\varepsilon}E_i\phi.
\endaligned
\end{equation}

We analyze each term of (\ref{3.24}).
First of all, by (\ref{3.19}) and  Lebesgue's theorem   we get
\begin{equation}\label{3.25}
\aligned
\int_{\Omega_\varepsilon}
e^{
\omega_{\mu_i}\left(y-\xi'_i\right)}
Z_{i0}\phi
\rightarrow
C_i\int_{\mathbb{R}^2}\frac{8(|z|^2-1)^2}{(|z|^2+1)^4}dz=\frac{8\pi}{3}C_i.
\endaligned
\end{equation}
For the second term in the left-hand side of (\ref{3.24}), we decompose
$$
\aligned
\int_{\Omega_\varepsilon}
\left(
W_{\xi'}-e^{
\omega_{\mu_i}\left(y-\xi'_i\right)}
\right)
Pu_i\phi
=\left[
\sum_{k=1}^m
\left(
\int_{B_{\mu_k|\log\varepsilon|^\tau}(\xi'_k)}
+
\int_{B_{d/\varepsilon}(\xi'_k)\setminus B_{\mu_k|\log\varepsilon|^\tau}(\xi'_k)}
\right)
+\int_{\Omega_\varepsilon\setminus
\bigcup_{k=1}^mB_{d/\varepsilon}(\xi'_k)
}
\right]\left(
W_{\xi'}-e^{
\omega_{\mu_i}\left(y-\xi'_i\right)}
\right)Pu_i\phi.
\endaligned
$$
Using (\ref{3.20}), (\ref{3.21}) and the expansion of $W_{\xi'}$ in (\ref{3.2}), we obtain
$$
\aligned
&
\quad\,\,\,
\int_{B_{\mu_i|\log\varepsilon|^\tau}(\xi'_i)}
\left(
W_{\xi'}-e^{
\omega_{\mu_i}\left(y-\xi'_i\right)}
\right)
Pu_i\phi
\\[1mm]
&
=\,\frac{p-1}{p}\frac1{\gamma^p}\frac{4\log(\varepsilon\mu_i)}{3}
\int_{B_{\mu_i|\log\varepsilon|^\tau}(0)}
\frac{8\mu_i^2}{(\mu_i^2+|z|^2)^2}
Z_0\left(\frac{z}{\mu_i}\right)
\widehat{\phi}_i\left(\frac{z}{\mu_i}\right)
\left[
1+\omega^{1}_{\mu_i}+\frac{1}{2}(\omega_{\mu_i})^2+2\omega_{\mu_i}
\right](z)
dz+O\left(\frac{1}{|\log\varepsilon|}\right)
\\[1mm]
&
=-\frac{(p-1)C_i}{3}
\int_{\mathbb{R}^2}
\frac{8\mu_i^2}{(\mu_i^2+|z|^2)^2}
\left[Z_0\left(\frac{z}{\mu_i}\right)\right
]^2
\left[
1+\omega^{1}_{\mu_i}+\frac{1}{2}(\omega_{\mu_i})^2+2\omega_{\mu_i}
\right](z)dz+o(1),
\endaligned
$$
since  Lebesgue's theorem  and (\ref{3.19}) give
$$
\aligned
&\quad\,\,
\int_{B_{\mu_i|\log\varepsilon|^\tau}(0)}
\frac{8\mu_i^2}{(\mu_i^2+|z|^2)^2}
Z_0\left(\frac{z}{\mu_i}\right)
\widehat{\phi}_i\left(\frac{z}{\mu_i}\right)
\left[
1+\omega^{1}_{\mu_i}+\frac{1}{2}(\omega_{\mu_i})^2+2\omega_{\mu_i}
\right](z)
dz\\[1mm]
&
=C_i\int_{\mathbb{R}^2}
\frac{8\mu_i^2}{(\mu_i^2+|z|^2)^2}
\left[Z_0\left(\frac{z}{\mu_i}\right)\right
]^2
\left[
1+\omega^{1}_{\mu_i}+\frac{1}{2}(\omega_{\mu_i})^2+2\omega_{\mu_i}
\right](z)dz+o\left(1\right).
\endaligned
$$
Thanks to  relation (\ref{3.26}) regarding $\omega^1_{\mu_i}$ stated in Lemma A.2
$$
\aligned
\int_{\mathbb{R}^2}
\frac{8\mu_i^2}{(\mu_i^2+|z|^2)^2}
\left[Z_0\left(\frac{z}{\mu_i}\right)\right
]^2
\left[
1+\omega^{1}_{\mu_i}+\frac{1}{2}(\omega_{\mu_i})^2+2\omega_{\mu_i}
\right](z)dz=8\pi,
\endaligned
$$
we find
$$
\aligned
\int_{B_{\mu_i|\log\varepsilon|^\tau}(\xi'_i)}
\left(
W_{\xi'}-e^{
\omega_{\mu_i}\left(y-\xi'_i\right)}
\right)
Pu_i\phi
=-\frac{8\pi }{3}(p-1)C_i
+o(1).
\endaligned
$$
Notice that (\ref{3.20})-(\ref{3.21}) imply that, as $\varepsilon\rightarrow0$,
\begin{equation*}\label{}
\aligned
Pu_i=O\left(|\log\varepsilon|\right)
\ \
\textrm{in}
\,\,\,
C^1\big(B_{d/\varepsilon}(\xi'_i)\setminus B_{\mu_i|\log\varepsilon|^\tau}(\xi'_i)\big),
\quad\ \,
\textrm{but}
\quad
Pu_i=-\frac{8\pi}{3}G (\varepsilon y,\xi_i)
+O\left(
\varepsilon
\right)
\ \
\textrm{in}
\,\,\,
C^1\big(\overline{\Omega}_\varepsilon\setminus B_{d/\varepsilon}(\xi'_i)\big).
\endaligned
\end{equation*}
Then by the estimate of $W_{\xi'}$ in (\ref{3.5}),
$$
\aligned
\left|
\int_{B_{d/\varepsilon}(\xi'_i)\setminus B_{\mu_i|\log\varepsilon|^\tau}(\xi'_i)}
\left(
W_{\xi'}-e^{
\omega_{\mu_i}\left(y-\xi'_i\right)}
\right)
Pu_i\phi
\right|
\leq
C|\log\varepsilon|
\int_{B_{d/\varepsilon}(\xi'_i)\setminus B_{\mu_i|\log\varepsilon|^\tau}(\xi'_i)}
\frac{8\mu_i^2}{(\mu_i^2+|y-\xi'_i|^2)^2}dy
=o\left(1\right),
\endaligned
$$
and for all $k\neq i$,
$$
\aligned
\left|
\int_{B_{d/\varepsilon}(\xi'_k)\setminus B_{\mu_k|\log\varepsilon|^\tau}(\xi'_k)}
\left(
W_{\xi'}-e^{
\omega_{\mu_i}\left(y-\xi'_i\right)}
\right)
Pu_i\phi
\right|
\leq
C
\int_{B_{d/\varepsilon}(\xi'_k)\setminus B_{\mu_k|\log\varepsilon|^\tau}(\xi'_k)}
\frac{8\mu_k^2}{(\mu_k^2+|y-\xi'_k|^2)^2}dy
=o\left(1\right).
\endaligned
$$
Moreover, by  the expansion of $W_{\xi'}$ in (\ref{3.2}), we have that for all
$k\neq i$,
$$
\aligned
\int_{B_{\mu_k|\log\varepsilon|^\tau}(\xi'_k)}
\left(
W_{\xi'}-e^{
\omega_{\mu_i}\left(y-\xi'_i\right)}
\right)
Pu_i\phi
=-\frac{8\pi}{3}G (\xi_k,\xi_i)
\int_{B_{\mu_k|\log\varepsilon|^\tau}(\xi'_k)}
\frac{8\mu_k^2}{(\mu_k^2+|y-\xi'_k|^2)^2}
\phi
dy+O\left(\frac{1}{|\log\varepsilon|}\right)
=o\left(1\right),
\endaligned
$$
since  Lebesgue's theorem  and (\ref{3.19}) deduce
$$
\aligned
\int_{B_{\mu_k|\log\varepsilon|^\tau}(\xi'_k)}
\frac{8\mu_k^2}{(\mu_k^2+|y-\xi'_k|^2)^2}
\phi(y)
dy
=\int_{B_{|\log\varepsilon|^\tau}(0)}
\frac{8}{(1+|z|^2)^2}
\widehat{\phi}_k(z)
dz\rightarrow
C_k\int_{\mathbb{R}^2}
\frac{8}{(1+|z|^2)^2}
\frac{|z|^2-1}{|z|^2+1}dz
=0.
\endaligned
$$
In addition, by the estimate of $W_{\xi'}$ in (\ref{3.6}) we get
\begin{eqnarray*}
\int_{
\Omega_\varepsilon\setminus
\bigcup_{k=1}^mB_{d/\varepsilon}(\xi'_k)
}
\left(
W_{\xi'}-e^{
\omega_{\mu_i}\left(y-\xi'_i\right)}
\right)
Pu_i\phi
\qquad\qquad\qquad\qquad\qquad\qquad\
\qquad\qquad\qquad\qquad\qquad\qquad\
\\[1mm]
=
\int_{
\Omega_\varepsilon\setminus
\bigcup_{k=1}^mB_{d/\varepsilon}(\xi'_k)
}
\left
\{
\left(\frac{O\left(\varepsilon^{(2+p)/p}\right)}
{|\log\varepsilon|^{p-1}}
+\frac{O\left(\varepsilon^{(2+p)/p}\right)}
{|\log\varepsilon|^{2(p-1)}}\right)
\exp\left[
-\frac{2-p}{p}|\log\varepsilon|+O\left(\frac{
1
}{|\log\varepsilon|^{p-1}}\right)\right]
+O\left(\varepsilon^4\right)
\right\}dy
=o\left(1\right).
\end{eqnarray*}
Hence we obtain
\begin{equation}\label{3.27}
\aligned
\int_{\Omega_\varepsilon}
\left(
W_{\xi'}-e^{
\omega_{\mu_i}\left(y-\xi'_i\right)}
\right)
Pu_i\phi
=-\frac{8\pi }{3}(p-1)C_i
+o(1).
\endaligned
\end{equation}
As for the right-hand side of (\ref{3.24}), we have that by
(\ref{3.14})  and (\ref{3.21}),
\begin{eqnarray}\label{3.28}
&&\left|\int_{\Omega_\varepsilon}Pu_ih\right|=O\left(\|h\|_{*}
\int_{\Omega_\varepsilon}\left[
\sum\limits_{k=1}^m\frac{\mu_k^\sigma}{(\mu_k^2+|y-\xi'_k|^2)^{(2+\sigma)/2}}
\right]|u_i|dy
\right)+O\left(
\|h\|_{*}
\right)
=O
\left(
|\log\varepsilon|
\|h\|_{*}
\right),
\end{eqnarray}
since $|u_i|=O\left(|\log\varepsilon|\right)$ in $\mathbb{R}^2$ and
$$
\aligned
\int_{B_{\mu_k|\log\varepsilon|^\tau}(\xi'_k)}
\frac{\mu_k^\sigma}{(\mu_k^2+|y-\xi'_k|^2)^{(2+\sigma)/2}}|u_i|dy
\leq
\int_{\mathbb{R}^2}
\frac{1}{(1+|z|^2)^{(2+\sigma)/2}}|u_i(\xi'_k+\mu_k z)|dz
=O\left(|\log\varepsilon|\right).
\endaligned
$$
By (\ref{3.21}) and (\ref{3.23})
we deduce
\begin{equation}\label{3.30}
\aligned
\int_{\Omega_\varepsilon}E_i\phi
=O\left(\varepsilon
\int_{\Omega_\varepsilon}
e^{
\omega_{\mu_i}\left(y-\xi'_i\right)
}
\big(| y-\xi'_i|+1\big)dy\right)
=O\left(\varepsilon\right).
\endaligned
\end{equation}
Finally, substituting  (\ref{3.25})-(\ref{3.30}) into (\ref{3.24}) and taking into account
the assumption condition (\ref{3.18}), we conclude
$$
\aligned
\frac{8\pi }{3}(2-p)C_i=o(1)
\,\quad\,\quad\textrm{for any}\,\,\,i=1,\ldots,m.
\endaligned
$$
Necessarily, $C_i=0$ by contradiction and
the claim is proved.

{\bf Step 5:} We establish the validity of the a priori estimate
\begin{equation}\label{3.31}
\aligned
\|\phi\|_{L^{\infty}(\Omega_\varepsilon)}
\leq C|\log\varepsilon|
\|h\|_{*}
\endaligned
\end{equation}
for solutions of
problem (\ref{3.9}) and $h\in C^{0,\alpha}(\overline{\Omega}_\varepsilon)$.
Step 4 gives
$$
\aligned
\|\phi\|_{L^{\infty}(\Omega_\varepsilon)}\leq C|\log\varepsilon|
\left(
\|h\|_{*}+\sum_{i=1}^m\sum_{j=1}^2|c_{ij}|\|e^{\omega_{\mu_i}\left(y-\xi'_i\right)}Z_{ij}\|_*
\right)\leq C|\log\varepsilon|\left(
\|h\|_{*}+\sum_{i=1}^m\sum_{j=1}^2|c_{ij}|
\right).
\endaligned
$$
As before, proceeding by contradiction as in Step $3$, we
can  suppose further that
\begin{equation}\label{3.32}
\aligned
\|\phi_n\|_{L^{\infty}(\Omega_{\varepsilon_n})}=1,
\,\,\ \,\quad\,\,|\log\varepsilon_n|\|h_n\|_{*}\rightarrow0,
\,\,\ \,\quad\,\,
|\log\varepsilon_n|\sum_{i=1}^m\sum_{j=1}^2|c_{ij}^n|\geq\delta>0
\,\,\,\quad\,\,\,
\textrm{as}\,\,\,n\rightarrow+\infty.
\endaligned
\end{equation}
We omit the dependence on $n$.
It suffices to estimate the size of the coefficients $c_{ij}$.
To this end, we first define  $PZ_{ij}$ as the projection
on $H_{0}^1(\Omega_\varepsilon)$ of $Z_{ij}$, precisely
\begin{equation}\label{3.33}
\aligned
-\Delta PZ_{ij}=-\Delta Z_{ij}=
e^{\omega_{\mu_i}\left(y-\xi'_i\right)}
Z_{ij}
\quad
\textrm{in}\,\,\,\,\Omega_\varepsilon,
\qquad\qquad
PZ_{ij}=0
\quad
\textrm{on}\,\,\,\partial\Omega_\varepsilon.
\endaligned
\end{equation}
Then for   $i=1,\ldots,m$ and $j=1,2$,
\begin{equation}\label{3.34}
\aligned
PZ_{ij}=Z_{ij}+8\pi\varepsilon\mu_i
\partial_{(\xi_i)_j}H(\varepsilon y,\xi_i)
+O\left(\varepsilon^3\right)
\qquad\quad
\textrm{in}
\,\,\,
C^1(\overline{\Omega}_\varepsilon),
\endaligned
\end{equation}
and
\begin{equation}\label{3.35}
\aligned
PZ_{ij}=8\pi\varepsilon\mu_i\partial_{(\xi_i)_j}G(\varepsilon y,\xi_i)
+O\left(\varepsilon^3
\right)
\qquad\quad
\textrm{in}
\,\,\,
C^1(\overline{\Omega}_\varepsilon\setminus B_{d/\varepsilon}(\xi'_i)).
\endaligned
\end{equation}
Let us claim that the following ``orthogonality'' relations hold: for any
$i,k=1,\ldots,m$ and $j,l=1,2$,
\begin{equation}\label{3.36}
\aligned
\int_{\Omega_\varepsilon}e^{\omega_{\mu_i}\left(y-\xi'_i\right)}
Z_{ij}PZ_{kl}=\left(64\int_{\mathbb{R}^2}
\frac{|z|^2}{(1+|z|^2)^4}
\right)\delta_{ik}\delta_{jl}
+O\left(\varepsilon^2\right),
\endaligned
\end{equation}
uniformly for $\xi=(\xi_1,\ldots,\xi_m)\in\mathcal{O}_d$,
where $\delta_{ik}$ and $\delta_{jl}$
denote the Kronecker's symbols.
Indeed,
by (\ref{3.34})-(\ref{3.35})
we get
$$
\aligned
\int_{\Omega_\varepsilon}e^{\omega_{\mu_i}\left(y-\xi'_i\right)}Z_{ij}PZ_{il}=&
\int_{B_{d/\varepsilon}(\xi'_i)}
e^{\omega_{\mu_i}\left(y-\xi'_i\right)}Z_{ij}
\Big[
Z_{il}+8\pi\varepsilon\mu_i
\partial_{(\xi_i)_l}H(\xi_i,\xi_i)
+O\left(\varepsilon^2|y-\xi'_i|
+\varepsilon^3\right)
\Big
]dy
+O\left(\varepsilon^4\right)\\
=&\int_{B_{d/(\varepsilon\mu_i)}(0)}
\frac{8}{(1+|z|^2)^2}\frac{4z_j}{1+|z|^2}
\left[
\frac{4z_l}{1+|z|^2}
+O\left(\varepsilon^2|z|\right)
\right]dz+O\left(\varepsilon^3\right)
\\[1mm]
=&\left(64\int_{\mathbb{R}^2}
\frac{|z|^2}{(1+|z|^2)^4}
\right)\delta_{jl}
+O\left(\varepsilon^2\right),
\endaligned
$$
but for $i\neq k$,
$$
\aligned
\int_{\Omega_\varepsilon}
e^{\omega_{\mu_i}\left(y-\xi'_i\right)}
Z_{ij}
P Z_{kl}=
\int_{B_{d/\varepsilon}(\xi'_i)}
e^{\omega_{\mu_i}\left(y-\xi'_i\right)}Z_{ij}
\Big[
8\pi\varepsilon\mu_k\partial_{(\xi_k)_l}G(\xi_i,\xi_k)
+O\left(\varepsilon^2|y-\xi'_i|
+\varepsilon^3
\right)
\Big]dy
+O\left(\varepsilon^3\right)
=O\left(\varepsilon^2\right).
\endaligned
$$
Next,
testing (\ref{3.9}) against
$PZ_{ij}$, $i=1,\ldots,m$ and $j=1,2$,  we obtain
\begin{equation}\label{3.38}
\aligned
\sum_{k=1}^m\sum_{l=1}^2c_{kl}\int_{\Omega_\varepsilon}
e^{\omega_{\mu_k}\left(y-\xi'_k\right)}
Z_{kl}
PZ_{ij}+\int_{\Omega_\varepsilon}hPZ_{ij}=
\int_{\Omega_\varepsilon}e^{\omega_{\mu_i}\left(y-\xi'_i\right)}Z_{ij}\phi
-\int_{\Omega_\varepsilon}W_{\xi'}\phi PZ_{ij}.
\endaligned
\end{equation}
From (\ref{3.14}) and (\ref{3.34}) we have
$$
\aligned
\left|\int_{\Omega}hPZ_{ij}\right|
=O\left(\|h\|_*\right),
\endaligned
$$
then  by (\ref{3.36}),
the left-hand side of (\ref{3.38}) takes the form
\begin{equation}\label{3.39}
\aligned
\textrm{L.H.S of (\ref{3.38})}=
Dc_{ij}+O\left(\varepsilon^2
\sum_{k=1}^m\sum_{l=1}^2 |c_{kl}|\right)+
O\left(\|h\|_*\right)
\endaligned
\end{equation}
with $D=64\int_{\mathbb{R}^2}
\frac{|z|^2}{(1+|z|^2)^4}$.
On the other hand, by   the estimates of $W_{\xi'}$
in (\ref{3.1})-(\ref{3.3})
the right-hand side of (\ref{3.38}) becomes
\begin{eqnarray}\label{3.40}
\textrm{R.H.S of (\ref{3.38})}
=
\int_{\Omega_\varepsilon}e^{\omega_{\mu_i}\left(y-\xi'_i\right)}Z_{ij}\phi
-\int_{B_{\mu_i|\log\varepsilon|^\tau}(\xi'_i)}W_{\xi'}\phi P Z_{ij}
+O\left(\frac{\|\phi\|_{L^\infty(\Omega_\varepsilon)}}{|\log\varepsilon|^2}\right)
\qquad\qquad\qquad\qquad\qquad\qquad\quad\,
\nonumber\\
=
-\int_{B_{\mu_i|\log\varepsilon|^\tau}(\xi'_i)}
\left(W_{\xi}-
e^{\omega_{\mu_i}\left(y-\xi'_i\right)}
\right)
\phi P Z_{ij}
+
\int_{\Omega_\varepsilon}e^{\omega_{\mu_i}\left(y-\xi'_i\right)}\big(Z_{ij}-P Z_{ij}\big)\phi
+O\left(\frac{\|\phi\|_{L^\infty(\Omega_\varepsilon)}}{|\log\varepsilon|^2}\right)
\quad\,\,
\nonumber\\
=-\frac{p-1}{p}\frac1{\gamma^p}
\int_{
B_{|\log\varepsilon|^\tau}(0)
}
\frac{32z_j}{(1+|z|^2)^{3}}\widehat{\phi}_i(z)
\left[1+\omega^{1}_{\mu_i}+\frac{1}{2}(\omega_{\mu_i})^2+2\omega_{\mu_i}\right]\left(\mu_i|z|\right)dz
+O\left(\frac{\|\phi\|_{L^\infty(\Omega_\varepsilon)}}{|\log\varepsilon|^2}\right)
\end{eqnarray}
because of (\ref{3.34})-(\ref{3.35}), where $\widehat{\phi}_i(z)=\phi(\mu_iz+\xi'_i)$.
Substituting  estimates (\ref{3.39})-(\ref{3.40}) into  (\ref{3.38}), we obtain
$$
\aligned
Dc_{ij}+O\left(\varepsilon^2
\sum_{k=1}^m\sum_{l=1}^2 |c_{kl}|\right)
=O\left(\|h\|_{*}+\frac{1}{|\log\varepsilon|}\|\phi\|_{L^\infty(\Omega_\varepsilon)}\right).
\endaligned
$$
Then
\begin{equation}\label{3.28}
\aligned
\sum_{k=1}^m\sum_{l=1}^2|c_{kl}|
=O\left(\frac{1}{|\log\varepsilon|}\|\phi\|_{L^\infty(\Omega_\varepsilon)}\right)+
O\left(\|h\|_{*}\right).
\endaligned
\end{equation}
From the first two assumptions in (\ref{3.32}) we give
$\sum_{k=1}^m\sum_{l=1}^2|c_{kl}|=o\left(1\right)$.
As in Step 4, we conclude  that for each $i=1,\ldots,m$,
$$
\aligned
\widehat{\phi}_i\rightarrow
C_i\frac{|z|^2-1}{|z|^2+1}
\,\,\,\,\,
\,\,\textrm{uniformly in}\,\,\,C_{loc}^0(\mathbb{R}^2),
\endaligned
$$
with some constants $C_i$.
In view of the oddness of the function  $\frac{32z_j}
{(1+|z|^{2})^3}$, $j=1,2$, by (\ref{8.2}) and
 Lebesgue's theorem we find
$$
\aligned
\int_{
B_{|\log\varepsilon|^\tau}(0)
}
\frac{32z_j}{(1+|z|^2)^{3}}\widehat{\phi}_i(z)
\left[1+\omega^{1}_{\mu_i}+\frac{1}{2}(\omega_{\mu_i})^2+2\omega_{\mu_i}\right]\left(\mu_i|z|\right)dz
\rightarrow
0.
\endaligned
$$
Substituting estimates (\ref{3.39})-(\ref{3.40}) into  (\ref{3.38}) again,  we have a better estimate
$$
\aligned
\sum_{k=1}^m\sum_{l=1}^2|c_{kl}|
=o\left(\frac{1}{|\log\varepsilon|}\right)+
O\left(\|h\|_{*}\right),
\endaligned
$$
which is impossible because of the last assumption in (\ref{3.32}).

{\bf Step 6:}  We prove the solvability  of
problem (\ref{3.9}). For this purpose, we introduce the subspace
of $ H^1_0(\Omega_\varepsilon)$ defined by
$$
\aligned
K_{\xi'}=\left\{\sum_{i=1}^m\sum_{j=1}^2 c_{ij} PZ_{ij}:
\,\,c_{ij}\in\mathbb{R}\,\,\,\,\textrm{for}\,\,\,
i=1,\ldots,m,
\,\,
j=1,2
\right\},
\endaligned
$$
and its orthogonal space
$$
\aligned
K_{\xi'}^{\bot}=\left\{\phi\in H_0^1(\Omega_\varepsilon):\,\,
\int_{\Omega_\varepsilon}e^{\omega_{\mu_i}\left(y-\xi'_i\right)}
Z_{ij}\phi=0
\,\,\,\,\textrm{for}\,\,\,
i=1,\ldots,m,
\,\,
j=1,2
\right\}.
\endaligned
$$
Let $\Pi_{\xi'}:  H^1_0(\Omega_\varepsilon)\rightarrow K_{\xi'}$
and $\Pi_{\xi'}^\bot=Id-\Pi_{\xi'}:  H^1_0(\Omega_\varepsilon)\rightarrow K^\bot_{\xi'}$
be
the corresponding
orthogonal projections such that
$$
\aligned
\Pi_{\xi'}\phi=\sum_{i=1}^m\sum_{j=1}^2 c_{ij} PZ_{ij},
\endaligned
$$
where the coefficients $c_{ij}$ are uniquely determined  by the system
$$
\aligned
\int_{\Omega_\varepsilon}e^{\omega_{\mu_i}\left(y-\xi'_k\right)}
Z_{kl}
\left(
\phi-
\sum_{i=1}^m\sum_{j=1}^2 c_{ij} PZ_{ij}
\right)=0\,\,\quad\,\,\textrm{for any}\,\,\,
k=1,\ldots,m,
\,\,\,l=1,2
\endaligned
$$
because of  (\ref{3.36}).
Problem (\ref{3.9}), expressed in a weak form, is equivalent to finding
$\phi\in K_{\xi}^{\bot}$ such that
$$
\aligned
\langle\phi,\,\psi\rangle_{H^1_0(\Omega_\varepsilon)}=\int_{\Omega_\varepsilon}\left(W_{\xi'}\phi
+h\right)\psi\,\,\,\ \quad\,\,\,\textrm{for all}
\,\,\,\psi\in K_{\xi'}^{\bot}.
\endaligned
$$
With the aid of Riesz's representation theorem,
this equation  can be rewritten in $K_{\xi}^{\bot}$ in the operator form
$$
\aligned
(Id-K)\phi=\tilde{h},
\endaligned
$$
where
$\tilde{h}=\Pi_{\xi'}^\bot(-\Delta)^{-1}h$
and
$K(\phi)=\Pi_{\xi'}^\bot(-\Delta)^{-1}(W_{\xi'}\phi)$ is a linear compact operator
in $K_{\xi'}^{\bot}$.
Fredholm's alternative
guarantees    unique solvability  of
this problem
for any $\tilde{h}\in K^\bot_\xi$ provided
that the homogeneous
equation $\phi=K(\phi)$ has only the trivial solution
in $K_{\xi'}^{\bot}$, which in turn follows from
the a priori estimate (\ref{3.31}) in Step 5.
Finally,  by elliptic regularity theory
the solution constructed in this way
belongs to $H^2(\Omega_\varepsilon)$.
Moreover,
by density of $C^{0,\alpha}(\overline{\Omega}_\varepsilon)$
in $(C(\overline{\Omega}_\varepsilon),\,\|\cdot\|_{L^{\infty}(\Omega_\varepsilon)})$,
we can approximate $h\in C(\overline{\Omega}_\varepsilon)$
by smooth functions and, by (\ref{3.31}) and
elliptic regularity theory, we obtain the validity of
(\ref{3.31}) also  for  $h\in C(\overline{\Omega}_\varepsilon)$
(not only for $h\in C^{0,\alpha}(\overline{\Omega}_\varepsilon)$).
\qquad\qquad\qquad\qquad\qquad\qquad\qquad\qquad\qquad\qquad\qquad
\qquad\qquad\qquad\qquad\qquad\qquad\qquad$\square$

\vspace{1mm}
\vspace{1mm}
\vspace{1mm}
\vspace{1mm}

\noindent{\bf Remark 3.3.}\,\,The operator $\mathcal{T}$ is differentiable  with
respect to the variables $\xi'=(\xi_1',\ldots,\xi_m')$.
Indeed, similar to
those used in \cite{DKM},  if  we fix $h\in C(\overline{\Omega}_\varepsilon)$ with
$\|h\|_*<\infty$  and set $\phi=T(h)$, then
by formally
computing the derivative of $\phi$ with respect to $\xi'=(\xi'_1,\ldots,\xi'_m)$
and using the delicate estimate $\|\partial_{(\xi'_{i})_j}W_{\xi'}\|_{*}=O\left(1\right)$
we obtain the a priori  estimate
\begin{equation*}\label{3.74}
\aligned
\|\partial_{(\xi'_{i})_j}\mathcal{T}(h)\|_{L^{\infty}(\Omega_\varepsilon)}\leq C|\log\varepsilon|^2\,\|h\|_{*},
\quad\,\,\,\,\,\forall\,\,\,
i=1,\ldots,m,\,\,\,j=1,2.
\endaligned
\end{equation*}

\vspace{1mm}
\vspace{1mm}
\vspace{1mm}
\vspace{1mm}

\noindent{\bf Remark 3.4.}\,\,Given $h\in C(\overline{\Omega}_\varepsilon)$ with
$\|h\|_*<\infty$, let $\phi$ be the unique solution of problem (\ref{3.9}) given by Proposition 3.2.
Multiplying  (\ref{3.9}) by $\phi$ and integrating by parts, we get
$$
\aligned
\|\phi\|_{H^1_0(\Omega_\varepsilon)}^2
=\int_{\Omega_\varepsilon}W_{\xi'}\phi^2
+\int_{\Omega_\varepsilon}h\phi.
\endaligned
$$
By Proposition 3.1 we find
\begin{equation*}\label{3.73}
\aligned
\|\phi\|_{H_0^1(\Omega_\varepsilon)}\leq C\big(\|\phi\|_{L^{\infty}(\Omega_\varepsilon)}+\|h\|_*\big).
\endaligned
\end{equation*}

\vspace{1mm}
\vspace{1mm}
\vspace{1mm}

\section{The nonlinear projected problem}
Consider the nonlinear  projected  problem: for any points
$\xi=(\xi_1,\ldots,\xi_m)\in\mathcal{O}_d$,
we find a function $\phi$ such that
\begin{equation}\label{4.1}
\left\{\aligned
&-\Delta
(V_{\xi'}+\phi) =f(V_{\xi'}+\phi)+
\sum\limits_{i=1}^m\sum\limits_{j=1}^{2}c_{ij}e^{
\omega_{\mu_i}\left(y-\xi'_i\right)}Z_{ij}
\,\,\ \,\textrm{in}\,\,\,\,\,\,\Omega_\varepsilon,\\
&
\phi=0\,\,
\qquad\qquad\qquad
\qquad\qquad\qquad\qquad
\qquad\qquad\qquad\quad
\textrm{on}\,\,\,\,\partial\Omega_{\varepsilon},\\[1mm]
&\int_{\Omega_\varepsilon}
e^{
\omega_{\mu_i}\left(y-\xi'_i\right)}
Z_{ij}\phi=0
\qquad\qquad\qquad\,\,\,
\forall\,\,\,i=1,\ldots,m,\,\,\,j=1, 2,
\endaligned\right.
\end{equation}
for some  coefficients
$c_{ij}$, $i=1,\ldots,m$ and $j=1,2$, where
the function
$f(\cdot)$ is given by
(\ref{2.27}). The following result holds.


\vspace{1mm}
\vspace{1mm}
\vspace{1mm}
\vspace{1mm}

\noindent{\bf Proposition 4.1.}\,\,{\it
There exist constants $C>0$ and  $\varepsilon_0>0$ such
that for any $0<\varepsilon<\varepsilon_0$ and any points
$\xi=(\xi_1,\ldots,\xi_m)\in\mathcal{O}_d$,
problem {\upshape(\ref{4.1})} admits
a unique solution
$\phi_{\xi'}$
for some coefficients $c_{ij}(\xi')$,
$i=1,\ldots,m$, $j=1,2$, such that
\begin{equation}\label{4.2}
\aligned
\|\phi_{\xi'}\|_{L^{\infty}(\Omega_\varepsilon)}\leq\frac{C}{|\log\varepsilon|^3},\,\,\quad\,\quad\,\,
\sum_{i=1}^m\sum_{j=1}^{2}|c_{ij}(\xi')|\leq\frac{C}{|\log\varepsilon|^4}\quad\,\quad\textrm{and}\quad\,\quad
\|\phi_{\xi'}\|_{H^1_0(\Omega_\varepsilon)}\leq\frac{C}{|\log\varepsilon|^3}.
\endaligned
\end{equation}
Furthermore, the map $\xi'\mapsto\phi_{\xi'}$ is a $C^1$-function in $C(\overline{\Omega}_\varepsilon)$ and $H^1_0(\Omega_\varepsilon)$,
precisely for any $i=1,\ldots,m$ and $j=1,2$,
\begin{equation}\label{4.3}
\aligned
\|\partial_{(\xi'_{i})_j}\phi_{\xi'}\|_{L^{\infty}(\Omega_\varepsilon)}\leq\frac{C}{|\log\varepsilon|^2},
\endaligned
\end{equation}
where $\xi':=(\xi'_1,\ldots,\xi'_m)=(\frac1{\varepsilon}\xi_1,\ldots,\frac1{\varepsilon}\xi_m)$.
}

\vspace{1mm}

\begin{proof}
Proposition $3.2$ and Remarks 3.3-3.4 allow us to apply the Contraction Mapping
Theorem
and the Implicit Function Theorem
to find a unique solution for problem (\ref{4.1})
satisfying (\ref{4.2})-(\ref{4.3}).
Since it is a standard procedure, we omit the details,
see Lemmas 4.1-4.2 in \cite{DKM} for a similar proof.
We just mention that
$\|N(\phi)\|_{*}\leq
C\|\phi\|_{L^{\infty}(\Omega_\varepsilon)}^2$
and
$\|\partial_{(\xi'_{i})_j}E_{\xi'}\|_{*}
\leq
C|\log\varepsilon|^{-3}$.
\end{proof}

\section{Variational reduction}
Since problem (\ref{4.1}) has been solved, we  find a solution of problem (\ref{2.29})
and hence to the original equation  (\ref{1.1}) if we detect some points  $\xi'$ such  that
\begin{equation}\label{5.1}
\aligned
c_{ij}(\xi')=0\,\quad\,\,\textrm{for all}\,\,\,i=1,\ldots,m,\,\,\,j=1,2.
\endaligned
\end{equation}
Recall that
$\lambda$ is assumed to be a  free  parameter. Let us consider the free  functional $J^p_\lambda$
associated to equation (\ref{1.1}), namely
\begin{equation}\label{5.10}
\aligned
J^p_\lambda(u)=\frac12\int_{\Omega} |\nabla u|^2-\frac{\lambda}{p}
\int_{\Omega}e^{|u|^p},\,\,\ \ \,\,\forall\,
u\in H_0^1(\Omega).
\endaligned
\end{equation}
Furthermore, we
take its finite-dimensional restriction
\begin{equation}\label{5.2}
\aligned
F_\lambda(\xi)= J^p_\lambda\left(\big(U_\xi+\widetilde{\phi}_{\xi}\big)(x)
\right)\,\ \,\,\ \ \,\forall\,\,\xi=(\xi_1,\ldots,\xi_m)\in\mathcal{O}_d,
\endaligned
\end{equation}
where
\begin{equation}\label{5.3}
\aligned
\big(U_\xi+\widetilde{\phi}_{\xi}\big)(x)=\gamma+\frac1{p\gamma^{p-1}}\big(V_{\xi'}+\phi_{\xi'}\big)
\left(\frac{x}{\varepsilon}\right),
\,\quad\,x\in\Omega,
\endaligned
\end{equation}
with $V_{\xi'}$ defined in (\ref{2.28}) and $\phi_{\xi'}$ the unique solution to problem
(\ref{4.1}) given by Proposition 4.1.
Let
\begin{equation}\label{5.4}
\aligned
I_\varepsilon(\upsilon)=\frac12\int_{\Omega_\varepsilon}
|\nabla \upsilon|^2dy
-\int_{\Omega_\varepsilon}e^{\gamma^p\left[\left|1+\frac{\upsilon}{p\gamma^p}\right|^p-1\right]}dy,
\,\quad\,
\forall
\,\upsilon\in H^1_0(\Omega_\varepsilon).
\endaligned
\end{equation}
By (\ref{1.5}),
\begin{equation}\label{5.5}
\aligned
I_\varepsilon\big(V_{\xi'}+\phi_{\xi'}\big)=
p^2\gamma^{2(p-1)}F_\lambda(\xi)
\qquad\,\,
\textrm{and}
\qquad\,\,
I_\varepsilon\big(V_{\xi'}+\phi_{\xi'}\big)-I_\varepsilon(V_{\xi'}
)=
p^2\gamma^{2(p-1)}\big[F_\lambda(\xi)
-
J_\lambda\left(U_\xi\right)\big].
\endaligned
\end{equation}

\vspace{1mm}
\vspace{1mm}

\noindent{\bf Proposition 5.1.}\,\,{\it The function $F_\lambda:\mathcal{O}_d\mapsto\mathbb{R}$
is of class $C^1$.
Moreover, for all $\lambda$ sufficiently small,
if  $D_{\xi}F_\lambda(\xi)=0$,  then $\xi'=\xi/\varepsilon$ satisfies {\upshape (\ref{5.1})},
that is,
$U_\xi+\widetilde{\phi}_{\xi}$ is a solution of equation {\upshape(\ref{1.1})}.
}

\vspace{1mm}

\begin{proof}
The function $F_\lambda$ is of class $C^1$
since the map
$\xi'\mapsto\phi_{\xi'}$ is a $C^1$-function in $C(\overline{\Omega}_\varepsilon)$ and $H^1_0(\Omega_\varepsilon)$.
Assume that  $D_{\xi}F_\lambda(\xi)=0$.
Since $\phi_{\xi'}$
solves problem (\ref{4.1}),
by (\ref{5.5}) we deduce that for any $k=1,\ldots,m$ and $l=1,2$,
\begin{eqnarray}\label{5.6}
&&0
=I'_\varepsilon\big(V_{\xi'}+\phi_{\xi'}\big)\partial_{(\xi'_k)_l}\big(V_{\xi'}+\phi_{\xi'}\big)\nonumber\\[1mm]
&&\,\,\,\,\,=\sum\limits_{i=1}^m\sum\limits_{j=1}^{2}c_{ij}(\xi')
\int_{\Omega_\varepsilon}e^{
\omega_{\mu_i}\left(y-\xi'_i\right)}
Z_{ij}
\partial_{(\xi'_k)_l}V_{\xi'}
-\sum\limits_{i=1}^m\sum\limits_{j=1}^{2}c_{ij}(\xi')\int_{\Omega_\varepsilon}\phi_{\xi'}
\partial_{(\xi'_k)_l}
\left(
e^{
\omega_{\mu_i}\left(y-\xi'_i\right)}
Z_{ij}
\right)
\end{eqnarray}
because of $\int_{\Omega_\varepsilon}e^{
\omega_{\mu_i}\left(y-\xi'_i\right)}
Z_{ij}\phi_{\xi'}=0$.
Recalling  that  $D_{\xi'}V_{\xi'}(y)=p\gamma^{p-1}D_{\xi'}U_\xi(\varepsilon y)$,
by the expression of $U_\xi$ in (\ref{2.16})
we obtain
\begin{equation}\label{5.32}
\aligned
\partial_{(\xi'_k)_l}V_{\xi'}(y)
=\sum_{i=1}^{m}a_i
\partial_{(\xi'_k)_l}\left[
\omega_{\mu_i}\left(y-\xi'_i\right)
+\sum_{j=1}^3\left(\frac{p-1}{p}\right)^j\frac{1}{\gamma^{jp}}
\omega^j_{\mu_{i}}
\left(y-\xi'_i\right)
+p\gamma^{p-1}H_i(\varepsilon y)\right].
\endaligned
\end{equation}
By (\ref{2.3}) and (\ref{3.10}),
\begin{equation}\label{5.29}
\aligned
\partial_{(\xi'_k)_l}\omega_{\mu_i}\left(y-\xi'_i\right)=
\frac1{\mu_i}\delta_{ki}
Z_{l}\left(
\frac{y-\xi'_i}{\mu_i}
\right)
+\left(\partial_{(\xi'_k)_l}\mu_i\right)\frac{d}{d\mu_i}\omega_{\mu_i}\left(y-\xi'_i\right),
\endaligned
\end{equation}
and for  $j=1, 2, 3$, by (\ref{2.13}) and (\ref{3.10}),
\begin{equation}\label{5.30}
\aligned
\partial_{(\xi'_k)_l}\omega^j_{\mu_{i}}
\left(y-\xi'_i\right)=
-\frac1{\mu_i}\delta_{ki}\left[\frac{D^j_{\mu_i}}{4}Z_{l}\left(
\frac{y-\xi'_i}{\mu_i}
\right)+
O\left(\frac{\mu_i^2}{|y-\xi'_i|^2+\mu_i^2}\right)
\right]
+
\left(\partial_{(\xi'_k)_l}\mu_i\right)\frac{d}{d\mu_i}\omega^j_{\mu_i}\left(y-\xi'_i\right),
\endaligned
\end{equation}
where $\delta_{ki}$ denotes the Kronecker's symbol.
In addition,
differentiating  equation (\ref{2.17}) of $H_i$'s with respect to the variable $(\xi'_k)_l$
and using harmonicity and the maximum principle, we can prove
\begin{equation}\label{5.31}
\aligned
\partial_{(\xi'_k)_l}\big[p\gamma^{p-1}H_i(\varepsilon y)\big]
=O\left(\varepsilon\right).
\endaligned
\end{equation}
Inserting (\ref{5.29})-(\ref{5.31}) into (\ref{5.32}) and
using the fact that
$|\partial_{(\xi'_k)_l}\mu_i|=O\left(\varepsilon\right)$
for any $i=1,\ldots,m$,
we have that
\begin{eqnarray}\label{5.7}
\partial_{(\xi'_k)_l}V_{\xi'}(y)=
\frac{a_k}{\mu_k}
\left[
1-\frac{1}{4}\sum_{j=1}^3\left(\frac{p-1}{p}\right)^j\frac{D^j_{\mu_k}}{\gamma^{jp}}
\right]
Z_{l}\left(
\frac{y-\xi'_k}{\mu_k}
\right)
+\sum_{j=1}^3\left(\frac{p-1}{p}\right)^j\frac{1}{\gamma^{jp}}
O\left(\frac{\mu_k}{|y-\xi'_k|^2+\mu_k^2}\right)
+
O\left(\varepsilon\right).
\end{eqnarray}
Notice that  by (\ref{2.3}) and (\ref{3.10}) for $j=1,2$,
\begin{eqnarray}\label{5.8}
\partial_{(\xi'_k)_l}
\left(
e^{
\omega_{\mu_i}\left(y-\xi'_i\right)}
Z_{ij}
\right)
=-4\mu_ie^{
\omega_{\mu_i}\left(y-\xi'_i\right)}
\left\{
\left[
\frac{\delta_{jl}}{|y-\xi'_i|^2+\mu_i^2}-
6\frac{(y-\xi'_i)_j(y-\xi'_k)_l}{(|y-\xi'_i|^2+\mu_i^2)^2}
\right]\delta_{ki}
+O\left(\varepsilon\right)
\right\}.
\end{eqnarray}
Therefore by (\ref{5.7}),  equations
(\ref{5.6}) can be rewritten as, for
each $k=1,\ldots,m$  and  $l=1,2$,
$$
\aligned
\frac{a_k}{\mu_k}\sum\limits_{i=1}^m\sum\limits_{j=1}^{2}
c_{ij}(\xi')
\int_{\Omega_\varepsilon}e^{
\omega_{\mu_i}\left(y-\xi'_i\right)}
Z_{ij}Z_{kl}+
\sum\limits_{i=1}^m\sum\limits_{j=1}^{2}|c_{ij}(\xi')|O\left(\frac{1}{|\log\varepsilon|}
+\|\phi_{\xi'}\|_{L^{\infty}(\Omega_\varepsilon)}\int_{\Omega_\varepsilon}
\left|\partial_{(\xi'_k)_l}
\left(
e^{
\omega_{\mu_i}\left(y-\xi'_i\right)}
Z_{ij}
\right)
\right|
\right)
=0,
\endaligned
$$
so that,  using (\ref{4.2}), (\ref{5.8}) and the argument in expansion (\ref{3.36}),
it follows that
$$
\aligned
\frac{64a_k}{\mu_k}\left(\int_{\mathbb{R}^2}
\frac{|z|^2}{(1+|z|^2)^4}
\right)
c_{kl}(\xi')
+
O\left(\frac{1}{|\log\varepsilon|}
\sum\limits_{i=1}^m\sum\limits_{j=1}^{2}|c_{ij}(\xi')|
\right)
=0,
\endaligned
$$
and hence $c_{kl}(\xi')=0$
for each $k=1,\ldots,m$  and  $l=1,2$.
\end{proof}

\vspace{1mm}
\vspace{1mm}

In order to solve for critical points of $F_\lambda$,
we  need first to obtain
its  expansion in terms of $\varphi_m(\xi)$
as $\varepsilon$ goes to zero.

\vspace{1mm}
\vspace{1mm}
\vspace{1mm}
\vspace{1mm}

\noindent{\bf Proposition 5.2.}\,\,{\it
With the choice
{\upshape(\ref{2.22})} for the parameters $\mu=(\mu_1,\ldots,\mu_m)$, there exists $\varepsilon_0>0$
such that for
any  $0<\varepsilon<\varepsilon_0$
and any points $\xi=(\xi_1,\ldots,\xi_m)\in\mathcal{O}_d$,
the following expansion uniformly holds
\begin{eqnarray}\label{5.9}
F_\lambda(\xi)
=\frac{4\pi}{p^2\gamma^{2(p-1)}}
\left
[
m
\Big(
4|\log\varepsilon|
-4+2\log8
\Big)
-8\pi \varphi_m(\xi)
+
O\left(\frac{1}{|\log\varepsilon|}\right)
\right
],
\end{eqnarray}
where
$\varphi_m(\xi)$  is given by
{\upshape(\ref{1.3})}.
}

\noindent

\begin{proof}
Taking into account $DI_\varepsilon(V_{\xi'}+\phi_{\xi'})[\phi_{\xi'}]=0$, a Taylor expansion and an integration by parts,
by (\ref{5.5}) we obtain
$$
\aligned
F_\lambda(\xi)
-
J^p_\lambda\left(U_\xi\right)
&=\frac{1}{p^2\gamma^{2(p-1)}}\int_0^1D^2I_\varepsilon(V_{\xi'}+t\phi_{\xi'})\phi_{\xi'}^2(1-t)dt\\
&=\frac{1}{p^2\gamma^{2(p-1)}}\int_0^1\left\{
\int_{\Omega_\varepsilon}
\big[f'(V_{\xi'})
-f'(V_{\xi'}+t\phi_{\xi'})\big]\phi_{\xi'}^2-
\big[E_{\xi'}+N(\phi_{\xi'})\big]
\phi_{\xi'}
\right\}(1-t)dt,
\endaligned
$$
so we get
\begin{equation}\label{5.11}
\aligned
F_\lambda(\xi)
-
J^p_\lambda\left(U_\xi\right)
=\frac{1}{p^2\gamma^{2(p-1)}}O\left(\frac{1}{|\log\varepsilon|^7}\right),
\endaligned
\end{equation}
in view of $\|\phi_{\xi'}\|_{L^{\infty}(\Omega_\varepsilon)}\leq C|\log\varepsilon|^{-3}$,
$\|E_{\xi'}\|_{*}\leq C|\log\varepsilon|^{-4}$ and
$\|N(\phi_{\xi'})\|_{*}\leq C|\log\varepsilon|^{-6}$
and estimate
(\ref{3.3}).
Next we expand
$$
\aligned
J^p_\lambda\left(U_\xi\right)=
\frac12\int_{\Omega}|\nabla U_\xi|^2
-\frac{\lambda}{p}\int_{\Omega}e^{|U_\xi|^p}:=I_A-I_B.
\endaligned
$$
From the definition  of $U_\xi$ in (\ref{2.16}) we get
$$
\aligned
I_A=
\frac12\int_{\Omega}(-\Delta  U_\xi)U_\xi dx
=\frac12\int_{\Omega}\Big(-\sum\limits_{i=1}^{m}a_i\Delta U_i
\Big)U_\xi dx
=\frac12
\left(
\sum\limits_{k=1}^{m}\int_{B_{d}(\xi_k)}
+\int_{\Omega\setminus\bigcup_{k=1}^m B_{d}(\xi_k)
}
\right)
\Big(-\sum_{i=1}^{m}a_i\Delta U_i\Big)U_\xi dx.
\endaligned
$$
Moreover, by  (\ref{2.8})-(\ref{2.9}),
$$
\aligned
-\sum\limits_{i=1}^{m}a_i\Delta U_i=
\frac1{p\gamma^{p-1}\varepsilon^2}
\sum\limits_{i=1}^{m}a_i
e^{\omega_{\mu_i}\left(\frac{x-\xi_i}{\varepsilon}\right)}\left[
1
+\sum_{j=1}^3\left(\frac{p-1}{p}\right)^j\frac{1}{\gamma^{jp}}
\big(
\omega^j_{\mu_i}-f^j_{\mu_i}
\big)\left(\frac{x-\xi_i}{\varepsilon}\right)
\right].
\endaligned
$$
Applying the expansions  of $U_\xi$ in (\ref{2.20})-(\ref{2.21}),   we can compute
\begin{eqnarray}\label{5.12}
I_A=\,
\frac{1}{2p^2\gamma^{2(p-1)}}
\left\{
\sum\limits_{k=1}^{m}
\frac{1}{\varepsilon^2}
\int_{B_{d}(\xi_k)}
e^{\omega_{\mu_k}\left(\frac{x-\xi_k}{\varepsilon}\right)}\left[
1
+\frac{p-1}{p}\frac{1}{\gamma^{p}}
\big(
\omega^1_{\mu_k}-f^1_{\mu_k}
\big)\left(\frac{x-\xi_k}{\varepsilon}\right)
\right]
\times
\left[
p\gamma^{p}+
\omega_{\mu_k}\left(\frac{x-\xi_k}{\varepsilon}\right)
\right]
dx
\right.
\nonumber
\\[1mm]
\left.
+\,
O\left(\frac{1}{|\log\varepsilon|}\right)
\right\}.
\qquad\qquad\qquad\qquad\qquad\qquad\qquad\qquad\qquad
\qquad\qquad\qquad\qquad\qquad\qquad\qquad\qquad\qquad
\qquad\,\,\,
\end{eqnarray}
Changing  variables $\varepsilon\mu_k z=x-\xi_k$
and using  the relation $p\gamma^{p}=-4\log\varepsilon$, by (\ref{2.3}) we obtain
$$
\aligned
\frac{1}{\varepsilon^2}
\int_{B_{d}(\xi_k)}
e^{\omega_{\mu_k}\left(\frac{x-\xi_k}{\varepsilon}\right)}
\left[
p\gamma^{p}+
\omega_{\mu_k}\left(\frac{x-\xi_k}{\varepsilon}\right)
\right]
dx
=
\int_{B_{d/(\varepsilon\mu_k)}(0)}
\frac{8}{(1+|z|^2)^2}
\left[\log\frac{8}{(1+|z|^2)^2}-
\log\left(\varepsilon^4\mu_k^2\right)
\right]dz.
\endaligned
$$
But
$$
\aligned
\int_{B_{d/(\varepsilon\mu_k)}(0)}\frac{8}{(1+|z|^2)^2}
=8\pi + O\left(\varepsilon^2\right),
\endaligned
$$
and
$$
\aligned
\int_{B_{d/(\varepsilon\mu_k)}(0)}
\frac{8}{(1+|z|^2)^2}
\log\frac{1}{(1+|z|^2)^2}=-16\pi+O(\varepsilon).
\endaligned
$$
Then
\begin{equation}\label{5.13}
\aligned
\frac{1}{\varepsilon^2}
\int_{B_{d}(\xi_k)}
e^{\omega_{\mu_k}\left(\frac{x-\xi_k}{\varepsilon}\right)}
\left[
p\gamma^{p}+
\omega_{\mu_k}\left(\frac{x-\xi_k}{\varepsilon}\right)
\right]
dx
=8\pi\big[\log8-\log\left(\varepsilon^4\mu_k^2\right)-2\big]+
O\left(\varepsilon
\right).
\endaligned
\end{equation}
Changing  variables $\varepsilon\mu_k z=x-\xi_k$ again, by (\ref{2.10}), (\ref{8.3}) and (\ref{8.2}) we obtain
$$
\aligned
&\frac{1}{\varepsilon^2}
\int_{B_{d}(\xi_k)}
e^{\omega_{\mu_k}\left(\frac{x-\xi_k}{\varepsilon}\right)}
\left[
\frac{p-1}{p}\frac{1}{\gamma^{p}}
\big(
\omega^1_{\mu_k}-f^1_{\mu_k}
\big)\left(\frac{x-\xi_k}{\varepsilon}\right)
\right]
\left[
p\gamma^{p}+
\omega_{\mu_k}\left(\frac{x-\xi_k}{\varepsilon}\right)
\right]
dx\\[1mm]
=&\,\int_{B_{d/(\varepsilon\mu_k)}(0)}
\frac{8(p-1)}{(1+|z|^2)^2}
\left\{
\left[
\frac12\big(\upsilon_{\infty}\big)^2
-\omega^0_{\infty}
\right](z)
+(1-2\log\mu_k)\left(
\frac{1-|z|^2}{|z|^2+1}\log8-\frac{2|z|^2}{|z|^2+1}
\right)
\right.\\[1mm]
&\left.
+\,
2(\log^2\mu_k-\log\mu_k)\frac{\,|z|^2-1\,}{\,|z|^2+1\,}
\right\}dz+O\left(\frac{1}{|\log\varepsilon|}\right).
\endaligned
$$
Since
$$
\aligned
\int_{B_{d/(\varepsilon\mu_k)}(0)}
\frac{8|z|^2}{(1+|z|^2)^3}
=4\pi +O(\varepsilon^2)
\qquad\quad
\textrm{and}
\qquad\quad
\int_{B_{d/(\varepsilon\mu_k)}(0)}\frac{8}{(1+|z|^2)^2}
=8\pi + O\left(\varepsilon^2\right),
\endaligned
$$
we give
\begin{eqnarray}\label{5.14}
\frac{1}{\varepsilon^2}
\int_{B_{d}(\xi_k)}
e^{\omega_{\mu_k}\left(\frac{x-\xi_k}{\varepsilon}\right)}
\left[
\frac{p-1}{p}\frac{1}{\gamma^{p}}
\big(
\omega^1_{\mu_k}-f^1_{\mu_k}
\big)\left(\frac{x-\xi_k}{\varepsilon}\right)
\right]
\left[
p\gamma^{p}+
\omega_{\mu_k}\left(\frac{x-\xi_k}{\varepsilon}\right)
\right]
dx
\quad
\nonumber
\\[1mm]
=8\pi(p-1)\left\{
\frac{1}{8\pi}\int_{\mathbb{R}^2}
\frac{8}{(1+|z|^2)^2}
\left[
\frac12\big(\upsilon_{\infty}\big)^2
-\omega^0_{\infty}
\right](z)dz
-1+2\log\mu_k\right\}
+O\left(\frac{1}{|\log\varepsilon|}\right).
\end{eqnarray}
Substituting   (\ref{5.13})-(\ref{5.14}) into  (\ref{5.12}) and using  relation (\ref{5.17}) in
Lemma A.1 as follows:
\begin{eqnarray*}
\frac{1}{8\pi}\int_{\mathbb{R}^2}
\frac{8}{(1+|z|^2)^2}
\left[
\frac12\big(\upsilon_{\infty}\big)^2
-\omega^0_{\infty}
\right](z)dz
=3-\log8,
\end{eqnarray*}
we find
\begin{equation}\label{5.15}
\aligned
I_A=\frac{4\pi}{p^2\gamma^{2(p-1)}}
\left
\{
2(p-2)\sum\limits_{k=1}^{m}\log\mu_k
+m
\Big[
4|\log\varepsilon|
+(p-2)(2-\log8)
\Big]
+
O\left(\frac{1}{|\log\varepsilon|}\right)
\right
\}.
\endaligned
\end{equation}
Regarding the expression $I_B$, by  the definitions of
$\varepsilon$ and $V_{\xi'}$  in  (\ref{1.5}) and   (\ref{2.28}),  respectively,
 we have
$$
\aligned
I_B
=\frac{1}{p^2\gamma^{2(p-1)}}
\left[\sum\limits_{i=1}^m
\left(\int_{B_{\mu_i|\log\varepsilon|^\tau}(\xi'_i)}
+\int_{
B_{d/\varepsilon}(\xi'_i)\setminus
B_{\mu_i|\log\varepsilon|^\tau}(\xi'_i)
}
\right)
+\int_{\Omega_\varepsilon\setminus\cup_{i=1}^mB_{d/\varepsilon}(\xi'_i)}
\right]
e^{\gamma^p\left(\left|
1+\frac{V_{\xi'}(y)}{p\gamma^p}
\right|^p-1\right)}dy.
\endaligned
$$
By (\ref{2.7}) and (\ref{2.36}),
$$
\aligned
\frac{1}{p^2\gamma^{2(p-1)}}
\int_{\Omega_\varepsilon\setminus\cup_{i=1}^mB_{d/\varepsilon}(\xi'_i)}
e^{\gamma^p\left(\left|
1+\frac{V_{\xi'}(y)}{p\gamma^p}
\right|^p-1\right)}dy
&=
\frac{1}{p^2\gamma^{2(p-1)}}O(\varepsilon^{\frac{4}{p}-2})\exp\left[
O\left(\frac{1}{|\log\varepsilon|^{p-1}}\right)
\right]
\\[1mm]
&=
\frac{1}{p^2\gamma^{2(p-1)}}O(\varepsilon^{\frac{2-p}{p}})\exp\left[
-\frac{2-p}{p}|\log\varepsilon|+
O\left(\frac{1}{|\log\varepsilon|^{p-1}}\right)
\right].
\endaligned
$$
By    (\ref{2.6}) and (\ref{2.46}),
$$
\aligned
\frac{1}{p^2\gamma^{2(p-1)}}
\int_{
B_{d/\varepsilon}(\xi'_i)\setminus
B_{\mu_i|\log\varepsilon|^\tau}(\xi'_i)
}
e^{\gamma^p\left(\left|
1+\frac{V_{\xi'}(y)}{p\gamma^p}
\right|^p-1\right)}dy
&
\leq\frac{D}{p^2\gamma^{2(p-1)}}
\int_{
B_{d/\varepsilon}(\xi'_i)\setminus
B_{\mu_i|\log\varepsilon|^\tau}(\xi'_i)
}
e^{\left[
1
+
O\left(\frac{1}{|\log\varepsilon|}\right)
\right]
\omega_{\mu_i}\left(y-\xi'_i\right)
}dy
\\[1mm]
&=\frac{1}{p^2\gamma^{2(p-1)}}O\left(\frac{1}{\,|\log\varepsilon|^{\tau}}\right).
\endaligned
$$
By (\ref{2.13})  and (\ref{2.381}),
$$
\aligned
\frac{1}{p^2\gamma^{2(p-1)}}
\int_{B_{\mu_i|\log\varepsilon|^\tau}(\xi'_i)}
e^{\gamma^p\left(\left|
1+\frac{V_{\xi'}(y)}{p\gamma^p}
\right|^p-1\right)}dy
&=
\frac{1}{p^2\gamma^{2(p-1)}}
\left(\int_{B_{\mu_i|\log\varepsilon|^\tau}(\xi'_i)}
e^{\omega_{\mu_i}\left(y-\xi'_i\right)}
dy
\right)
\left[1
+
O\left(\frac{1}{|\log\varepsilon|}\right)
\right]\\[1mm]
&=\frac{8\pi}{p^2\gamma^{2(p-1)}}\left[
1
+
O\left(\frac{1}{|\log\varepsilon|}\right)
\right].
\endaligned
$$
Then
\begin{equation}\label{5.16}
\aligned
I_B=\frac{8m\pi}{p^2\gamma^{2(p-1)}}
\left[
1
+
O\left(\frac{1}{|\log\varepsilon|}\right)
\right].
\endaligned
\end{equation}
Hence by (\ref{5.11}), (\ref{5.15}) and (\ref{5.16}), we obtain
$$
\aligned
F_\lambda(\xi)=\frac{4\pi}{p^2\gamma^{2(p-1)}}
\left
\{
2(p-2)\sum\limits_{i=1}^{m}\log\mu_i
+m
\Big[
4|\log\varepsilon|-2
+(p-2)(2-\log8)
\Big]
+
O\left(\frac{1}{|\log\varepsilon|}\right)
\right
\},
\endaligned
$$
which, together with the expansion of $\mu_i$ in (\ref{2.24}), gives (\ref{5.9}).
\end{proof}

\vspace{1mm}

Next, we need to prove that  the expansion of
$F_\lambda(\xi)$ in terms of $\varphi_m(\xi)$ holds in a $C^1$-sense.

\vspace{1mm}
\vspace{1mm}
\vspace{1mm}

\noindent{\bf Proposition 5.3.}\,\,{\it
There exists $\varepsilon_0>0$
such that for
any  $0<\varepsilon<\varepsilon_0$
and any points $\xi=(\xi_1,\ldots,\xi_m)\in\mathcal{O}_d$,
the following expansion uniformly holds
\begin{equation}\label{5.27}
\aligned
\nabla_{(\xi_k)_l}F_\lambda(\xi)=-\frac{32\pi^2}{p^2\gamma^{2(p-1)}}
\left\{\nabla_{(\xi_k)_l}\varphi_m(\xi_1,\ldots,\xi_m)
+
O\left(\frac{1}{|\log\varepsilon|}\right)
\right
\},
\endaligned
\end{equation}
where  $k=1,\ldots,m$ and $l=1,2$.
}

\begin{proof}
Observe that for any  $k=1,\ldots,m$ and $l=1,2$,
\begin{equation}\label{5.18}
\aligned
\partial_{(\xi_k)_l}F_\lambda(\xi)=
\frac{1}{\varepsilon p^2\gamma^{2(p-1)}}
\partial_{(\xi'_k)_l}
I_\varepsilon\big(V_{\xi'}+\phi_{\xi'}\big)
=
\frac{1}{\varepsilon p^2\gamma^{2(p-1)}}
\left\{
-\int_{\Omega_\varepsilon}
\big[
\Delta \upsilon_{\xi'}+
f(\upsilon_{\xi'})
\big]
\left(\partial_{(\xi'_k)_l}V_{\xi'}
+\partial_{(\xi'_k)_l}\phi_{\xi'}
\right)
\right\},
\endaligned
\end{equation}
where $\upsilon_{\xi'}=V_{\xi'}+\phi_{\xi'}$.
Let $\eta$ be a  radial smooth, non-increasing cut-off function
such that  $0\leq\eta\leq1$, $\eta=1$ for $|x|\leq  d$,
and $\eta=0$  for $|x|\geq 2 d$.
Since $\phi_{\xi'}$
solves problem (\ref{4.1}),  we deduce
$$
\aligned
-\int_{\Omega_\varepsilon}
\big[
\Delta \upsilon_{\xi'}+
f(\upsilon_{\xi'})
\big]
\partial_{(\xi'_k)_l}\phi_{\xi'}
&=\sum_{i=1}^m\sum_{j=1}^2c_{ij}(\xi')\int_{\Omega_\varepsilon}
e^{
\omega_{\mu_i}\left(y-\xi'_i\right)}Z_{ij}
\partial_{(\xi'_k)_l}\phi_{\xi'}=-\sum_{i=1}^m\sum_{j=1}^2c_{ij}(\xi')
\int_{\Omega_\varepsilon}
\partial_{(\xi'_k)_l}\left(e^{
\omega_{\mu_i}\left(y-\xi'_i\right)}Z_{ij}\right)\phi_{\xi'}
\quad\quad\quad\quad\quad\quad\quad\quad\,
\quad\quad\quad\quad\quad\quad\quad\\
&=\sum_{i=1}^m\sum_{j=1}^2c_{ij}(\xi')\int_{\Omega_\varepsilon}
\partial_{y_l}\left(e^{
\omega_{\mu_i}\left(y-\xi'_i\right)}Z_{ij}\right)
\eta(\varepsilon y-\xi_k)\phi_{\xi'}
\quad\quad\quad\quad\,\,\,\,\,\,\quad
\quad\quad\quad\quad\quad\quad\quad\\
&\quad-\sum_{i=1}^m\sum_{j=1}^2c_{ij}(\xi')\int_{\Omega_\varepsilon}
\left[\partial_{(\xi'_k)_l}\left(e^{
\omega_{\mu_i}\left(y-\xi'_i\right)}Z_{ij}\right)
+\eta(\varepsilon y-\xi_k)\partial_{y_l}\left(e^{
\omega_{\mu_i}\left(y-\xi'_i\right)}Z_{ij}\right)
\right]\phi_{\xi'}.
\endaligned
$$
Since $|\log\varepsilon|^2\,\phi_{\xi'}(y)\rightarrow0$
in $C(\overline{\Omega}_\varepsilon)$,  by
elliptic regularity  we get
\begin{equation}\label{5.19}
\aligned
|\log\varepsilon|^2\,\phi_{\xi'}(y)\rightarrow0
\,\quad\quad\,
\textrm{uniformly in}\,\,\,\,
C^1\big(\overline{\Omega}_\varepsilon\setminus\bigcup_{i=1}^mB_{d/\varepsilon}(\xi'_i)\big).
\endaligned
\end{equation}
Note that $\eta(\varepsilon y-\xi_k)\equiv0$ near $\po_\varepsilon$.
From an integration by parts of the partial derivative in $y_l$
and (\ref{4.2})
we  get
$$
\aligned
\int_{\Omega_\varepsilon}
\partial_{y_l}\left(e^{
\omega_{\mu_i}\left(y-\xi'_i\right)}Z_{ij}\right)
\eta(\varepsilon y-\xi_k)\phi_{\xi'}
&=\int_{\Omega_\varepsilon}
\partial_{y_l}\left(e^{
\omega_{\mu_i}\left(y-\xi'_i\right)}Z_{ij}
\eta(\varepsilon y-\xi_k)\phi_{\xi'}\right)-
\int_{\Omega_\varepsilon}
e^{
\omega_{\mu_i}\left(y-\xi'_i\right)}Z_{ij}
\partial_{y_l}\big(\eta(\varepsilon y-\xi_k)\phi_{\xi'}\big)\\
&=-
\int_{\Omega_\varepsilon}
e^{
\omega_{\mu_i}\left(y-\xi'_i\right)}Z_{ij}
\partial_{y_l}\big(\eta(\varepsilon y-\xi_k)\phi_{\xi'}\big)\\
&=-\int_{B_{d/\varepsilon}(\xi'_k)}
e^{
\omega_{\mu_i}\left(y-\xi'_i\right)}Z_{ij}
\partial_{y_l}\phi_{\xi'}+o\left(\frac{\varepsilon^3}{|\log\varepsilon|^2}\right).
\endaligned
$$
Similar to (\ref{5.8}),  by (\ref{2.3}) and (\ref{3.10}) for $j=1,2$ we get
$$
\aligned
\partial_{(\xi'_k)_l}&\left(
e^{
\omega_{\mu_i}\left(y-\xi'_i\right)
}
Z_{ij}\right)
+\eta(\varepsilon y-\xi_k)\partial_{y_l}\left(e^{
\omega_{\mu_i}\left(y-\xi'_i\right)}Z_{ij}\right)\\[1mm]
=&\,4\mu_ie^{
\omega_{\mu_i}\left(y-\xi'_i\right)}
\left\{
\frac{\eta(\varepsilon y-\xi_k)-\delta_{ik}}{|y-\xi'_i|^2+\mu_i^2}\delta_{jl}
+6
\frac{(y-\xi'_i)_j(y-\xi'_k)_l}{(|y-\xi'_i|^2+\mu_i^2)^2}
\big(\delta_{ik}-\eta(\varepsilon y-\xi_k)\big)
+
\frac{3}{4\mu_i}Z_{i0}Z_{ij}
\partial_{(\xi'_k)_l}\log\mu_i
\right\}\\[1mm]
=&\,3e^{
\omega_{\mu_i}\left(y-\xi'_i\right)}Z_{i0}Z_{ij}
\partial_{(\xi'_k)_l}\log\mu_i+O\left(\varepsilon^6\right),
\endaligned
$$
which, together with  (\ref{4.2}) and the fact that
$|Z_{i0}Z_{ij}|\leq2$ and $|\partial_{(\xi'_k)_l}\log\mu_i|=O\left(\varepsilon\right)$
for any $i=1,\ldots,m$, gives
$$
\aligned
\left|
\sum_{i=1}^m\sum_{j=1}^2c_{ij}(\xi')\int_{\Omega_\varepsilon}
\left[\partial_{(\xi'_k)_l}\left(e^{
\omega_{\mu_i}\left(y-\xi'_i\right)}Z_{ij}\right)
+\eta(\varepsilon y-\xi_k)\partial_{y_l}\left(e^{
\omega_{\mu_i}\left(y-\xi'_i\right)}Z_{ij}\right)
\right]\phi_{\xi'}
\right|
=O\left(\frac{\varepsilon}{|\log\varepsilon|^7}\right).
\endaligned
$$
Hence
\begin{eqnarray}\label{5.20}
-\int_{\Omega_\varepsilon}\big[
\Delta \upsilon_{\xi'}+
f(\upsilon_{\xi'})
\big]\partial_{(\xi_k)_l}\phi_{\xi'}
=-\sum_{i=1}^m\sum_{j=1}^2c_{ij}(\xi')\int_{B_{d/\varepsilon}(\xi'_k)}
e^{
\omega_{\mu_i}\left(y-\xi'_i\right)}Z_{ij}
\partial_{y_l}\phi_{\xi'}
+
O\left(\frac{\varepsilon}{|\log\varepsilon|^7}\right)
\nonumber\\
=\int_{B_{d/\varepsilon}(\xi'_k)}
\big[
\Delta \upsilon_{\xi'}+
f(\upsilon_{\xi'})
\big]
\partial_{y_l}\phi_{\xi'}
+O\left(\frac{\varepsilon}{|\log\varepsilon|^7}\right).
\quad\qquad\qquad\,\,\,\,
\end{eqnarray}
Taking into account
$(D_{\xi'}+D_{y})V_{\xi'}(y)=p\gamma^{p-1}(D_{\xi'}+D_{y})\big[U_\xi(\varepsilon y)\big]$,
by the expression of $U_\xi$ in (\ref{2.16}) we obtain
\begin{equation}\label{5.33}
\aligned
\partial_{(\xi'_k)_l}V_{\xi'}
+\partial_{y_l}V_{\xi'}
=\sum_{i=1}^{m}a_i\left\{
\Big(\partial_{(\xi'_k)_l}
+
\partial_{y_l}\Big)
\left[
\omega_{\mu_i}\left(y-\xi'_i\right)
+\sum_{j=1}^3\left(\frac{p-1}{p}\right)^j\frac{1}{\gamma^{jp}}
\omega^j_{\mu_{i}}
\left(y-\xi'_i\right)
+p\gamma^{p-1}H_i(\varepsilon y)\right]\right\}.
\qquad
\endaligned
\end{equation}
Clearly,
by (\ref{5.29}),
\begin{equation}\label{5.34}
\aligned
\Big(\partial_{(\xi'_k)_l}
+
\partial_{y_l}\Big)\big[\omega_{\mu_i}\left(y-\xi'_i\right)\big]
=
\frac1{\mu_i}\big(\delta_{ki}-1\big)
Z_{l}\left(
\frac{y-\xi'_i}{\mu_i}
\right)
+\left(\partial_{(\xi'_k)_l}\mu_i\right)\frac{d}{d\mu_i}\omega_{\mu_i}\left(y-\xi'_i\right),
\endaligned
\end{equation}
and for $j=1,2,3$, by (\ref{5.30}),
\begin{equation}\label{5.35}
\aligned
\Big(\partial_{(\xi'_k)_l}
+
\partial_{y_l}\Big)\big[\omega^j_{\mu_i}\left(y-\xi'_i\right)\big]
=
-\frac1{\mu_i}\big(\delta_{ki}-1\big)
\left[\frac{D^j_{\mu_i}}{4}Z_{l}\left(
\frac{y-\xi'_i}{\mu_i}
\right)+
O\left(\frac{\mu_i^2}{|y-\xi'_i|^2+\mu_i^2}\right)
\right]
+\left(\partial_{(\xi'_k)_l}\mu_i\right)\frac{d}{d\mu_i}\omega^j_{\mu_i}\left(y-\xi'_i\right),
\endaligned
\end{equation}
and by  (\ref{2.19}) and (\ref{5.31}),
\begin{equation}\label{5.36}
\aligned
\Big(\partial_{(\xi'_k)_l}
+
\partial_{y_l}\Big)\big[p\gamma^{p-1}H_i(\varepsilon y)\big]
=O\left(\varepsilon\right).
\endaligned
\end{equation}
Since
$|\partial_{(\xi'_k)_l}\mu_i|=O\left(\varepsilon\right)$
for all $i=1,\ldots,m$,
by inserting (\ref{5.34})-(\ref{5.36}) into (\ref{5.33}) we can derive that for  $|y-\xi'_k|\leq
2d/\varepsilon$,
$$
\aligned
\eta(\varepsilon y-\xi_k)\left[\partial_{(\xi'_k)_l}V_{\xi'}
+\partial_{y_l}V_{\xi'}
\right]
=O\left(\varepsilon\right).
\endaligned
$$
On the other hand,
by (\ref{5.7}) we have  that
for   $|y-\xi'_k|\geq d/\varepsilon$,
$$
\aligned
\big[1-\eta(\varepsilon y-\xi_k)\big]\partial_{(\xi'_k)_l}V_{\xi'}
=O\left(\varepsilon\right).
\endaligned
$$
So,
$$
\aligned
\partial_{(\xi'_k)_l}V_{\xi'}
+\eta(\varepsilon y-\xi_k)\partial_{y_l}V_{\xi'}
=\big[1-\eta(\varepsilon y-\xi_k)\big]\partial_{(\xi'_k)_l}V_{\xi'}
+\eta(\varepsilon y-\xi_k)\left[\partial_{(\xi'_k)_l}V_{\xi'}
+\partial_{y_l}V_{\xi'}
\right]
=O\left(\varepsilon\right).
\endaligned
$$
By  (\ref{4.2}),
\begin{eqnarray}\label{5.22}
-\int_{\Omega_\varepsilon}
\big[
\Delta \upsilon_{\xi'}+
f(\upsilon_{\xi'})
\big]
\partial_{(\xi'_k)_l}V_{\xi'}
=\sum_{i=1}^m\sum_{j=1}^2c_{ij}(\xi')\int_{\Omega_\varepsilon}
e^{
\omega_{\mu_i}\left(y-\xi'_i\right)}Z_{ij}
\left[
\partial_{(\xi'_k)_l}V_{\xi'}
+\eta(\varepsilon y-\xi_k)\partial_{y_l}V_{\xi'}
\right]
\nonumber
\\
-\sum_{i=1}^m\sum_{j=1}^2c_{ij}(\xi')\int_{\Omega_\varepsilon}
e^{
\omega_{\mu_i}\left(y-\xi'_i\right)}Z_{ij}
\eta(\varepsilon y-\xi_k)\partial_{y_l}V_{\xi'}
\quad\quad\quad\quad
\,\,\,\,\,\,
\nonumber\\
=-\sum_{i=1}^m\sum_{j=1}^2c_{ij}(\xi')\int_{B_{d/\varepsilon}(\xi'_k)}
e^{
\omega_{\mu_i}\left(y-\xi'_i\right)}
Z_{ij}
\partial_{y_l}V_{\xi'}
+O\left(\frac{\varepsilon}{|\log\varepsilon|^4}\right)
\,\,
\nonumber\\
=\int_{B_{d/\varepsilon}(\xi'_k)}
\big[
\Delta \upsilon_{\xi'}+
f(\upsilon_{\xi'})
\big]
\partial_{y_l}V_{\xi'}+O\left(\frac{\varepsilon}{|\log\varepsilon|^4}\right).
\quad\quad\quad
\quad\qquad
\end{eqnarray}
Substituting (\ref{5.20}) and (\ref{5.22}) into (\ref{5.18}), we conclude
\begin{equation}\label{5.23}
\aligned
\partial_{(\xi_k)_l}F_\lambda(\xi)=
\frac{1}{\varepsilon p^2\gamma^{2(p-1)}}
\left\{
\int_{B_{d/\varepsilon}(\xi'_k)}
\big[
\Delta \upsilon_{\xi'}+
f(\upsilon_{\xi'})
\big]
\partial_{y_l}\upsilon_{\xi'}
+O\left(\frac{\varepsilon}{|\log\varepsilon|^4}\right)
\right\}.
\endaligned
\end{equation}

Integrating by parts of the gradient
operator $\nabla$
and the partial derivative
 $\partial_{y_l}$ respectively,
we
obtain  the
Pohozaev-type identities: for any  $B\subset\Omega_\varepsilon$ and for
any function $\upsilon$,
\begin{equation}\label{5.24}
\aligned
\int_{B}\Delta \upsilon \partial_{y_l} \upsilon
=\int_{\partial B}\left(
\partial_\nu \upsilon \partial_{y_l} \upsilon
-\frac12|\nabla \upsilon|^2\nu_l
\right),
\quad\ \quad
\quad\ \quad
\int_{B}f(\upsilon )\partial_{y_l} \upsilon
=\int_{\partial B}
e^{\gamma^p\left(\left|1+\frac{\upsilon}{p\gamma^p}\right|^p-1\right)}\nu_l,
\,\quad
\endaligned
\end{equation}
where $\nu_l(y)$ denotes the $l$-th component of
the outer unit normal vector to
$\partial B$ at $y\in\partial B$.
Let
$$
\aligned
\psi_k( y)=H (\varepsilon y,\xi_k)+\sum_{i=1,\,i\neq k}^m
a_{i}a_{k}
G (\varepsilon y,\xi_i)
\,\quad
\textrm{for all}
\,\,\,k=1,\ldots,m.
\endaligned
$$
From (\ref{2.20}), (\ref{2.28}),  (\ref{5.19}) and the fact
that $\upsilon_{\xi'}=V_{\xi'}+\phi_{\xi'}$,  we can derive that
the follows expansions uniformly hold:
\begin{equation}\label{5.25}
\aligned
\upsilon_{\xi'}(y)=8\pi\sum_{i=1}^ma_iG(\varepsilon y,\xi_i)-p\gamma^p+O\left(\frac1{|\log\varepsilon|}\right)
=8\pi a_k\left[
\frac1{2\pi}\log\frac1{|\varepsilon y-\xi_k|}+\psi_k(y)
\right]-p\gamma^p
+O\left(\frac1{|\log\varepsilon|}\right),
\endaligned
\end{equation}
but
\begin{equation}\label{5.26}
\aligned
\nabla\upsilon_{\xi'}(y)
=8\pi a_k\left[
-\frac{\varepsilon}{2\pi}\frac{\varepsilon y-\xi_k}{|\varepsilon y-\xi_k|^2}+\nabla\psi_k(y)
\right]
+
O\left(\frac{\varepsilon}{|\log\varepsilon|}\right),
\endaligned
\end{equation}
 in
$C\big(\overline{\Omega}_\varepsilon\setminus\bigcup_{i=1}^mB_{d/\varepsilon}(\xi'_i)\big)$.
Applying (\ref{5.24}) on $B=B_{d/\varepsilon}(\xi'_k)$
for each $k=1,\ldots,m$ and using (\ref{5.25})-(\ref{5.26}), we conclude
$$
\aligned
\int_{B_{d/\varepsilon}(\xi'_k)}\big[
\Delta \upsilon_{\xi'}+
f(\upsilon_{\xi'})
\big]
\partial_{y_l}\upsilon_{\xi'}
=&\int_{\partial B_{d/\varepsilon}(\xi'_k)}
\left[
\partial_\nu \upsilon_{\xi'}
\partial_{y_l}
\upsilon_{\xi'}
-\frac12|\nabla \upsilon_{\xi'}|^2\nu_l+
e^{\gamma^p\left(\left|1+\frac{\upsilon}{p\gamma^p}\right|^p-1\right)}\nu_l
\right]\\
=&\,64\pi^2\left[\int_{\partial B_{d/\varepsilon}(\xi'_k)}\left(
-\frac{\varepsilon}{2\pi d}+\partial_\nu \psi_k
\right)\left(-\frac{\varepsilon}{2\pi d}\nu_l
+\partial_{y_l} \psi_k
\right)\right.\\
&\left.-\frac12\int_{\partial B_{d/\varepsilon}(\xi'_k)}\left|-\frac{\varepsilon}{2\pi d}\nu
+\nabla \psi_k
\right|^2\nu_l\right]
+O\left(\frac{\varepsilon}{|\log\varepsilon|}\right)\\
=&\,64\pi^2\int_{\partial B_{d/\varepsilon}(\xi'_k)}
\left[-\frac{\varepsilon}{2\pi d}\partial_{y_l}\psi_k+\left(
\partial_{\nu}\psi_k\partial_{y_l}\psi_k
-\frac12|\nabla\psi_k|^2\nu_l
\right)
\right]
+O\left(\frac{\varepsilon}{|\log\varepsilon|}\right)\\[1mm]
=&
-64\pi^2
\partial_{y_l} \psi_k(\xi'_k)
+O\left(\frac{\varepsilon}{|\log\varepsilon|}\right),
\endaligned
$$
because $\psi_k$ is a  harmonic function  on the ball  $B_{d/\varepsilon}(\xi'_k)$ such that
$$
\aligned
\frac{\varepsilon}{2\pi d}\int_{\partial B_{d/\varepsilon}(\xi'_k)}
\partial_{y_l}\psi_k
=\partial_{y_l}\psi_k(\xi'_k),
\endaligned
$$
and by (\ref{5.24}),
$$
\aligned
\int_{\partial B_{d/\varepsilon}(\xi'_k)}
\left(
\partial_{\nu}\psi_k\partial_{y_l}\psi_k
-\frac12|\nabla\psi_k|^2\nu_l
\right)
=\int_{B_{d/\varepsilon}(\xi'_k)}\Delta \psi_k\partial_{y_l}\psi_k=0.
\endaligned
$$
Accordingly, by (\ref{5.23}) we find
$$
\aligned
\partial_{(\xi_k)_l}F_\lambda(\xi)
=-
\frac{64\pi^2}{\varepsilon p^2\gamma^{2(p-1)}}
\left\{
\partial_{y_l} \psi_k(\xi'_k)
+O\left(\frac{\varepsilon}{|\log\varepsilon|}\right)
\right
\}
=
-\frac{32\pi^2}{p^2\gamma^{2(p-1)}}
\left\{
\partial_{(\xi_k)_l} \varphi_m(\xi)
+
O\left(\frac{1}{|\log\varepsilon|}\right)
\right
\},
\endaligned
$$
in view of $\partial_{y_l} \psi_k(\xi'_k)=\frac12\varepsilon\,\partial_{(\xi_k)_l} \varphi_m(\xi)$.
This completes the proof.
\end{proof}

\vspace {1mm}

\section{Proofs of theorems}
\noindent {\bf Definition 6.1.}
We say that $\xi^*=(\xi^*_1,\ldots,\xi^*_m)$ is a $C^1$-stable critical point
of $\varphi_m:\mathcal{F}_m(\Omega)\rightarrow \mathbb{R}$
if for any sequence of functions
$\Phi_n:\mathcal{F}_m(\Omega)\rightarrow \mathbb{R}$
such that $\Phi_n\rightarrow \varphi_m$
uniformly in $C_{loc}^1\big(\mathcal{F}_m(\Omega)\big)$,
$\Phi_n$ has at least one critical point $\xi^n=(\xi^n_1,\ldots,\xi^n_m)$
such that $\xi^n\rightarrow\xi^*$ as $n\rightarrow+\infty$.
Specially, $\xi^*$ is a $C^1$-stable critical point
of $\varphi_m$ if either one of the following conditions is satisfied:\\
\indent {\upshape (i)}\,\,
  $\xi^*$ is an isolated local maximum or minimum  point of $\varphi_m$; \\
\indent {\upshape (ii)}\,\,the Brouwer degree $\deg(\nabla\varphi_m, B_\varepsilon(\xi^*),0 )\neq0$
for any $\varepsilon>0$ small enough.

\vspace {1mm}
\vspace {1mm}
\vspace {1mm}

\noindent{\bf Lemma 6.2 (see Lemma $2.2$ of \cite{BP3}
and Lemma $2$ in \cite{H}).}\,\,\,{\it
Let $u$ be a solution  of
$-\Delta u=g$ in $\Omega$,
$u=0$
on $\po$. If $\Lambda$
is a neighbourhood of $\po$, then
$$
\aligned
\|\nabla u\|_{C^{0,\alpha}(\Lambda')}
\leq
C
(
\|g\|_{L^{1}(\Omega)}
+
\|g\|_{L^{\infty}(\Lambda)}
),
\endaligned
$$
where $\alpha\in(0,1)$
and $\Lambda'\subset\subset\Lambda$
is a neighbourhood of $\po$.
}

\vspace {1mm}
\vspace {1mm}
\vspace {1mm}

\noindent {\bf Proof of Theorem 1.1.}
According to Proposition  5.1, the function
$u_{\lambda}=U_{\xi}+\widetilde{\phi}_{\xi}$
is a solution to problem (\ref{1.1}) if we adjust
$\xi=(\xi_1,\ldots,\xi_m)\in\mathcal{O}_d\subset\mathcal{F}_m(\Omega)$   so that it is
a critical point of $F_\lambda$ defined in (\ref{5.2}).
This is equivalent to showing that
\begin{equation}\label{6.1}
\aligned
\widetilde{F}_\lambda(\xi)=
-\frac{\,p^2\gamma^{2(p-1)}\,}{32\pi^2}F_\lambda(\xi)
+\frac{m}{8\pi}
\Big(
4|\log\varepsilon|
-4+2\log8
\Big).
\endaligned
\end{equation}
has a critical point
$\xi^\varepsilon=(\xi^\varepsilon_1,\ldots,\xi^\varepsilon_{m})\in\mathcal{O}_d$.
Propositions 5.2-5.3 imply that for any small but fixed $d>0$, as $\lambda\rightarrow0$,
\begin{equation}\label{6.2}
\aligned
\widetilde{F}_\lambda(\xi)=\varphi_m(\xi)
+
O\left(\frac{1}{|\log\varepsilon|}\right)
\qquad
\textrm{uniformly in}
\,\,\,\,
C^1(\mathcal{O}_d).
\endaligned
\end{equation}
Since  $\xi^*=(\xi^*_1,\ldots,\xi^*_m)$ is a $C^1$-stable critical point
of $\varphi_m:\mathcal{F}_m(\Omega)\rightarrow \mathbb{R}$,
by Definition 6.1 it follows that there exists
at least one critical point $\xi^\varepsilon$ of $\widetilde{F}_\lambda$
in
$\mathcal{O}_d$
such that along a subsequence, $\xi^\varepsilon\rightarrow\xi^*$
as $\lambda\rightarrow0$. The function
$u_{\lambda}=U_{\xi^\varepsilon}+\widetilde{\phi}_{\xi^\varepsilon}$
is therefore a weak solution to problem (\ref{1.1})
with the qualitative property
(\ref{1.8})
which easily follows by  (\ref{2.20}),
(\ref{2.21})
and the fact that $\widetilde{\phi}_{\xi^\varepsilon}$
is a higher order term in $u_{\lambda}$
because of (\ref{5.3}) and the first estimate of  $\phi_{(\xi^\varepsilon)'}$
in (\ref{4.2}).

\vspace{1mm}

\noindent {\bf Proof of (\ref{1.4}).}\,
From
(\ref{2.20}), (\ref{2.21}), (\ref{4.2}) and (\ref{5.3}) we
obtain that estimate (\ref{1.4}) holds pointwise in
$\overline{\Omega}\setminus\{\xi^*_1,\ldots,\xi^*_m\}$. We
will try to prove
that
\begin{equation}\label{6.7}
\aligned
\big\|
\lambda u_\lambda|u_\lambda|^{p-2}e^{|u_\lambda|^p}
\big
\|_{L^{1}(\Omega)}
\leq
\frac{C}{p\gamma^{p-1} },
\endaligned
\end{equation}
and for a given neighbourhood $\Lambda$ of $\po$,
\begin{equation}\label{6.8}
\aligned
\big
\|\lambda u_\lambda|u_\lambda|^{p-2}e^{|u_\lambda|^p}
\big
\|_{L^{\infty}(\Lambda)}
\leq
\frac{C}{p\gamma^{p-1} }.
\endaligned
\end{equation}
Once these estimates are established, according to
Lemma $6.2$
we get that
for $\Lambda'\subset\subset
\overline{\Omega}\setminus\{\xi^*_1,\ldots,\xi^*_m\}$,
$\|\nabla u_\lambda\|_{C^{0,\alpha}(\Lambda')}\leq
\frac{C}{p\gamma^{p-1} }$, hence estimate (\ref{1.4})  follows by
(\ref{2.7})  and
the Ascoli-Arzel\'{a} Theorem.
Notice that if we set $\upsilon_\lambda(x)=p\gamma^{p-1}u_\lambda(x)-p\gamma^{p}$,
then
$\upsilon_\lambda(x)=\big(V_{(\xi^\varepsilon)'}+\phi_{(\xi^\varepsilon)'}\big)
(\frac{x}{\varepsilon} )$
and
\begin{equation}\label{6.9}
\aligned
\lambda u_\lambda|u_\lambda|^{p-2}e^{|u_\lambda|^p}
=
\lambda\gamma^{p-1}e^{\gamma^p}
\left(1+\frac{\upsilon_\lambda}{p\gamma^p}\right)
\left|1+\frac{\upsilon_\lambda}{p\gamma^p}\right|^{p-2}
e^{\gamma^p\left(\left|1+\frac{\upsilon_\lambda}{p\gamma^p}\right|^p-1\right)}
=\lambda\gamma^{p-1}e^{\gamma^p}f\big(V_{(\xi^\varepsilon)'}+\phi_{(\xi^\varepsilon)'}\big)
\Big(\frac{x}{\varepsilon}\Big),
\endaligned
\end{equation}
where the function $f$ is defined in (\ref{2.27}).
Recall that $\|\phi_{(\xi^\varepsilon)'}\|_{L^{\infty}(\Omega_\varepsilon)}\leq\frac{C}{|\log\varepsilon|^3}$.
Similar to the consideration of
those expansions in (\ref{2.37}), (\ref{2.39}), (\ref{2.43})  and (\ref{2.46}),
we can compute  that
if $|y-(\xi_i^\varepsilon)'|\geq d/\varepsilon$ for all $i=1,\ldots,m$,
\begin{equation}\label{6.10}
\aligned
\left|f\big(V_{(\xi^\varepsilon)'}+\phi_{(\xi^\varepsilon)'}\big)
\right|(y)
=\frac{O(\varepsilon^{\frac{2+p}{p}})}{\,|\log\varepsilon|^{p-1}\,}
\exp\left[
-\frac{2-p}{p}|\log\varepsilon|+
O\left(\frac{1}{|\log\varepsilon|^{p-1}}\right)
\right],
\endaligned
\end{equation}
if
$|y-(\xi_i^\varepsilon)'|<\mu_i|\log\varepsilon|^\tau$ with any $\tau\geq10$ large but fixed,
\begin{equation}\label{6.11}
\aligned
f\big(V_{(\xi^\varepsilon)'}+\phi_{(\xi^\varepsilon)'}\big)
(y)
=a_i e^{\omega_{\mu_i}\left(y-(\xi_i^\varepsilon)'\right)}
\left[1+
O
\left(\frac{\log^2(\mu_i+|y-(\xi_i^\varepsilon)'|)}{|\log\varepsilon|}\right)
\right],
\endaligned
\end{equation}
and if $\mu_i|\log\varepsilon|^\tau\leq
|y-(\xi_i^\varepsilon)'|\leq d/\varepsilon^\theta$
 with any
 $\theta<1$ but close enough to $1$,
\begin{equation}\label{6.12}
\aligned
\left|f\big(V_{(\xi^\varepsilon)'}+\phi_{(\xi^\varepsilon)'}\big)
\right|(y)
\leq
D\left(1+\frac{1}{|\log\varepsilon|^{p-1}}\right)e^{\left[
1
+
O\left(\frac{1}{|\log\varepsilon|}\right)
\right]
\omega_{\mu_i}\left(y-(\xi_i^\varepsilon)'\right)
},
\endaligned
\end{equation}
where  $D>0$ is a constant, independent of any $\theta<1$.
Therefore,  by using relations (\ref{1.5})-(\ref{2.7}) we can easily derive that
$$
\aligned
\big\|
\lambda u_\lambda|u_\lambda|^{p-2}e^{|u_\lambda|^p}
\big
\|_{L^{1}(\Omega)}
=&
\left[
\int_{\Omega_\varepsilon\setminus\cup_{i=1}^m B_{d/\varepsilon}\big((\xi_i^\varepsilon)'\big)}
+\sum\limits_{i=1}^m
\left(\int_{B_{\mu_i|\log\varepsilon|^\tau}\big((\xi_i^\varepsilon)'\big)}
+\int_{
B_{d/\varepsilon}\big((\xi_i^\varepsilon)'\big)\setminus
B_{\mu_i|\log\varepsilon|^\tau}\big((\xi_i^\varepsilon)'\big)
}
\right)
\right]
\\[1.5mm]
&
\,
\quad\,\,\,
\lambda\gamma^{p-1}e^{\gamma^p}
\varepsilon^2
\left|f\big(V_{(\xi^\varepsilon)'}+\phi_{(\xi^\varepsilon)'}\big)
\right|(y)dy\\[0.1mm]
=&\,
\lambda\gamma^{p-1}e^{\gamma^p}
\varepsilon^2
\left[
\sum\limits_{i=1}^m
\int_{B_{\mu_i|\log\varepsilon|^\tau}\big((\xi_i^\varepsilon)'\big)}
e^{\omega_{\mu_i}\left(y-(\xi_i^\varepsilon)'\right)}
dy
+
O
\left(\frac{1}{|\log\varepsilon|}\right)
\right]\\[0.1mm]
=&\,
\frac{1}{p\gamma^{p-1} }
\left[
8m\pi
+
O
\left(\frac{1}{|\log\varepsilon|}\right)
\right].
\endaligned
$$
Moreover,
$$
\aligned
\big
\|\lambda u_\lambda|u_\lambda|^{p-2}e^{|u_\lambda|^p}
\big
\|_{L^{\infty}\big(\Omega\setminus\cup_{i=1}^m B_{d}(\xi_i^\varepsilon)\big)}
=\frac{1}{\varepsilon^2p\gamma^{p-1} }
\frac{O(\varepsilon^{\frac{2+p}{p}})}{\,|\log\varepsilon|^{p-1}\,}
\exp\left[
-\frac{2-p}{p}|\log\varepsilon|+
O\left(\frac{1}{|\log\varepsilon|^{p-1}}\right)
\right]
=
o\left(\frac{1}{p\gamma^{p-1} }\right).
\endaligned
$$

\vspace{1mm}

\noindent {\bf Proof of (\ref{1.6}).}\,
From (\ref{6.9})-(\ref{6.12}) we conclude that for any $\psi\in C_{c}(\oo)$,
$$
\aligned
\int_{\Omega}
\lambda u_\lambda|u_\lambda|^{p-2}e^{|u_\lambda|^p}
\psi
dx
=&
\left[
\int_{\Omega_\varepsilon\setminus\cup_{i=1}^m B_{d/\varepsilon}\big((\xi_i^\varepsilon)'\big)}
+\sum\limits_{i=1}^m
\left(\int_{B_{\mu_i|\log\varepsilon|^\tau}\big((\xi_i^\varepsilon)'\big)}
+\int_{
B_{d/\varepsilon}\big((\xi_i^\varepsilon)'\big)\setminus
B_{\mu_i|\log\varepsilon|^\tau}\big((\xi_i^\varepsilon)'\big)
}
\right)
\right]
\\[1.5mm]
&
\,
\quad\,\,\,
\lambda\gamma^{p-1}e^{\gamma^p}
\varepsilon^2
f\big(V_{(\xi^\varepsilon)'}+\phi_{(\xi^\varepsilon)'}\big)
(y)\,\psi(\varepsilon y)dy\\[0.1mm]
=&\,
\lambda\gamma^{p-1}e^{\gamma^p}
\varepsilon^2
\left[\,
\sum\limits_{i=1}^ma_i
\int_{B_{\mu_i|\log\varepsilon|^\tau}\big((\xi_i^\varepsilon)'\big)}
e^{\omega_{\mu_i}\left(y-(\xi_i^\varepsilon)'\right)}
\,\psi(\varepsilon y)dy
+
O
\left(\frac{1}{|\log\varepsilon|}\right)
\right]\\[0.1mm]
=&\,
\frac{1}{p\gamma^{p-1} }
\left[\,8\pi
\sum\limits_{i=1}^ma_i
\psi(\xi_i^*)
+
o
\left(1\right)
\right].
\endaligned
$$

\vspace{1mm}

\noindent {\bf Proof of (\ref{1.7})-(\ref{1.7b}).}\,
Since $\widetilde{\phi}_{\xi^\varepsilon}$
is a higher order term in $u_{\lambda}$,
by using (\ref{2.381}) and the same calculus   as expansion
(\ref{5.16}) we get
$$
\aligned
\frac{\lambda p}{2}\int_{\Omega}
\big(e^{|u_\lambda|^p}-1\big)dx=&\frac{1}{2\gamma^{2(p-1)}}
\left\{\sum_{i=1}^m
\int_{B_{\mu_i|\log\varepsilon|^\tau}((\xi_i^\varepsilon)')}
e^{\omega_{\mu_i}\left(y-(\xi_i^\varepsilon)'\right)}
\left[1+\frac{p-1}{p}\frac1{\gamma^p}
B_1
+\left(\frac{p-1}{p}\right)^2\frac1{\gamma^{2p}}
\left(B_2+\frac12(B_1)^2\right)
\right]
dy
\right.
\\[1mm]
&+\left.
O\left(\frac{1}{|\log\varepsilon|^3}\right)
\right\}.
\endaligned
$$
Then
$$
\aligned
&\left(\frac{\lambda p}{2}\int_{\Omega}
\big(e^{|u_\lambda|^p}-1\big)dx
\right)^{\frac{2-p}{p}}
=
\left(
\frac{1}{2\gamma^{2(p-1)}}
\right)^{\frac{2-p}{p}}\left[1+
O\left(\frac{1}{|\log\varepsilon|^3}\right)
\right]
\\
&
\qquad\qquad
\times
\left\{
\sum_{i=1}^m
\underbrace{\int_{\mathbb{R}^2}
e^{\omega_{\mu_i}\left(y-(\xi_i^\varepsilon)'\right)}
\left[1+\frac{p-1}{p}\frac1{\gamma^p}
B_1
+\left(\frac{p-1}{p}\right)^2\frac1{\gamma^{2p}}
\left(B_2+\frac12(B_1)^2\right)
\right]
dy}\limits_{e_i}\right\}^{\frac{2-p}{p}}.
\endaligned
$$
Similar  to these arguments, by (\ref{2.380})-(\ref{2.381}) we obtain
$$
\aligned
&\qquad\qquad\qquad\qquad
\left(\frac{\lambda p}{2}\int_{\Omega}|u_\lambda|^p e^{|u_\lambda|^p}dx
\right)^{\frac{2(p-1)}{p}}
=\left(\frac{1}{2\gamma^{p-2}}
\right)^{\frac{2(p-1)}{p}}
\left[1+
O\left(\frac{1}{|\log\varepsilon|^3}\right)
\right]
\\
&
\times\left\{\sum_{i=1}^m
\underbrace{\int_{\mathbb{R}^2}
e^{\omega_{\mu_i}\left(y-(\xi_i^\varepsilon)'\right)}
\left[1
+\frac{1}{p\gamma^p}
\Big(\,pA_1+(p-1)B_1\Big)
+\frac{p-1}{p\gamma^{2p}}
\big(A_1B_1+B_1\big)
+
\left(\frac{p-1}{p}\right)^2\frac1{\gamma^{2p}}
\left(B_2+\frac12(B_1)^2\right)
\right]
dy}\limits_{f_i}
\right\}^{\frac{2(p-1)}{p}}.
\endaligned
$$
Thanks to the vanishing  identity of first order stated by Corollary A.8
\begin{eqnarray}\label{8.18}
\frac{\,2-p\,}{p}
\int_{\mathbb{R}^2}
e^{\omega_{\mu_i}\left(y-(\xi_i^\varepsilon)'\right)}
B_1
dy
+
\frac{2}{\,p\,}\int_{\mathbb{R}^2}
e^{\omega_{\mu_i}\left(y-(\xi_i^\varepsilon)'\right)}
\big[\,p A_1+(p-1)B_1\big]
dy
\,\equiv\,
0,
\end{eqnarray}
by the Taylor expansion we can compute
\begin{eqnarray}\label{6.6}
\beta_\lambda=
\left(\frac{\lambda p}{2}
\int_{\Omega}
\big(e^{|u_\lambda|^p}-1\big)dx
\right)^{\frac{2-p}{p}}
\left(\frac{\lambda p}{2}
\int_{\Omega}
|u_\lambda|^p
e^{|u_\lambda|^p}dx
\right)^{\frac{2(p-1)}{p}}
\qquad\qquad\qquad\qquad\qquad\qquad\qquad\qquad\quad\quad\quad \,\,\,
&&
\nonumber
\\[1mm]
=
4\pi m
\left\{1+
\frac{(p-1)^2}{ p^2\gamma^{2p}}\frac{1}{8\pi m}\sum_{i=1}^m
\left[\int_{\mathbb{R}^2}
e^{\omega_{\mu_i}\left(y-(\xi_i^\varepsilon)'\right)}
\Big(B_2+\frac12(B_1)^2\Big)
dy
+2
\int_{\mathbb{R}^2}
e^{\omega_{\mu_i}\left(y-(\xi_i^\varepsilon)'\right)}
(A_1B_1+B_1)
dy
\right]
\right.
&&
\nonumber
\\[1mm]
+\left.
\frac{(p-1)(p-2)}{ p^2\gamma^{2p}}\frac{1}{(16\pi m)^2}
\left(\sum_{i=1}^m\int_{\mathbb{R}^2}
e^{\omega_{\mu_i}\left(y-(\xi_i^\varepsilon)'\right)}
B_1
dy
\right)^2\right\}
\left[1+
O\left(\frac{1}{|\log\varepsilon|^3}\right)
\right].
\qquad\qquad\qquad\qquad\quad\,\,\,
&&
\end{eqnarray}
Obviously, by the definitions of $\varepsilon$ and $\gamma$
in  (\ref{1.5})-(\ref{2.7}) we arrive at
\begin{equation*}\label{6.16}
\aligned
\beta_\lambda=
4\pi m\left[
1
+
O\left(\frac{1}{|\log\varepsilon|^2}\right)
\right]\rightarrow 4\pi m.
\endaligned
\end{equation*}
If $0<p<1$,  using the inequality
$$
\aligned
\left(\sum_{i=1}^m\int_{\mathbb{R}^2}
e^{\omega_{\mu_i}\left(y-(\xi_i^\varepsilon)'\right)}
B_1
dy
\right)^2\leq
m\sum_{i=1}^m\left(\int_{\mathbb{R}^2}
e^{\omega_{\mu_i}\left(y-(\xi_i^\varepsilon)'\right)}
B_1
dy
\right)^2,
\endaligned
$$
by (\ref{6.6}) we find
\begin{eqnarray*}\label{6.14}
\beta_\lambda\leq
4\pi m
\left\{1+
\frac{p-1}{\,p^2\gamma^{2p}\,m\,}
\sum_{i=1}^m
\left[\frac{p-1}{8\pi}
\left(\int_{\mathbb{R}^2}
e^{\omega_{\mu_i}\left(y-(\xi_i^\varepsilon)'\right)}
\Big(B_2+\frac12(B_1)^2\Big)
dy
+2
\int_{\mathbb{R}^2}
e^{\omega_{\mu_i}\left(y-(\xi_i^\varepsilon)'\right)}
(A_1B_1+B_1)
dy
\right)
\right.
\right.
&&
\nonumber
\\[1mm]
\left.\left.
+
\frac{p-2}{\,(16\pi)^2\,}
\left(
\int_{\mathbb{R}^2}
e^{\omega_{\mu_i}\left(y-(\xi_i^\varepsilon)'\right)}
B_1
dy
\right)^2\right]
\right\}
\left[1+
O\left(\frac{1}{|\log\varepsilon|^3}\right)
\right]
\qquad\qquad\qquad\qquad\quad\quad
\qquad\qquad\qquad\qquad\qquad
&&\nonumber
\\[1mm]
=\,4\pi m
\left\{1+
\frac{4(p-1)}{\,p^2\gamma^{2p}\,}
\left[1+
O\left(\frac{1}{|\log\varepsilon|}\right)
\right]\right\}
<4\pi m,
\qquad\qquad\qquad\qquad
\qquad\qquad\qquad\qquad\qquad
\quad\qquad\qquad\qquad\,\,
&&
\end{eqnarray*}
where the above equality is due to  the identity  of second order  stated by Corollary B.6
\begin{eqnarray}\label{9.10}
\frac{p-1}{8\pi}
\left[\int_{\mathbb{R}^2}
e^{\omega_{\mu_i}\left(y-(\xi_i^\varepsilon)'\right)}
\Big(B_2+\frac12(B_1)^2\Big)
dy
+2
\int_{\mathbb{R}^2}
e^{\omega_{\mu_i}\left(y-(\xi_i^\varepsilon)'\right)}
(A_1B_1+B_1)
dy
\right]
+\frac{p-2}{(16\pi)^2}
\left(
\int_{\mathbb{R}^2}
e^{\omega_{\mu_i}\left(y-(\xi_i^\varepsilon)'\right)}
B_1
dy
\right)^2\equiv4.
\end{eqnarray}
While if $1<p<2$,  by  H\"{o}lder's inequality for vectors in $\mathbb{R}_{+}^m$ we  conclude
\begin{eqnarray*}
\beta_\lambda=
\frac{1}{2}\left(
\sum_{i=1}^{m}e_i
\right)^{\frac{2-p}{p}}
\left(
\sum_{i=1}^{m}f_i
\right)^{\frac{2(p-1)}{p}}
\left[1+
O\left(\frac{1}{|\log\varepsilon|^3}\right)
\right]
\geq
\frac{1}{2}
\left(
\sum_{i=1}^m
e_i^{\frac{2-p}{p}}f_i^{\frac{2(p-1)}{p}}
\right)
\left[1+
O\left(\frac{1}{|\log\varepsilon|^3}\right)
\right].
\end{eqnarray*}
Applying the Taylor expansion and using the identities (\ref{8.18}) and (\ref{9.10}) again, we
can compute
\begin{eqnarray*}
\left(
\frac{e_i}{8\pi}
\right)^{\frac{2-p}{p}}
\left(
\frac{f_i}{8\pi}
\right)^{\frac{2(p-1)}{p}}
=
1+
\frac{(p-1)^2}{ p^2\gamma^{2p}}\frac{1}{8\pi}
\left[\int_{\mathbb{R}^2}
e^{\omega_{\mu_i}\left(y-(\xi_i^\varepsilon)'\right)}
\Big(B_2+\frac12(B_1)^2\Big)
dy
+2
\int_{\mathbb{R}^2}
e^{\omega_{\mu_i}\left(y-(\xi_i^\varepsilon)'\right)}
(A_1B_1+B_1)
dy
\right]
&&
\nonumber
\\[1mm]
+\,
\frac{(p-1)(p-2)}{ p^2\gamma^{2p}}\frac{1}{(16\pi)^2}
\left(\int_{\mathbb{R}^2}
e^{\omega_{\mu_i}\left(y-(\xi_i^\varepsilon)'\right)}
B_1
dy
\right)^2+
O\left(\frac{1}{|\log\varepsilon|^3}\right)
\qquad\qquad\qquad\qquad\quad\quad
&&\nonumber
\\[1mm]
=1+
\frac{\,4(p-1)\,}{ p^2\gamma^{2p}}
+
O\left(\frac{1}{|\log\varepsilon|^3}\right).
\qquad\qquad\qquad\qquad\qquad\qquad\qquad\qquad\qquad
\qquad\qquad\quad\qquad\,
&&
\end{eqnarray*}
Hence for $1<p<2$,
\begin{eqnarray*}\label{6.17}
\beta_\lambda\geq
4\pi m
\left\{1+
\frac{4(p-1)}{\,p^2\gamma^{2p}\,}
\left[1+
O\left(\frac{1}{|\log\varepsilon|}\right)
\right]\right\}
>4\pi m.
\end{eqnarray*}

\vspace{1mm}
\vspace{1mm}
\vspace{1mm}

\noindent {\bf Proof of (\ref{1.9}).}\,
Let us first claim  that
 for all
points $\xi=(\xi_1,\ldots,\xi_m)\in \mathcal{F}_m(\Omega)$,
the normal derivative
\begin{equation}\label{6.13}
\aligned
\frac{\partial}{\partial\nu}\left[\sum\limits_{i=1}^{m}
 a_iG(x,\xi_i)\right]\not\equiv0
 \,
 \quad
 \,
 \textrm{on}
  \,\,\,\,
 \po.
\endaligned
\end{equation}
By contradiction, we suppose
that there exists an $m$-tuple $\xi=(\xi_1,\ldots,\xi_m)\in \mathcal{F}_m(\Omega)$
such that
$\sum_{i=1}^{m}
 a_i\frac{\partial G}{\partial\nu}(x,\xi_i)=0$  for all $x\in\po$.
Applying Green's formula
to any harmonic function $\psi$ in $\Omega$, we find
$$
\aligned
\sum\limits_{i=1}^{m}
 a_i\psi(\xi_i)
=
-\int_{\Omega}\left[\sum\limits_{i=1}^{m}
 a_iG(x,\xi_i)\right]\Delta\psi dx
+\int_{\po}
\left\{
\left[\sum\limits_{i=1}^{m}
 a_iG(x,\xi_i)\right]
 \frac{\partial\psi}{\partial\nu}
 -
 \psi
\left[\sum\limits_{i=1}^{m}
 a_i\frac{\partial G}{\partial\nu}
 (x,\xi_i)\right]
 \right\}
 d\sigma_x\equiv0.
\endaligned
$$
Then
$$
\aligned
\sum\limits_{i=1}^{m}
 a_i\xi_i=0,
 \qquad\quad
\sum\limits_{i=1}^{m}
 a_i(\xi_i)_1(\xi_i)_2=0
\qquad\quad
  \textrm{and}
\qquad\quad
\sum\limits_{i=1}^{m}
 a_i\log|\xi_i-q|=0
 \quad
\textrm{for all}
\,\,\,
 q\in\mathbb{R}^2\setminus\Omega.
\endaligned
$$
Since  $\sum_{i=1}^ma_i=0$
with $a_i\in\{-1,1\}$, we find that
 these identities should be absurd for any
points $\xi=(\xi_1,\ldots,\xi_m)\in \mathcal{F}_m(\Omega)$.

Assume that there exists a sequence $\lambda_n\rightarrow0$
such that
$\overline{\{x\in\Omega:\,u_{\lambda_n}(x)=0\}}\cap\po=\emptyset$.
Clearly,
 $\frac{\partial u_{\lambda_n}}{\partial\nu}$
does not change sign on $\po$.
Thus by (\ref{1.4}),
 $\sum_{i=1}^{m}
 a_i\frac{\partial G}{\partial\nu}(x,\xi_i^*)$
does not change sign on $\po$. On the other hand,
a simple calculus gives
\begin{equation}\label{6.15}
\aligned
\int_{\po}
\frac{\partial}{\partial\nu}\left[\sum\limits_{i=1}^{m}
 a_iG(x,\xi_i^*)\right]
 d\sigma_x
 =
 \sum\limits_{i=1}^{m}
 a_i
 \int_{\Omega}\Delta G(x,\xi_i^*) dx=-\sum\limits_{i=1}^{m}
 a_i=0.
\endaligned
\end{equation}
This, together with (\ref{6.13}), implies that
$\sum_{i=1}^{m}
 a_i\frac{\partial G}{\partial\nu}(x,\xi_i^*)$
 must change sign on $\po$, which is a contradiction.
\qquad\qquad\qquad\qquad$\square$

\vspace{1mm}
\vspace{1mm}
\vspace{1mm}
\vspace{1mm}

\noindent {\bf Proof of Theorem 1.2.}
Fix $m=2$, $a_1=1$ and $a_2=-1$.
Then the function $\varphi_2$
defined in (\ref{1.3}) becomes
$$
\aligned
\varphi_2(\xi)
=
H(\xi_1,\xi_1)+
H(\xi_2,\xi_2)
-2
  G(\xi_1,\xi_2)
\,\quad
\,\,\,\,\,\,
\textrm{for all}
\,\,\,\,
\xi=(\xi_1,\xi_2)\in\mathcal{F}_2(\Omega),
\endaligned
$$
where $
\mathcal{F}_2(\Omega)=\{\xi=(\xi_1,\xi_2)\in\Omega\times\Omega
:\,\xi_1\neq \xi_2\}$.

\vspace{1mm}

\noindent {\bf Proof of Theorem 1.2(i).}
According to Proposition 5.1,
we need to prove that if $\lambda$ is small enough,
the function $F_{\lambda}$ has at least
$\cat\big(\mathcal{C}_2(\Omega)\big)$
pairs of critical points. By (\ref{6.1}) it
reduces to  prove that
$\widetilde{F}_{\lambda}$ has at least
$\cat\big(\mathcal{C}_2(\Omega)\big)$
pairs of critical points.
Notice that
\begin{equation}\label{6.3}
\aligned
\varphi_2(\xi)
\rightarrow
  -\infty
\,\quad
\,
\textrm{as}
\,\,\,\,
\xi
\rightarrow\partial\mathcal{F}_2(\Omega).
\endaligned
\end{equation}
Recall that $\mathcal{C}_2(\Omega)$
denotes the quotient manifold of $\mathcal{F}_2(\Omega)$
modulo the equivalence $(\xi_1,\xi_2)\sim(\xi_2,\xi_1)$.
Under the map $(\xi_1,\xi_2)\rightarrow(\xi_2,\xi_1)$,
we get $U_{\xi}\rightarrow-U_{\xi}$ and
$\widetilde{\phi}_{\xi}\rightarrow-\widetilde{\phi}_{\xi}$,
and further  $\widetilde{F}_{\lambda}(\xi_1,\xi_2)=\widetilde{F}_{\lambda}(\xi_2,\xi_1)$
for  any  $(\xi_1,\xi_2)\in\mathcal{F}_2(\Omega)$.
The induced functions $\widehat{\widetilde{F}}_{\lambda}$,
$\widehat{\varphi}_2:\mathcal{C}_2(\Omega)\rightarrow\mathbb{R}$
are well defined. Setting
$k:=\cat\big(\mathcal{C}_2(\Omega)\big)$,
we observe that there exists a compact subset $K_0\subset\mathcal{C}_2(\Omega)$
such that $\cat(K_0)=k$.
From (\ref{6.3})
we see that
the
 upper level set $\widehat{\varphi}_2^a=\{\xi\in\mathcal{C}_2(\Omega):\,
\widehat{\varphi}_2(\xi)\geq a\}$ is compact for any $a\in\mathbb{R}$.
We take $a<\min_{K_0}\widehat{\varphi}_2$ and consider
 $\widehat{\widetilde{F}}_{\lambda}$ on the
compact manifold $K=\widehat{\varphi}_2^a$ with boundary $B=\widehat{\varphi}_2^{-1}(a)$.
Clearly,
by (\ref{6.2}) we have  that
$\widehat{\widetilde{F}}_{\lambda}\rightarrow\widehat{\varphi}_2$, $C^1$-uniformly
 on  compact subset  of $\mathcal{C}_2(\Omega)$.
If $\lambda$ is small enough, it follows that
$\max_{B}\widehat{\widetilde{F}}_{\lambda}<\min_{K_0}\widehat{\widetilde{F}}_{\lambda}$.
Standard critical point theory implies that $\widehat{\widetilde{F}}_{\lambda}$
has at least $k$
distinct  critical points in $K$. From
(\ref{6.2}) and Palais's principle of symmetric criticality (see \cite{P}),
 $\widetilde{F}_{\lambda}$ has at least
$k$
pairs $(\xi_{1}^{i,\varepsilon},\,\xi_{2}^{i,\varepsilon})$,
$(\xi_{2}^{i,\varepsilon},\,\xi_{1}^{i,\varepsilon})$
of critical points with
$i=1,\ldots,k$
such that each
 $\xi^{i,\varepsilon}=(\xi_{1}^{i,\varepsilon},\,\xi_{2}^{i,\varepsilon})$
converges  along a subsequence towards a critical  point
$\xi^i=(\xi_{1}^i,\,\xi_{2}^i)$
of $\varphi_2$ in $\mathcal{F}_2(\Omega)$.

\vspace{1mm}

\noindent {\bf Proof of Theorem 1.2(ii).}
Let $u_{\lambda}=U_{\xi^\varepsilon}+\widetilde{\phi}_{\xi^\varepsilon}$
be any one of weak solutions to problem (\ref{1.1})
found in Theorem 1.2(i), where
 $\xi^\varepsilon=(\xi_{1}^{\varepsilon},\,\xi_{2}^{\varepsilon})$
is a critical point of $F_\lambda$ in $\mathcal{O}_d$
such that   it converges along a subsequence
towards a critical  point
$\xi^*=(\xi_{1}^*,\,\xi_{2}^*)$
of $\varphi_2$ in $\mathcal{F}_2(\Omega)$.
Observe that
for $|x-\xi_{i}^{\varepsilon}|\leq r$ with
 $0<r<d$ and  $i=1,2$, by
(\ref{2.7})
and (\ref{2.13}),
$$
\aligned
p\gamma^{p}+
\omega_{\mu_i}\left(\frac{x-\xi_i^\varepsilon}{\varepsilon}\right)
+\sum_{j=1}^3\left(\frac{p-1}{p}\right)^j\frac{1}{\gamma^{jp}}\omega^j_{\mu_{i}}\left(\frac{x-\xi_i^\varepsilon}{\varepsilon}\right)
&\geq p\gamma^{p}+\log\frac{8}{\mu_i^2}
+\left[
-2+
\frac{p-1}{p}\frac{D^1_{\mu_i}}{2\gamma^{p}}
\right]\log\left(1+\frac{r^2}{\varepsilon^2\mu_i^2}\right)+O\left(\frac{1}{|\log\varepsilon|}\right)\\
&=\log\frac{8}{\mu_i^2}
-4\log\frac{ r}{\mu_i}
+\frac{(p-1)D^1_{\mu_i}}4
+O\left(\frac{1}{|\log\varepsilon|}\right).
\endaligned
$$
Then by  (\ref{2.6}), (\ref{2.15}), (\ref{2.21}), (\ref{4.2}) and (\ref{5.3}),
it is easily checked that,
choosing $r>0$ smaller if  necessary,
there exists $\delta>0$ such that
$p\gamma^{p-1} u_\lambda(x)>\delta$ for any  $x\in B_r(\xi_1^\varepsilon)$,
$p\gamma^{p-1} u_\lambda(x)<-\delta$ for any  $x\in B_r(\xi_2^\varepsilon)$.
Moreover,
$\dist(B_r(\xi_1^\varepsilon),\,B_r(\xi_2^\varepsilon))\geq \delta$,
$\xi_1^*\in B_r(\xi_1^\varepsilon)$
and $\xi_2^*\in B_r(\xi_2^\varepsilon)$
for any $\varepsilon$ small enough.
So the set $\Omega\setminus\{x\in\Omega:\,u_{\lambda}(x)=0\}$
has at least two connected components.
Set  $\hat{u}_\lambda(x)=p\gamma^{p-1} u_\lambda(x)$
and
$\Omega_r=\Omega\setminus\big[B_r(\xi_1^\varepsilon)\cup B_r(\xi_2^\varepsilon)
\big]$. By contradiction, we assume that there exists a third connected component
$\Lambda_\varepsilon$ of  $\Omega\setminus\{x\in\Omega:\,u_{\lambda}(x)=0\}$.
Then $\Lambda_\varepsilon\subset\subset\Omega_r$ and $\hat{u}_\lambda\in H_0^1(\Lambda_\varepsilon)$
is a weak solution of the equation
\begin{equation}\label{6.4}
\aligned
-\Delta \hat{u}_\lambda=\frac{\lambda}{\,(p\gamma^{p-1})^{p-2}\, }
\hat{u}_\lambda|\hat{u}_\lambda|^{p-2}e^{\frac{1}{\,(p\gamma^{p-1})^{p}\,} |\hat{u}_\lambda|^p}
\,\,\,\ \ \,\,
\textrm{in}\,\,\,\,\,
\Lambda_\varepsilon.
\endaligned
\end{equation}
Since  $\Lambda_\varepsilon\subset\subset\Omega_r$, by
(\ref{2.20}), (\ref{2.21}), (\ref{4.2}) and (\ref{5.3}) we
find that as $\varepsilon$ tends to zero,
\begin{equation}\label{6.5}
\aligned
\hat{u}_\lambda
\rightarrow
8\pi \big[ G(x,\xi_1^*)-G(x,\xi_2^*)\big]
\,\,\,
\quad
\textrm{uniformly over}
\,\,\,\,
\overline{\Lambda}_\varepsilon,
\endaligned
\end{equation}
so $\sup_{\Lambda_\varepsilon}|\hat{u}_\lambda|\leq C<+\infty$
for any $\varepsilon$ small enough.
Testing equation (\ref{6.4}) against $\hat{u}_\lambda$
and using the definitions of
$\varepsilon$ and $\gamma$ in (\ref{1.5})-(\ref{2.7}), we
can directly check that  for any $0<p<2$,
$$
\aligned
\|\hat{u}_\lambda\|^2_{H_0^1(\Lambda_\varepsilon)}
=
\frac{\lambda}{\,(p\gamma^{p-1})^{p-2}\, }
\int_{\Lambda_\varepsilon}
|\hat{u}_\lambda|^{p}e^{\frac{1}{\,(p\gamma^{p-1})^{p}\,} |\hat{u}_\lambda|^p}
=\frac{O(\varepsilon^{\frac{2-p}{p}})}{|\log\varepsilon|^{p-1}}
\exp\left[
-\frac{2-p}{p}|\log\varepsilon|+
O\left(\frac{1}{|\log\varepsilon|^{p-1}}\right)
\right]
\rightarrow0
\,\,\,
\quad
\textrm{as}
\,\,\,\,
\varepsilon\rightarrow0.
\endaligned
$$
From the compactness of Sobolev embedding
$H_0^1(\Lambda_\varepsilon)\hookrightarrow
L^2(\Lambda_\varepsilon)$, it follows that $\hat{u}_\lambda\rightarrow0$
in $L^2(\Lambda_\varepsilon)$. On the other hand, from
(\ref{6.5}) and  Lebesgue's theorem, we obtain
$\hat{u}_\lambda
\rightarrow
8\pi \big[ G(x,\xi_1^*)-G(x,\xi_2^*)\big]$
in  $L^2(\Lambda_\varepsilon)$, which is a contradiction.

\vspace{1mm}

\noindent {\bf Proof of Theorem 1.2(iii).}
This part  is an immediate consequence of (\ref{1.9}) because of
$a_1+a_2=0$.
\,\qquad\qquad\qquad\qquad\qquad$\square$

\vspace {1mm}
\vspace {1mm}
\vspace {1mm}
\vspace {1mm}

\noindent {\bf Proof of Theorem 1.3.}
For the special case $m=3$ or $m=4$ and $a_i=(-1)^{i+1}$,
$i=1,\ldots,m$, we consider the existence of
 critical points of
the function  $\widetilde{F}_\lambda$
as in (\ref{6.1}).
From (\ref{6.2}) it follows that
$\widetilde{F}_{\lambda}\rightarrow\varphi_m$, $C^1$-uniformly
 on any compact subset $\mathcal{O}_d$ of $\mathcal{F}_m(\Omega)$.
According to Theorems 2.2-2.3 in \cite{BP2},
 $\widetilde{F}_{\lambda}$ has a critical point
$\xi^\varepsilon=(\xi^\varepsilon_1,\ldots,\xi^\varepsilon_{m})\in\mathcal{O}_d$
which converges along a subsequence towards a critical point
$\xi^*=(\xi^*_1,\ldots,\xi^*_{m})$ of $\varphi_m$ in $\mathcal{F}_m(\Omega)$.
\,\,\,\qquad\qquad\qquad\qquad\qquad\qquad\qquad\qquad\qquad\qquad\qquad$\square$

\vspace {1mm}
\vspace {1mm}
\vspace {1mm}
\vspace {1mm}

\noindent {\bf Proof of Theorem 1.4.}
For the general case $m\geq1$  and $a_i=(-1)^{i+1}$,
$i=1,\ldots,m$, we consider the existence of
 critical points of
the function  $\widetilde{F}_\lambda$
as in (\ref{6.1}).
Set $\Omega^S=\Omega\cap(\mathbb{R}\times\{0\})\neq \emptyset$.
From (\ref{6.2}) it follows that
$\widetilde{F}_{\lambda}\rightarrow\varphi_m$, $C^1$-uniformly
 on  any compact subset of $\mathcal{F}_m(\Omega^S)$.
Since
$\Omega$ is symmetric with respect to the reflection at
$\mathbb{R}\times\{0\}$,
by Theorem 3.3 in \cite{BPW}
we get that
 $\widetilde{F}_{\lambda}$ has a critical point
$\xi^\varepsilon=(\xi^\varepsilon_1,\ldots,\xi^\varepsilon_{m})\in\mathcal{F}_m(\Omega^S)$
which converges along a subsequence towards a critical point
$\xi^*=(\xi^*_1,\ldots,\xi^*_{m})\in\mathcal{F}_m(\Omega^S)$ of $\varphi_m$
with $\xi_{i}^*=(t_i,\,0)$ and
$t_1<t_2<\cdots<t_m$.
\,\,\,\qquad\qquad\qquad\qquad\qquad\qquad\qquad\qquad\qquad
\qquad\qquad\qquad\qquad\qquad\,\,$\square$

\vspace{1mm}
\vspace{1mm}
\vspace{1mm}
\vspace{1mm}

\section{Proof of the vanishing identity  (\ref{8.18}) of first order}
According to \cite{CI,EMP,MM}, for a radial function $f(y)=f(|y|)$ there exists a unique radial solution
\begin{equation}\label{8.1}
\aligned
\omega(r)=\frac{1-r^2}{1+r^2}
\left(
\int_{0}^{r}\frac{\phi_f(s)-\phi_f(1)}{(s-1)^2}ds
+\phi_{f}(1)\frac{r}{1-r}
\right)
\endaligned
\end{equation}
for the equation
$$
\aligned
\Delta\omega+\frac{8}{(1+|y|^2)^2}\omega=\frac{8}{(1+|y|^2)^2}f(y)
\qquad
\textrm{in}\,\,\ \,\,\mathbb{R}^2,
\qquad
\omega(0)=0,
\endaligned
$$
where
$$
\aligned
\phi_f(s)=8\left(\frac{s^2+1}{s^2-1}\right)^2\frac{(s-1)^2}{s}
\int_{0}^{s}t\frac{1-t^2}{(t^2+1)^3}f(t)dt
\quad\,
\textrm{for}\,\,\,\,s\neq1,
\qquad
\textrm{but}\ \,\,\,\,
\phi_f(1)=\lim_{s\rightarrow1}\phi_f(s).
\endaligned
$$
Furthermore, if $f$ is the smooth function with at most logarithmic growth at infinity,
then a derect computation shows that
\begin{equation}\label{7.1}
\aligned
\omega(r)=\frac{D_{f}}{2}\log\left(1+r^2\right)+C_{f}+O\left(\frac{1}{1+r}\right),
\,\,\qquad\,\,
\partial_{r}\omega(r)=\frac{r\,D_{f}}{\,1+r^2\,}+O\left(\frac{1}{1+r^2}\right)
\quad\,\,\,\textrm{as}\,\,\,r\rightarrow+\infty,
\endaligned
\end{equation}
where
$$
\aligned
D_{f}=\frac{1}{2\pi}\int_{\mathbb{R}^2}
\Delta\omega(y) dy
\qquad\qquad
\textrm{and}
\qquad\qquad
D_{f}=8\int_{0}^{+\infty}t\frac{t^2-1}{(t^2+1)^3}f(t)dt.
\endaligned
$$

Making the change to variables $z=\mu_iy$, we denote
\begin{equation}\label{8.3}
\aligned
\widetilde{\omega}_{\mu_i}(y):=\omega_{\mu_i}(\mu_iy),
\qquad
\widetilde{\omega}^1_{\mu_i}(y):=\omega^1_{\mu_i}(\mu_iy),
\qquad
\widetilde{f}^j_{\mu_i}(y):=f^j_{\mu_i}(\mu_iy),
\qquad
\upsilon_{\infty}(y):=\widetilde{\omega}_{\mu_i}(y)+2\log\mu_i=\log\frac{8}{(1+|y|^2)^2}
\endaligned
\end{equation}
with $j\in\{1,\,2\}$.
Let $\omega^0_{\infty}$,  $\omega^1_{\infty}$ and $\omega^2_{\infty}$  be the
radial  solutions of
$$
\aligned
\Delta\omega^j_{\infty}+\frac{8}{(1+|y|^2)^2}\omega^j_{\infty}=\frac{8}{(1+|y|^2)^2}f_j(y)
\quad\,
\textrm{in}\,\,\ \,\,\mathbb{R}^2,
\qquad
j=0,1,2,
\endaligned
$$
where
$$
\aligned
f_0(y)=\frac12\big(\upsilon_{\infty}(y)\big)^2,
\qquad\qquad
f_1(y)=\upsilon_{\infty}(y),
\qquad\qquad
f_2(y)=1.
\endaligned
$$
Obviously,
$$
\aligned
\omega^2_{\infty}(y)=
1-Z_0(y)=\frac{2}{|y|^2+1}.
\endaligned
$$
Applying the formulas (\ref{8.1})-(\ref{7.1}) and replacing $\omega(r)$ with $\omega(r)-C_{f}Z_0(r)$, we can compute
\begin{equation}\label{8.4}
\aligned
\omega^0_{\infty}(y)=&\,\frac12\big(\upsilon_{\infty}(y)\big)^2+6\log(|y|^2+1)
+\frac{2\log8-10}{|y|^2+1}+\frac{|y|^2-1}{|y|^2+1}
\left[4\int_{|y|^2}^{+\infty}\frac{\log(s+1)}{s(s+1)}ds
-2\log^2(|y|^2+1)-\frac12\log^28\right],
\endaligned
\end{equation}
and
\begin{equation}\label{8.5}
\aligned
\omega^1_{\infty}(y)=
\frac{|y|^2-1}{|y|^2+1}
\left\{
\frac{2}{|y|^2-1}\big[
\upsilon_{\infty}(y)
+|y|^2
\big]
+
\upsilon_{\infty}(y)
-
\log8
-2
\right\}.
\endaligned
\end{equation}
By (\ref{2.10}) we obtain
$$
\aligned
\widetilde{f}^1_{\mu_i}(y)=
-\big[
f_0(y)
+
(1-2\log\mu_i)f_1(y)+2(\log^2\mu_i-\log\mu_i)f_2(y)
\big],
\endaligned
$$
and hence
\begin{equation}\label{8.6}
\aligned
\widetilde{\omega}^1_{\mu_i}(y)=
-\omega^0_{\infty}(y)
-
(1-2\log\mu_i)\omega^1_{\infty}(y)-4(\log^2\mu_i-\log\mu_i)\frac{\,1\,}{\,|y|^2+1\,}.
\endaligned
\end{equation}

\vspace {1mm}
\vspace {1mm}
\vspace {1mm}
\vspace {1mm}

\noindent{\bf Lemma A.1.}\,\,\,{\it
\begin{eqnarray}\label{5.17}
\frac{1}{8\pi}\int_{\mathbb{R}^2}
\frac{8}{(1+|z|^2)^2}
\left[
\frac12\big(\upsilon_{\infty}\big)^2
-\omega^0_{\infty}
\right](z)dz
=3-\log8.
\end{eqnarray}
}

\noindent{\it Proof.}
Applying the divergence theorem and (\ref{8.4}), we deduce
$$
\aligned
\qquad\,\,
\frac{1}{8\pi}\int_{\mathbb{R}^2}
\frac{8}{(1+|z|^2)^2}
\left[
\frac12\big(\upsilon_{\infty}\big)^2
-\omega^0_{\infty}
\right](z)dz
=
\frac{1}{8\pi}\int_{\mathbb{R}^2}
\Delta\omega^0_{\infty}
=
\frac{1}{8\pi}
\lim_{r\rightarrow\infty}
\int_{B_r(0)}
\Delta\omega^0_{\infty}
=\frac{1}{4}
\lim_{r\rightarrow\infty}r\big(\omega^0_{\infty}\big)'(r)
=3-\log8.
\qquad\,
\square
\endaligned
$$

\vspace {1mm}
\vspace {1mm}
\vspace {1mm}
\vspace {1mm}

\noindent{\bf Lemma A.2.}\,\,\,{\it
\begin{equation}\label{3.26}
\aligned
\int_{\mathbb{R}^2}
\frac{8\mu_i^2}{(\mu_i^2+|z|^2)^2}
\left[Z_0\left(\frac{z}{\mu_i}\right)\right
]^2
\left[
1+\omega^{1}_{\mu_i}+\frac{1}{2}(\omega_{\mu_i})^2+2\omega_{\mu_i}
\right](z)dz=8\pi.
\endaligned
\end{equation}
}

\noindent{\it Proof.}
From relations (\ref{8.3})-(\ref{8.6})  we get
\begin{eqnarray}\label{8.2}
\left[
1+\omega^{1}_{\mu_i}+\frac{1}{2}(\omega_{\mu_i})^2+2\omega_{\mu_i}
\right](\mu_i y)=\left[
\frac12\big(\upsilon_{\infty}\big)^2
+\upsilon_{\infty}-\omega_{\infty}^0
\right](y)+
\Big[
2\log^2\mu_i
-2\log\mu_i
+(\log8+1
)
(1-2\log\mu_i
)
\Big]
\frac{|y|^{2}-1}{|y|^{2}+1}.
\end{eqnarray}
Then
$$
\aligned
&\int_{\mathbb{R}^2}
\frac{8\mu_i^2}{(\mu_i^2+|z|^2)^2}
\left[Z_0\left(\frac{z}{\mu_i}\right)\right
]^2
\left[
1+\omega^{1}_{\mu_i}+\frac{1}{2}(\omega_{\mu_i})^2+2\omega_{\mu_i}
\right](z)dz\\[2mm]
=&
\int_{\mathbb{R}^2}
\frac{8(|y|^2-1)^2}{(|y|^2+1)^4}\left\{\left[
\frac12\big(\upsilon_{\infty}\big)^2
+\upsilon_{\infty}-\omega^0_{\infty}
\right](y)+\Big[
2\log^2\mu_i
-2\log\mu_i
+(\log8+1
)
(1-2\log\mu_i
)
\Big]
\frac{|y|^{2}-1}{|y|^{2}+1}
\right\}dy.
\endaligned
$$
Notice that
$$
\aligned
\int_{\mathbb{R}^2}
\frac{8(|y|^2-1)^2}{(|y|^2+1)^4}dy=\frac{8\pi}{3}
\qquad\quad
\textrm{and}
\qquad\quad
\int_{\mathbb{R}^2}
\frac{8(|y|^2-1)^2}{(|y|^2+1)^4}\frac{|y|^2}{|y|^2+1}dy=\frac{4\pi}{3}.
\endaligned
$$
Applying the explicit expression of $\omega^0_{\infty}(y)$ in   (\ref{8.4})
and  integrating by parts,  we can compute
$$
\aligned
\int_{\mathbb{R}^2}
\frac{8(|y|^2-1)^2}{(|y|^2+1)^4}\left[
\frac12\big(\upsilon_{\infty}\big)^2
+\upsilon_{\infty}-\omega^0_{\infty}
\right](y)dy=8\pi,
\endaligned
$$
(also see \cite{EMP} on Page $50$). Hence
from all these computations  we can easily deduce that (\ref{3.26}) holds.
\qquad\qquad\qquad\qquad\qquad
$\square$

\vspace {1mm}
\vspace {1mm}
\vspace {1mm}

\noindent{\bf Lemma A.3.}\,\,\,{\it
\begin{eqnarray}
\int_{\mathbb{R}^2}
e^{\omega_{\mu_i}}
\omega_{\mu_i}
=8\pi\left(\log8-2\log\mu_i-2\right),
\qquad\qquad\qquad\qquad
\qquad\qquad\qquad\qquad
&&
\label{8.8}\\
\int_{\mathbb{R}^2}
e^{\omega_{\mu_i}}
(\omega_{\mu_i})^2
=8\pi\left[4\log^2\mu_i-4(\log8-2)\log\mu_i
+\log^28-4\log8+8
\right],
\qquad\qquad\qquad\qquad
&&
\label{8.9}\\
\int_{\mathbb{R}^2}
e^{\omega_{\mu_i}}
(\omega_{\mu_i})^3
=8\pi\big[-8\log^3\mu_i
+12(\log8-2)\log^2\mu_i
-6(\log^2 8-4\log8 +8)\log\mu_i
+\log^38-6\log^28+24\log8+48
\big].
&&
\label{8.10}
\end{eqnarray}
}

\noindent{\it Proof.}
A direct computation gives
$$
\aligned
\int_{\mathbb{R}^2}
\frac{1}{(1+|z|^2)^2}\log(1+|z|^{2})
=\pi,
\,\qquad\,
\int_{\mathbb{R}^2}
\frac{1}{(1+|z|^2)^2}\log^2(1+|z|^{2})
=2\pi,
\,\qquad\,
\int_{\mathbb{R}^2}
\frac{1}{(1+|z|^2)^2}\log^3(1+|z|^{2})
=6\pi.
\endaligned
$$
Since $\upsilon_{\infty}(z)=\log8-2\log(1+|z|^2)$,  we   get
$$
\aligned
\int_{\mathbb{R}^2}
\frac{8}{(1+|z|^2)^2}
\upsilon_{\infty}(z)dz
&=8\pi(\log8-2),
\\
\int_{\mathbb{R}^2}
\frac{8}{(1+|z|^2)^2}\big(\upsilon_{\infty}\big)^2(z)
dz
&=8\pi(\log^28-4\log8+8),
\\
\int_{\mathbb{R}^2}
\frac{8}{(1+|z|^2)^2}
\big(\upsilon_{\infty}\big)^3(z)
dz
&=8\pi(\log^38-6\log^28+24\log8+48).
\endaligned
$$
By changing  variables $y=\mu_i z$ and using the equality
$\omega_{\mu_i}(\mu_iz)=\upsilon_{\infty}(z)-2\log\mu_i$,
we can easily deduce (\ref{8.8})-(\ref{8.10}).
\qquad\,$\square$

\vspace {1mm}
\vspace {1mm}
\vspace {1mm}
\vspace {1mm}

\noindent{\bf Lemma A.4.}\,\,\,{\it
\begin{equation}\label{8.7}
\aligned
\int_{\mathbb{R}^2}
e^{\omega_{\mu_i}\left(y-\xi'_i\right)}
B_1
dy
=\int_{\mathbb{R}^2}
e^{\omega_{\mu_i}\left(y-\xi'_i\right)}
\left[\omega^{1}_{\mu_i}+\frac{1}{2}(\omega_{\mu_i})^2\right](y-\xi'_i)
dy
=16\pi\big(
2\log\mu_i
-\log8+2
\big).
\endaligned
\end{equation}
}

\noindent{\it Proof.}
Using the change of  variables $\mu_i z=y-\xi'_i$, by (\ref{8.3}) and (\ref{8.2}) we obtain
$$
\aligned
&\,\int_{\mathbb{R}^2}
e^{\omega_{\mu_i}\left(y-\xi'_i\right)}
\left[\omega^{1}_{\mu_i}+\frac{1}{2}(\omega_{\mu_i})^2\right](y-\xi'_i)
dy
\\
=&\int_{\mathbb{R}^2}
\frac{8}{(1+|z|^2)^2}
\left\{
\left[
\frac12\big(\upsilon_{\infty}\big)^2
-\omega^0_{\infty}
\right](z)
+\Big[
2\log^2\mu_i
-2\log\mu_i
+(\log8+1
)
(1-2\log\mu_i
)
\Big]
\frac{|z|^{2}-1}{|z|^{2}+1}
\right.\\[1mm]
&
\left.\,
+\,4\log\mu_i-1-\log\frac{8}{(1+|z|^2)^2}
\right\}dz.
\endaligned
$$
Note that
$$
\aligned
\int_{\mathbb{R}^2}
\frac{8}{(1+|z|^2)^2}\frac{|z|^{2}-1}{|z|^{2}+1}
dz=0
\qquad
\quad
\textrm{and}
\quad
\qquad
\int_{\mathbb{R}^2}
\frac{8}{(1+|z|^2)^2}
\log\frac{8}{(1+|z|^2)^2}
dz=8\pi(\log8-2).
\endaligned
$$
Then by  (\ref{5.17}),  we have that (\ref{8.7}) holds.
\qquad\qquad\qquad\qquad\qquad\qquad\qquad\qquad
\qquad\qquad
\qquad\qquad\qquad\qquad\qquad\qquad\qquad\qquad$\square$

\vspace {1mm}
\vspace {1mm}
\vspace {1mm}
\vspace {1mm}

\noindent{\bf Lemma  A.5.}\,\,\,{\it
\begin{equation}\label{8.11}
\aligned
\int_{\mathbb{R}^2}
e^{\omega_{\mu_i}}
\omega^{1}_{\mu_i}
=4\pi\left[-4\log^2\mu_i+4(\log8)\log\mu_i
-\log^28\right].
\endaligned
\end{equation}
}

\noindent{\it Proof.}
This is an immediate consequence of (\ref{8.9}) and (\ref{8.7}).
\qquad\qquad\qquad\qquad\qquad\qquad\qquad\qquad\qquad\qquad
\qquad\qquad\qquad\,\,$\square$

\vspace {1mm}
\vspace {1mm}
\vspace {1mm}
\vspace {1mm}

\noindent{\bf Lemma A.6.}\,\,\,{\it
\begin{eqnarray}\label{8.12}
\int_{\mathbb{R}^2}
e^{\omega_{\mu_i}}
\omega^{1}_{\mu_i}\omega_{\mu_i}
=4\pi
\left[
8\log^3\mu_i
+(4-12\log8)\log^2\mu_i
+(6\log^28-4\log8+24)\log\mu_i
-\log^38+\log^28-12\log8-64
\right].
\end{eqnarray}
}

\noindent{\it Proof.}
Testing the equation (\ref{2.9}) for $j=1$ against $\omega_{\mu_i}$ and
the equation $-\Delta\omega_{\mu_i}=e^{\omega_{\mu_i}}$ against $\omega^1_{\mu_i}$, respectively, we find
\begin{equation}\label{8.13}
\aligned
\int_{\mathbb{R}^2}
e^{\omega_{\mu_i}}
\omega^{1}_{\mu_i}
\omega_{\mu_i}
=
\int_{\mathbb{R}^2}
\big(
\omega^1_{\mu_i}\Delta
\omega_{\mu_i}
-
\omega_{\mu_i}\Delta
\omega^1_{\mu_i}
\big)
+
\int_{\mathbb{R}^2}
e^{\omega_{\mu_i}}
\omega^1_{\mu_i}
-\int_{\mathbb{R}^2}
e^{\omega_{\mu_i}}
\left[
(\omega_{\mu_i})^2+\frac{1}{2}(\omega_{\mu_i})^3
\right].
\endaligned
\end{equation}
For any $r>1$ large enough,
by making the change of  variables $y=\mu_i z$ and using  the divergence theorem  we compute
$$
\aligned
\int_{B_{r\mu_i}(0)}
\big(
\omega^1_{\mu_i}\Delta
\omega_{\mu_i}
-
\omega_{\mu_i}\Delta
\omega^1_{\mu_i}
\big)dy
=&
\int_{B_{r}(0)}
\left[
\widetilde{\omega}^1_{\mu_i}\Delta
\upsilon_{\infty}
-
\big(\upsilon_{\infty}-2\log\mu_i\big)\Delta
\widetilde{\omega}^1_{\mu_i}
\right]dz\\
=&\,
2\pi
r\left[
\big(
\upsilon_{\infty}\big)'\,
\widetilde{\omega}^1_{\mu_i}
-
\big(
\widetilde{\omega}^1_{\mu_i}
\big)'\,
\big(\upsilon_{\infty}-2\log\mu_i\big)
\right](r)\\
=&\,
4\pi
\big[
-8\log^2\mu_i
+(8\log8-8)\log\mu_i
+4\log8-2\log^28
\big]+O\left(
\frac{\log r}{r}\right),
\endaligned
$$
where the last equality is due to the expansions
$$
\aligned
\widetilde{\omega}^1_{\mu_i}(r)
=(4\log8-8-8\log\mu_i)\log r
+
O\left(
\frac{\log r}{r^2}\right),
\qquad\qquad
\big(
\widetilde{\omega}^1_{\mu_i}
\big)'(r)
=\frac{1}{r}(4\log8-8-8\log\mu_i)
+O\left(
\frac{\log r}{r^2}\right).
\endaligned
$$
Then
\begin{equation}\label{8.14}
\aligned
\int_{\mathbb{R}^2}
\big(
\omega^1_{\mu_i}\Delta
\omega_{\mu_i}
-
\omega_{\mu_i}\Delta
\omega^1_{\mu_i}
\big)
=4\pi
\big[
-8\log^2\mu_i
+(8\log8-8)\log\mu_i
+4\log8-2\log^28
\big].
\endaligned
\end{equation}
As a consequence, substituting (\ref{8.9}), (\ref{8.10}) , (\ref{8.11}), (\ref{8.14})  and into (\ref{8.13}),
we can deduce (\ref{8.12}).
\qquad\qquad\qquad\qquad\qquad$\square$

\vspace {1mm}
\vspace {1mm}
\vspace {1mm}
\vspace {1mm}

\noindent{\bf Lemma  A.7.}\,\,\,{\it
\begin{eqnarray}
\int_{\mathbb{R}^2}
e^{\omega_{\mu_i}\left(y-\xi'_i\right)}
\big[\,p A_1+(p-1)B_1\big]
dy
=\int_{\mathbb{R}^2}
e^{\omega_{\mu_i}\left(y-\xi'_i\right)}
\left\{
p\,\omega_{\mu_i}
+(p-1)\left[\omega^{1}_{\mu_i}+\frac{1}{2}(\omega_{\mu_i})^2\right]
\right\}(y-\xi'_i)
dy
&&
\nonumber\\[1mm]
=8\pi(p-2)\big(\,
2\log\mu_i
-\log8+2\,
\big),
\qquad\qquad\qquad\qquad\qquad\qquad\,
&&
\label{8.17}
\end{eqnarray}
and
\begin{eqnarray}
\int_{\mathbb{R}^2}
e^{\omega_{\mu_i}\left(y-\xi'_i\right)}
(A_1B_1+B_1)
dy
=\int_{\mathbb{R}^2}
e^{\omega_{\mu_i}\left(y-\xi'_i\right)}
\left\{
\omega_{\mu_i}
\left[\omega^{1}_{\mu_i}+\frac{1}{2}(\omega_{\mu_i})^2\right]
+\left[\omega^{1}_{\mu_i}+\frac{1}{2}(\omega_{\mu_i})^2\right]
\right\}(y-\xi'_i)
dy
&&
\nonumber\\[1.5mm]
=2\pi\big[
-40\log^2\mu_i+(
40\log8-32)\log\mu_i
-10\log^28
+16\log8-16
\big],
\,\,\,\,
&&
\label{8.16}
\end{eqnarray}
and
\begin{eqnarray}
\int_{\mathbb{R}^2}
e^{\omega_{\mu_i}\left(y-\xi'_i\right)}
(A_2+A_1B_1)
dy
=\int_{\mathbb{R}^2}
e^{\omega_{\mu_i}\left(y-\xi'_i\right)}
\left\{
\left[\omega^{1}_{\mu_i}+\frac{p-2}{2(p-1)}(\omega_{\mu_i})^2\right]
+\omega_{\mu_i}
\left[\omega^{1}_{\mu_i}+\frac{1}{2}(\omega_{\mu_i})^2\right]
\right\}(y-\xi'_i)
dy
&&
\nonumber\\[1mm]
=2\pi
\left\{
-\left(\frac{8}{p-1}
+40\right)
\log^2\mu_i
+\left[
\left(
\frac{8}{p-1}+40
\right)\log8
-\frac{16}{p-1}-32
\right]\log\mu_i
\quad\,
\right.
&&
\nonumber\\[1mm]
\left.
-
\left(\frac{2}{p-1}
+10\right)\log^28
+\left(
\frac{8}{p-1}+16
\right)\log8
-\frac{16}{p-1}-16
\right\}.
\qquad\qquad\qquad
\,
&&
\label{8.15}
\end{eqnarray}
}


\noindent{\it Proof.}
These are  an immediate consequence of (\ref{8.8})-(\ref{8.12}).
\qquad\qquad\qquad\qquad\qquad
\qquad\qquad\qquad\qquad\qquad
\qquad\qquad\qquad\quad\,\,\,$\square$

\vspace {1mm}
\vspace {1mm}
\vspace {1mm}
\vspace {1mm}

\noindent{\bf Corollary  A.8.}\,\,\,{\it
For any $0<p\leq2$,
\begin{eqnarray*}
\frac{\,2-p\,}{p}
\int_{\mathbb{R}^2}
e^{\omega_{\mu_i}\left(y-\xi'_i\right)}
B_1
dy
+
\frac{2}{\,p\,}\int_{\mathbb{R}^2}
e^{\omega_{\mu_i}\left(y-\xi'_i\right)}
\big[\,p A_1+(p-1)B_1\big]
dy
\,\equiv\,
0.
\end{eqnarray*}
}

\section{Proof of the identity (\ref{9.10}) of second order}
Let
$$
\aligned
\psi_0(t)=8t\frac{t^2-1}{(t^2+1)^3},
\qquad\qquad
\eta_0(t)=\log(1+t^2),
\qquad\qquad
\zeta_0(t)=
\frac{1}{t^2+1},
\qquad\qquad
\theta_0(t)=
\int_{t^2}^{+\infty}\frac{\log(s+1)}{s(s+1)}ds.
\endaligned
$$
Some straightforward but very tedious computations give
$$
\aligned
\int_{0}^{+\infty}\psi_0dt
&=0,
\qquad\quad
\int_{0}^{+\infty}\psi_0\eta_0dt
=2,
\qquad\quad
\int_{0}^{+\infty}\psi_0\eta^2_0dt
=6,
\qquad\quad
\int_{0}^{+\infty}\psi_0\eta^3_0dt
=21,
\qquad\quad
\int_{0}^{+\infty}\psi_0\eta^4_0dt
=90,
\endaligned
$$
and
$$
\aligned
\int_{0}^{+\infty}\psi_0\zeta_0dt
&=-\frac23,
\qquad\quad\qquad\qquad
\int_{0}^{+\infty}\psi_0\zeta_0\eta_0dt
=\frac19,
\endaligned
$$
$$
\aligned
\int_{0}^{+\infty}\psi_0\zeta_0\eta^2_0dt
=\frac{11}{27},
\qquad\qquad\quad\qquad
\int_{0}^{+\infty}\psi_0\zeta_0\eta^3_0dt
=\frac{49}{54},
\qquad\qquad\quad\qquad
\int_{0}^{+\infty}\psi_0\zeta_0\eta^4_0dt
=\frac{179}{81},
\endaligned
$$
and
$$
\aligned
\int_{0}^{+\infty}\psi_0\zeta^2_0dt
&=-\frac23,
\qquad\quad\qquad\qquad
\int_{0}^{+\infty}\psi_0\zeta^2_0\eta_0dt
=-\frac{1}{18},
\endaligned
$$
$$
\aligned
\int_{0}^{+\infty}\psi_0\zeta^2_0\eta^2_0dt
=\frac{5}{108},
\qquad\qquad\quad\qquad
\int_{0}^{+\infty}\psi_0\zeta^2_0\eta^3_0dt
=\frac{47}{432},
\qquad\qquad\quad\qquad
\int_{0}^{+\infty}\psi_0\zeta^2_0\eta^4_0dt
=\frac{269}{1296}.
\endaligned
$$
In particular, integrating by parts, we can compute
$$
\aligned
\int_{0}^{+\infty}\psi_0\theta_0dt
=-1,
\qquad\quad\qquad\qquad
\int_{0}^{+\infty}\psi_0\zeta_0\theta_0dt
=-\frac{61}{54},
\qquad\quad\qquad\qquad
\int_{0}^{+\infty}\psi_0\eta_0\theta_0dt
=\frac{1}{2},
\endaligned
$$
and
$$
\aligned
\int_{0}^{+\infty}\psi_0\zeta_0^2\theta_0dt
=-\frac{223}{216},
\,\,\qquad\quad\quad\,\,
\int_{0}^{+\infty}\psi_0\eta^2_0\theta_0dt
=11-8\zeta(3),
\,\,\qquad\quad\quad\,\,
\int_{0}^{+\infty}\psi_0\eta_0\zeta_0\theta_0dt
=\frac{4}{3}\zeta(3)-\frac{179}{108},
\endaligned
$$
and
$$
\aligned
\int_{0}^{+\infty}\psi_0\eta^2_0\zeta_0\theta_0dt
=4\zeta(4)-\frac{4}{9}\zeta(3)-\frac{589}{162},
\,\,\qquad\,\,\,\,
\,\,\,\,\qquad\,\,
\int_{0}^{+\infty}\psi_0\eta^2_0\zeta^2_0\theta_0dt
=4\zeta(4)+\frac{2}{9}\zeta(3)-\frac{11893}{2592},
\endaligned
$$
where
$\zeta(3)$ and $\zeta(4)$   are  two positive irrational
numbers defined in Ap\'{e}ry's constants
 (see \cite{CDT}) and
$\zeta(\cdot)$ denotes the famous
Euler-Riemann zeta function
$$
\aligned
\zeta(z)=\sum_{n=1}^{+\infty}\frac{1}{n^z}
=\frac{1}{\Gamma(z)}\int_{0}^{+\infty}\frac{t^{z-1}}{e^{t}-1}dt
=\frac{1}{\Gamma(z)}\int_{0}^{+\infty}\frac{\log^{z-1}(s+1)}{s(s+1)}ds,
\qquad
\textrm{for any}
\,\,\,
z\in\mathbb{C}
\,\,\,
\textrm{and}
\,\,\,
Re(z)>1.
\endaligned
$$
Furthermore, integrating by  parts twice,  we can compute
$$
\aligned
\int_{0}^{+\infty}\psi_0\theta^2_0dt
=8\zeta(3)-\frac{23}{2},
\qquad\quad\qquad
\int_{0}^{+\infty}\psi_0\zeta_0\theta^2_0dt
=\frac{68}{9}\zeta(3)-\frac{3517}{324},
\qquad\quad\qquad
\int_{0}^{+\infty}\psi_0\zeta^2_0\theta^2_0dt
=\frac{62}{9}\zeta(3)-\frac{51127}{5184}.
\endaligned
$$

\vspace{1mm}
\vspace{1mm}

\noindent{\bf Lemma B.1.}\,\,\,{\it
\begin{eqnarray}
\int_{0}^{+\infty}\psi_0(\widetilde{\omega}_{\mu_i})^2dt
=-8(\log8-2\log\mu_i)+24,
\qquad\qquad\qquad\qquad\qquad
&&\label{9.1}\\
\int_{0}^{+\infty}\psi_0(\widetilde{\omega}_{\mu_i})^3dt
=-12(\log8-2\log\mu_i)^2+72(\log8-2\log\mu_i)-168,
\qquad\qquad\qquad
&&\label{9.2}\\
\int_{0}^{+\infty}\psi_0(\widetilde{\omega}_{\mu_i})^4dt
=-16(\log8-2\log\mu_i)^3+144(\log8-2\log\mu_i)^2-672(\log8-2\log\mu_i)+1440.
&&\label{9.3}
\end{eqnarray}
}

\noindent{\it Proof.}
These are an immediate consequence of the above integral computations because of
$\widetilde{\omega}_{\mu_i}=\log8-2\log\mu_i-2\eta_0$.
\qquad$\square$

\vspace{1mm}
\vspace{1mm}
\vspace{1mm}

\noindent{\bf Lemma B.2.}\,\,\,{\it
\begin{eqnarray}
\int_{0}^{+\infty}\psi_0\widetilde{\omega}^1_{\mu_i}dt
=\frac{8}{3}\log^2\mu_i
-\Big(\frac{8}{3}\log8+\frac{40}{3}\Big)\log\mu_i
+\frac{2}{3}\log^28+\frac{20}{3}\log8
-20,
\qquad\qquad\qquad
&&\label{9.4}\\[1mm]
\int_{0}^{+\infty}\psi_0\widetilde{\omega}^1_{\mu_i}\widetilde{\omega}_{\mu_i}dt
=-\frac{16}{3}\log^3\mu_i
+\Big(
8\log8+
\frac{248}{9}
\Big)\log^2\mu_i
+
\Big(
\frac{776}{9}
-\frac{248}{9}\log8
-4\log^28
\Big)\log\mu_i
&&\nonumber\\
+\frac{2}{3}\log^38
+\frac{62}{9}\log^28
-\frac{388}{9}\log8
-\frac{64}{3}\zeta(3)+
84,
\qquad\qquad\qquad\qquad
\qquad\ \ \,
&&\label{9.5}\\[1mm]
\int_{0}^{+\infty}\psi_0\widetilde{\omega}^1_{\mu_i}(\widetilde{\omega}_{\mu_i})^2dt
=\frac{32}{3}\log^4\mu_i
-\Big(\frac{64}{3}\log8+\frac{512}{9}\Big)\log^3\mu_i
+\Big(16\log^28+\frac{256}{3}\log8-\frac{7328}{27}\Big)\log^2\mu_i
\qquad
&&\nonumber\\
+\Big(
-\frac{16}{3}\log^38-\frac{128}{3}\log^28+\frac{7328}{27}\log8-\frac{17792}{27}+\frac{256}{3}\zeta(3)
\Big)\log\mu_i
+\frac{2}{3}\log^48
\quad\,\,\,\,
&&\nonumber\\
+\frac{64}{9}\log^38
-\frac{1832}{27}\log^28
+\Big(
\frac{8896}{27}-\frac{128}{3}\zeta(3)\Big)\log8
-\frac{1952}{3}+\frac{1024}{9}\zeta(3)+128\zeta(4).
&&\label{9.6}
\end{eqnarray}
}

\noindent{\it Proof.}
From the expressions of $\omega^0_{\infty}$ and  $\omega^1_{\infty}$ in (\ref{8.4}) and
(\ref{8.5}), respectively, we get
$$
\aligned
\omega^0_{\infty}=(\log^28+2\log8-10)\zeta_0+4\zeta_0\eta_0^2
+(6-2\log8)\eta_0+4\theta_0
-8\zeta_0\theta_0,
\endaligned
$$
$$
\aligned
\omega^1_{\infty}=2(1+\log8)\zeta_0
-2\eta_0.
\endaligned
$$
Hence by (\ref{8.6}),
$$
\aligned
\widetilde{\omega}^1_{\mu_i}=
&
\,\,
(10-2\log8-\log^28)\zeta_0-4\zeta_0\eta_0^2
+2(\log8-3)\eta_0-4\theta_0
+8\zeta_0\theta_0\\
&
+2(2\log\mu_i-1)
\Big[
(1+\log8)\zeta_0
-\eta_0
\Big]
+4(\log\mu_i-\log^2\mu_i)
\zeta_0,
\endaligned
$$
and
$$
\aligned
\widetilde{\omega}^1_{\mu_i}\widetilde{\omega}_{\mu_i}=&
\,\,
(\log8-2\log\mu_i)\widetilde{\omega}^1_{\mu_i}
-2(10-2\log8-\log^28)\zeta_0\eta_0
+8\zeta_0\eta_0^3
-4(\log8-3)\eta^2_0+8\eta_0\theta_0
\\
&
-16\zeta_0\eta_0\theta_0-4(2\log\mu_i-1)
\Big[
(1+\log8)\zeta_0\eta_0
-\eta^2_0
\Big]
-8(\log\mu_i-\log^2\mu_i)
\zeta_0\eta_0,
\endaligned
$$
and
$$
\aligned
\widetilde{\omega}^1_{\mu_i}(\widetilde{\omega}_{\mu_i})^2=
&
\,\,
(\log8-2\log\mu_i)^2\widetilde{\omega}^1_{\mu_i}
+4(10-2\log8-\log^28)\zeta_0\eta_0^2-16\zeta_0\eta_0^4
+8(\log8-3)\eta^3_0-16\eta^2_0\theta_0
\\
&
+32\zeta_0\eta^2_0\theta_0+8(2\log\mu_i-1)
\Big[
(1+\log8)\zeta_0\eta^2_0
-\eta^3_0
\Big]
+16(\log\mu_i-\log^2\mu_i)
\zeta_0\eta^2_0\\
&
+4(2\log\mu_i-\log8)
\Big\{
(10-2\log8-\log^28)\zeta_0\eta_0
-4\zeta_0\eta_0^3
+2(\log8-3)\eta^2_0
-4\eta_0\theta_0
\\
&
+8\zeta_0\eta_0\theta_0+2(2\log\mu_i-1)
\Big[
(1+\log8)\zeta_0\eta_0
-\eta^2_0
\Big]
+4(\log\mu_i-\log^2\mu_i)
\zeta_0\eta_0
\Big\}.
\endaligned
$$
These equalities combined with the previous integral computations can derive that
(\ref{9.4})-(\ref{9.6}) hold.
\qquad\qquad\qquad\qquad\qquad$\square$

\vspace{1mm}
\vspace{1mm}
\vspace{1mm}
\vspace{1mm}

\noindent{\bf Lemma B.3.}\,\,\,{\it
\begin{eqnarray}
\int_{0}^{+\infty}\psi_0(\widetilde{\omega}^1_{\mu_i})^2dt
=-\frac{32}{3}\log^4\mu_i
+\Big(
\frac{64}{3}\log8
+\frac{416}{9}
\Big)
\log^3\mu_i
+\Big(-16\log^28
-\frac{208}{3}\log8
+\frac{3344}{27}
\Big)\log^2\mu_i
\nonumber\\[1mm]
+\Big(
\frac{16}{3}\log^38
+\frac{104}{3}\log^28
-\frac{3344}{27}\log8
-\frac{256}{3}\zeta(3)
+\frac{4880}{27}
\Big)\log\mu_i
-\frac{2}{3}\log^48
\qquad\quad\,
\nonumber\\[1mm]
-\frac{52}{9}\log^38
+\frac{836}{27}\log^28
+\Big(\frac{128}{3}\zeta(3)-
\frac{2440}{27}
\Big)\log8
-128\zeta(4)-\frac{256}{9}\zeta(3)
+\frac{584}{3}.
\quad\,\,
&&\label{9.7}
\end{eqnarray}
}

\noindent{\it Proof.}
Notice that
$$
\aligned
(\widetilde{\omega}^1_{\mu_i})^2=
&
\,\,
(10-2\log8-\log^28)^2\zeta^2_0+16\zeta^2_0\eta_0^4
+4(\log8-3)^2\eta_0^2+16\theta^2_0
+64\zeta^2_0\theta^2_0
+16(\log\mu_i-\log^2\mu_i)^2
\zeta^2_0\\
&
+4(2\log\mu_i-1)^2
\Big[
(1+\log8)^2\zeta_0^2
+\eta^2_0
-2(1+\log8)\zeta_0\eta_0
\Big]
-8(10-2\log8-\log^28)\zeta^2_0\eta_0^2
\\
&
+4(10-2\log8-\log^28)(\log8-3)\zeta_0\eta_0
-8(10-2\log8-\log^28)\zeta_0\theta_0
+16(10-2\log8-\log^28)\zeta^2_0\theta_0\\
&
+4(10-2\log8-\log^28)(2\log\mu_i-1)
\Big[
(1+\log8)\zeta^2_0
-\zeta_0\eta_0
\Big]\\
&
+8(10-2\log8-\log^28)(\log\mu_i-\log^2\mu_i)
\zeta^2_0
-16(\log8-3)\zeta_0\eta^3_0+32\zeta_0\eta^2_0\theta_0
-64\zeta^2_0\eta^2_0\theta_0\\
&
-16(2\log\mu_i-1)
\Big[
(1+\log8)\zeta^2_0\eta_0^2
-\zeta_0\eta_0^3
\Big]
-32(\log\mu_i-\log^2\mu_i)
\zeta_0^2\eta_0^2\\
&
-16(\log8-3)\eta_0\theta_0
+32(\log8-3)\zeta_0\eta_0\theta_0
+8(\log8-3)(2\log\mu_i-1)
\Big[
(1+\log8)\zeta_0\eta_0
-\eta^2_0
\Big]\\
&
+16(\log8-3)(\log\mu_i-\log^2\mu_i)
\zeta_0\eta_0
-64\zeta_0\theta^2_0
-16(2\log\mu_i-1)
\Big[
(1+\log8)\zeta_0\theta_0
-\eta_0\theta_0
\Big]\\
&
-32(\log\mu_i-\log^2\mu_i)
\zeta_0\theta_0
+32(2\log\mu_i-1)
\Big[
(1+\log8)\zeta^2_0\theta_0
-\zeta_0\eta_0\theta_0
\Big]\\
&
+64(\log\mu_i-\log^2\mu_i)
\zeta^2_0\theta_0
+16(2\log\mu_i-1)
(\log\mu_i-\log^2\mu_i)
\Big[
(1+\log8)\zeta^2_0
-\zeta_0\eta_0
\Big],
\endaligned
$$
This equality combined with the previous integral computations  can derive that
(\ref{9.7}) holds.
\qquad\qquad\qquad\qquad\qquad\qquad\qquad$\square$

\vspace{1mm}
\vspace{1mm}
\vspace{1mm}
\vspace{1mm}

\noindent{\bf Theorem B.4.}\,\,\,{\it
\begin{eqnarray}\label{9.8}
D^2_{\mu_i}=\int_{0}^{+\infty}\psi_0\widetilde{f}^2_{\mu_i}dt
=-\Big(
\frac{8}{p-1}
+24
\Big)\log^2\mu_i
+
\Big[
\Big(
\frac{8}{p-1}
+24
\Big)\log8
-
\frac{16}{p-1}
\Big]
\log\mu_i
-\Big(
\frac{2}{p-1}+6\Big)\log^28
&&\nonumber\\[1mm]
+\,
\frac{8}{p-1}
\log8
-\frac{16}{p-1}.
\qquad\qquad\qquad\qquad\qquad\qquad\qquad
\qquad\qquad\qquad\qquad\qquad\quad\,
&&
\end{eqnarray}
}

\noindent{\it Proof.}
From (\ref{2.11}) and (\ref{8.3}) we get
$$
\aligned
\widetilde{f}^2_{\mu_i}=-
\frac{1}{2}(\widetilde{\omega}^1_{\mu_i})^2
-
\frac{1}{2}\widetilde{\omega}^1_{\mu_i}(\widetilde{\omega}_{\mu_i})^2
-
2\widetilde{\omega}^1_{\mu_i}\widetilde{\omega}_{\mu_i}
-
\widetilde{\omega}^1_{\mu_i}
-
\frac{1}{8}(\widetilde{\omega}_{\mu_i})^4
-
\frac{4p-5}{6(p-1)}(\widetilde{\omega}_{\mu_i})^3
-
\frac{p-2}{2(p-1)}(\widetilde{\omega}_{\mu_i})^2.
\endaligned
$$
Applying (\ref{9.1})-(\ref{9.7}) and  the second definition  of $D^2_{\mu_i}$ in (\ref{2.14}),
we can derive that (\ref{9.8}) holds.
\qquad\qquad\qquad\qquad\qquad\qquad\quad$\square$

\vspace {1mm}
\vspace {1mm}
\vspace {1mm}
\vspace {1mm}

\noindent{\bf Lemma  B.5.}\,\,\,{\it
\begin{eqnarray}
\int_{\mathbb{R}^2}
e^{\omega_{\mu_i}\left(y-\xi'_i\right)}
\left(B_2+\frac12(B_1)^2\right)
dy
=\int_{\mathbb{R}^2}
e^{\omega_{\mu_i}\left(y-\xi'_i\right)}
\left\{
\left[
\omega^{2}_{\mu_i}+\omega_{\mu_i}\omega^{1}_{\mu_i}+\frac{p-2}{6(p-1)}(\omega_{\mu_i})^3\right]
+
\frac12\left[\omega^{1}_{\mu_i}+\frac{1}{2}(\omega_{\mu_i})^2\right]^2
\right\}(y-\xi'_i)
dy
&&\nonumber\\[1mm]
=2\pi\left\{
\Big(
\frac{16}{p-1}
+64
\Big)\log^2\mu_i
+\Big[
\frac{32}{p-1}
+32
-
\Big(
\frac{16}{p-1}\
+
64
\Big)\log8
\Big]
\log\mu_i
\qquad\qquad\qquad\,\,
\right.
&&\nonumber\\[1mm]
\left.
+\,\Big(
\frac{4}{p-1}
+16
\Big)\log^28
-
\Big(
\frac{16}{p-1}
+
16
\Big)\log8
+\frac{32}{p-1}
+16
\right\}.
\ \,\qquad\qquad\qquad\qquad\quad\quad
&&\label{9.9}
\end{eqnarray}
}

\noindent{\it Proof.}
From (\ref{2.9}) and (\ref{2.11})  we find
$$
\aligned
\Delta\left[\omega^2_{\mu_i}(y-\xi'_i)\right]
=e^{\omega_{\mu_i}\left(y-\xi'_i\right)}\big(f^2_{\mu_i}
-\omega^2_{\mu_i}
\big)(y-\xi'_i)
=
-e^{\omega_{\mu_i}\left(y-\xi'_i\right)}
\left[
A_2+A_1B_1+
B_2+\frac12(B_1)^2
\right]
\qquad
\textrm{in}\,\,\ \,\,\mathbb{R}^2,
\endaligned
$$
Using the the first definition  of $D^2_{\mu_i}$ in (\ref{2.14}), we obtain
$$
\aligned
\int_{\mathbb{R}^2}
e^{\omega_{\mu_i}\left(y-\xi'_i\right)}
\Big(B_2+\frac12(B_1)^2\Big)
dy
=-2\pi D^2_{\mu_i}
-\int_{\mathbb{R}^2}e^{\omega_{\mu_i}\left(y-\xi'_i\right)}
\left(
A_2+A_1B_1
\right)dy.
\endaligned
$$
Hence by
 (\ref{8.15}) and (\ref{9.8}),
we can derive that (\ref{9.9}) holds.
\qquad\qquad\qquad\qquad\qquad\qquad\qquad\qquad\qquad
\qquad\qquad\qquad\qquad\qquad$\square$

\vspace {1mm}
\vspace {1mm}
\vspace {1mm}
\vspace {1mm}

\noindent{\bf Corollary  B.6.}\,\,\,{\it
For any $0<p\leq2$,
\begin{eqnarray*}
\frac{p-1}{8\pi}
\left[\int_{\mathbb{R}^2}
e^{\omega_{\mu_i}\left(y-\xi'_i\right)}
\Big(B_2+\frac12(B_1)^2\Big)
dy
+2
\int_{\mathbb{R}^2}
e^{\omega_{\mu_i}\left(y-\xi'_i\right)}
(A_1B_1+B_1)
dy
\right]
+\frac{p-2}{\,(16\pi)^2\,}
\left(
\int_{\mathbb{R}^2}
e^{\omega_{\mu_i}\left(y-\xi'_i\right)}
B_1
dy
\right)^2\,\equiv\,4.
\end{eqnarray*}
}

\noindent{\it Proof.}
This is an immediate consequence of (\ref{8.7}),  (\ref{8.16})  and (\ref{9.9}).
\qquad\qquad\qquad\qquad\qquad\qquad\qquad\qquad\qquad\qquad\qquad\quad\,\,$\square$

\end{document}